\RequirePackage{rotating}
\documentclass[preprint,authoryear,3p,11pt]{elsarticle}
\usepackage{amssymb, amsmath, textcomp, gensymb, makeidx, epsfig, graphicx, graphics, subfig, url, multirow, multicol, booktabs, array, bm, rotating, setspace, tabularx, paralist, afterpage, hyperref, enumitem, siunitx, longtable, lscape}
\usepackage{amscd,verbatim,color,algpseudocode}
\usepackage{algorithm}
\usepackage[utf8]{inputenc}
\usepackage[bottom]{footmisc}
\usepackage[strings]{underscore}
\usepackage[section]{placeins}
\usepackage[table]{xcolor}
\usepackage{hhline}
\usepackage{amsthm} 
\usepackage{thmtools, thm-restate}
\hypersetup{colorlinks=true, linkcolor=blue, anchorcolor =red,citecolor =blue,filecolor = red,urlcolor = red,
            pdfauthor=author}

\graphicspath{{figures/}}
\DeclareMathAlphabet{\mathcal}{OMS}{cmsy}{m}{n}

\newtheorem{definition}{Definition}

\newtheorem{proposition}{Proposition}

\newtheorem{lemma}{Lemma}
\newtheorem{conjecture}{Conjecture}

\makeatletter
\renewcommand*\env@matrix[1][c]{\hskip -\arraycolsep
  \let\@ifnextchar\new@ifnextchar
  \array{*\c@MaxMatrixCols #1}}
\makeatother

\newcommand{\costf}{Q}

\newcommand{\reduced}{\mathcal{Z}}
\newcommand{\zeromean}{\mathcal{N}}

\newcommand{\Nneq}{M}

\newcommand{\ones}{J}

\newcommand{\Tr}{\mathsf{T}}

\newcommand{\disunif}[2]{\mathcal{U} \{#1, #2\}}
\newcommand{\conunif}[2]{\mathcal{U} [#1, #2]}

\DeclareMathOperator{\trace}{tr}

\journal{European Journal of Operational Research}

\begin{document}
	
	\begin{frontmatter}
	
	\title{Instance Space Analysis for the Quadratic Assignment Problem \tnoteref{t1}}
		\tnotetext[t1]{\textsuperscript{\textcopyright}\ 2025. This manuscript version is made available under the CC-BY-NC-ND 4.0 license \url{http://creativecommons.org/licenses/by-nc-nd/4.0/}}
		\author[unimelb]{Jeffrey Christiansen\corref{cor}}
		\ead{jeffreyc1@unimelb.edu.au}
		\author[unimelb]{Kate Smith-Miles}
		\ead{smith-miles@unimelb.edu.au}

 	\cortext[cor]{Corresponding author}
 	\address[unimelb]{ARC Training Centre in Optimisation Technologies, Integrated Methodologies and Applications (OPTIMA), School of Mathematics and Statistics, The University of Melbourne, Parkville, Victoria, 3010, Australia.}

        \begin{abstract}
            For any optimisation problem where diverse algorithmic approaches are available, the task of predicting algorithm performance and selecting the algorithm most likely to perform well on a given instance holds great practical interest. However, if our test instances do not adequately represent the potential problem space, we may be misled about the strengths and weaknesses of an algorithm. In this paper we consider algorithm prediction and selection for the Quadratic Assignment Problem (QAP). We identify and eliminate superficial differences between QAP instances which do not affect problem difficulty and propose several new features for quantifying the characteristics of a particular instance. We then apply Instance Space Analysis to compare the performance of evolutionary and ant colony-based algorithms. Using this analysis, we identify limitations of the existing instance libraries and generators which obscure the true performance profile of these algorithms. Finally, we propose new instance classes which expand the instance space and support new insights into the properties of the problem and algorithms.
		\end{abstract}

	
	\end{frontmatter}
	
	
	\section{Introduction} \label{sec:introduction}
        
        The Quadratic Assignment Problem (QAP) is typically understood as a facility location problem, in which facilities are assigned to locations with the goal of minimising the total cost. The original formulation of the QAP, proposed by \citet{Koopmans1957Original}, defines the cost of a solution as the sum of two terms. The first term establishes the cost of individual assignments, using a problem data matrix which defines the cost associated with each possible combination of one facility and one location. The QAP formulation which we employ in this work simplifies the original formulation by eliminating this term and its associated problem data, allowing us to focus on the more complicated second term. The second term defines the cost associated with pairs of assignments using a flow matrix $F = [f_{ij}]$, which defines the relationships between the facilities, and a distance matrix $D = [d_{rs}]$, which defines the relationships between the locations. $f_{ij}$ represents the required commodity flow from facility $i$ to facility $j$, and $d_{rs}$ represents the cost of transportation for one unit flow from location $r$ to location $s$. Then the quantity $f_{ij} d_{rs}$ represents the cost associated with assigning facility $i$ to location $r$ and assigning facility $j$ to location $s$. However, we may rewrite this quantity as $d_{rs} f_{ij}$ and instead interpret it as the cost of assigning location $r$ to facility $i$ and location $s$ to facility $j$. 
        Therefore the two matrices $F$ and $D$ play a symmetrical role in the definition of the cost function. The cost of assigning a facility to a location, independent of other assignments, may still appear in this formulation as terms of the form $f_{ii}d_{rr}$.
        For a summary of previous work on the QAP we refer the reader to \citet{Loiola2007657}. 
        
        The QAP is an NP-hard problem, and it is very difficult to find and verify optimal solutions for some QAP instances of even moderate size. The well-known QAPLIB instance library \citep{Burkard1997QAPLIB} contains several instances of size $n \leq 50$ which have not been solved to optimality. However, metaheuristics have been successful at producing high-quality solutions even for much larger QAP instances. Several local search methods have been proposed, such as Iterated Local Search \citep{Stutzle2006ILS}, Breakout Local Search \citep{BenlicBLS}, and Tabu Search methods \citep{Drezner2005Tabu,James2009Tabu}. Other algorithmic approaches include ant systems \citep{Talbi2001ant,Demirel2006ant,Lopez2018MMAS}, evolutionary algorithms \citep{Mise2004ga,BenlicBMA}, 
        and greedy randomised adaptive search \citep{Oliveira2004grasp}. These other approaches may also use local searches as a subalgorithm to improve their solutions \citep{BenlicBMA,Lopez2018MMAS}. The performance of these heuristics depends on the characteristic of the instances being solved, with no one algorithm being clearly superior in general \citep{Stutzle2006ILS,Dantas2020landscape}. Therefore, a method for quickly and reliably predicting the performance of algorithms on a given instance, then choosing the most suitable algorithm, would result in higher-quality solutions and/or shorter computational time.

        Early efforts to predict QAP instance difficulty focused on the role of the flow matrix. \citet{Vollmann1966Dom}
        proposed that a problem with few relatively large flows would be simpler and hence not require sophisticated techniques to obtain a good result. They referred to such a problem as having high flow dominance, but did not provide an explicit quantitative definition of this term. Subsequent work by several authors produced a variety of candidate formal definitions, but none of these resulted in a satisfactory predictor of problem difficulty on their own \citep{Herroelen1985}. \citet{Mojena1976} applied flow dominance alongside other features measuring the flow distribution and number of clusters of facilities to predict the performance of a branch and bound algorithm on instances with a fixed distance matrix.

        Later work has also taken into account the role of the distance matrix in determining problem difficulty. In view of the similar role which the distance and flow matrices take in the definition of a QAP instance, the dominance of the distance matrix has been considered \citep{Stutzle2006ILS}. Other features measured for both the distance and flow matrices include sparsity and asymmetry \citep{Dantas2020landscape}. 
        To measure properties of the problem as a whole rather than just one matrix, a simple approach is to measure a feature for both matrices then combine or compare them in some way; for example \citet{Stutzle2006ILS} considered the maximum sparsity of the two matrices. 
        More complicated features measured using both matrices are also possible. \citet{Angel2001landscape} proposed an method for deterministically calculating the autocorrelation coefficient and hence the ruggedness of the fitness landscape of a symmetric QAP instance, and later applied it to predicting instance difficulty \citep{Angel2002hardness}. By applying elementary landscape decomposition, \citet{Chicano2012698} generalised this result to asymmetric QAP instances.   

        In contrast to these deterministic measurements, an alternative approach to quantifying the properties of a QAP instance employs random sampling to survey the fitness landscape, often combined with some sort of local exploration or search procedure to explore the landscape near these sampled points. 
        \citet{Pitzer2013FLA} and \citet{Werth20232108} applied random, up-down and neutral walks to explore the fitness landscape from different perspectives. The Iterated Local Search algorithm proposed by \citet{Stutzle2006ILS} was immediately used to measure the fitness-distance correlation of QAP instances in the same paper. Later work by \citet{Thomson2023} applied ILS to investigate local optima networks. Notably, the studies applying ILS to find local minima focused their analysis on subsets of QAPLIB mostly containing instances of modest size. When analysing a wider range of QAPLIB instances, including instances of size $n=343$, \citet{Dantas2020landscape} noted that the computational time required to conduct the landscape analysis was a point of concern. Accordingly, they applied the simpler Best Improvement Local Search to find local minima, and explicitly examined the effect of modifying the sample size on feature reliability. 

        Algorithm selection for the QAP has been previously studied using a variety of learning methods. \citet{Smith-Miles2008} contrasted supervised and unsupervised learning approaches using neural networks, and also used self-organising maps to project the feature data into a $2$D space to support further insights into the problem and the qualities of the algorithms. This analysis employed a dataset previously compiled by \citet{Stutzle2006ILS}, consisting of 28 QAPLIB instances (with maximum size $n=100$), 4 deterministic features, 5 fitness landscape features, and performance data from 3 algorithms. \citet{Pitzer2013FLA} utilised a dataset of 117 QAPLIB instances and a set of features focusing on fitness landscape analysis using random sampling, initially identifying similar instances using k-means clustering and projecting to a two-dimensional space using Neighbourhood Component Analysis. They then compared the performance of one-variable rule learner, linear regression, and support vector machine (SVM) based approaches to an algorithm selection problem between two algorithms, and found that of these options a variant of the SVM approach using sequential minimal optimisation achieved the best results. \citet{Beham20171471} applied a k-nearest neighbour-based approach, constructing a selector to choose between 7 distinct algorithms based on 17 fitness landscape features measured using a variety of local exploration strategies. This analysis was restricted to 47 instances but included instances defined by \cite{DreznerTaillard2005} which are not present in QAPLIB, including instances of size $n=343$. \citet{Dantas2020landscape} used a Random Forest model to select from 3 algorithms based on 7 deterministic features and 7 fitness landscape features, with an emphasis on the tradeoff between the time spent evaluating the features and the efficiency gained by selecting the most appropriate algorithm; this work included 152 instances from QAPLIB and \citet{DreznerTaillard2005}. When applying any of these machine learning techniques to algorithm selection, a large collection of interesting and diverse test instances is desirable. While QAPLIB does contain 136 instances from 15 different sources, many of those instances are too small to yield useful performance data when comparing metaheuristics running on modern hardware.

        In this paper we apply Instance Space Analysis (ISA) \citep{Smith-Miles2023MATILDA}
        to build greater understanding of the QAP, and in particular to augment the quantity and quality of QAP instances available for future studies. Based on the Algorithm Selection Problem framework proposed by \citet{Rice76}, ISA is a methodology for investigating the relationships between benchmark instance sets, their measurable features, and the performance of algorithms which (approximately or exactly) solve them. Given feature and algorithm performance data for instances in a training set, ISA selects a subset of the features which are most correlated with algorithm performance, then chooses a projection of these feature values into a $2$D instance space which produces observable trends in the algorithm performance data. This projection may be used directly to construct an algorithm selector by identifying regions in the $2$D instance space where each algorithm is observed to typically perform well or poorly on the training set of instances, then predicting the performance of the algorithms on a new instance based on the region in which it falls. Visual inspection of the $2$D space may provide insights into the strengths and weaknesses of algorithms which are not immediately apparent when considering the high-dimensional feature data. 
        
        The 2$D$ space may also be applied to identify combinations of characteristics which are over- or under-represented in the current training set. When these flaws are resolved, an algorithm selector based on the improved training set can be expected to produce a more reliable result. For example, the results of ISA have been used to augment existing instance sets for black-box optimisation \citep{Munoz19a}, the 0-1 Knapsack Problem \citep{Smith-Miles2021KP}, curriculum-based course timetabling \citep{DeCoster2022TT}, and multiobjective problems \citep{Yap2022mobj}. On the other hand, an enhanced ISA framework was applied to the maximum flow problem by \citet{Alipour2023maxflow} to remove instances in over-represented regions of the instance space, reducing bias in the training set. In each of these studies the improved training set resulting from these changes was then used to produce greater insights into algorithm performance. Our goal in this paper is to apply ISA to not only investigate the strengths and weaknesses of metaheuristics for solving the QAP, but also to uncover the limitations of QAPLIB and other existing QAP benchmark instances, and address these limitations with new instances.

        The remainder of the paper is structured as follows. In Section~\ref{sec:isa} we present the components of Instance Space Analysis for the QAP. Section~\ref{sec:problem} defines the QAP formulation used in this work, and Section~\ref{sec:instances} lists the instances present in our initial instance subset. Sections~\ref{sec:algorithms} and~\ref{sec:performance} detail the algorithms being compared and the metric used to assess their performance. Section~\ref{sec:all-theory} contains theoretical results about the QAP which are used to justify the feature definitions in Section~\ref{sec:features}. With the initial instance subset, performance data, and feature space in place, we proceed to apply the ISA tools and interpret the results in Section~\ref{sec:initial}. Using the insights from this analysis, we propose additional instance classes to expand the instance subset, and re-apply ISA to this augmented set in Section~\ref{sec:more-instances}. Section~\ref{sec:conclusion} contains our concluding remarks.
 
	\section{Instance Space Analysis Components} \label{sec:isa}
    In this section we define the instance subset, algorithm data, and feature data to which we apply ISA. For a detailed technical description of ISA we refer the reader to \citep{Smith-Miles2023MATILDA}. The instance data files and a link to the code for instance generation, feature measurement, etc. may be found at \url{https://matilda.unimelb.edu.au/matilda/problems/opt/qap}. 

	\subsection{Problem Space} \label{sec:problem}
	
	We consider the following definition of the QAP. Given the index set $N = \left\{1,2,\dots,n\right\}$ and two $n \times n$ real-valued problem data matrices $A = \left[a_{ij}\right]$ and $B = \left[b_{ij}\right]$, find a permutation $\varphi : N \to N$ such that the cost function 
	\begin{equation} \label{eq:QAP-cost-function}
	    \costf_{A,B}(\varphi):=\sum_{i \in N }\sum_{j \in N} a_{ij} b_{\varphi(i)\varphi(j)}
	\end{equation}
	is minimised. We denote the set of possible permutations of $N$ as $S_n$. Using this notation the QAP instance with data matrices $A$ and $B$, hereafter denoted $QAP(A,B)$, can be written in the form
	\begin{equation} \label{eq:QAP-formulation}
	    \min_{\varphi \in S_n} \costf_{A,B}(\varphi).
	\end{equation}
    This formulation of the QAP is agnostic with respect to application; either of $A$ or $B$ may represent distances between locations, flows between facilities, or any other suitable quantity.
	
	The algorithm implementations and library test instances employed in this work encompass only the subset of QAP instances in which all elements of $A$ and $B$ are non-negative integers. Our instance space analysis is accordingly restricted to this subset. However, our QAP formulation is compatible with non-integer and non-positive problem data, which will be useful in the theoretical analysis and feature definitions
    which follow.
    We include instances where the problem data matrices, and hence the QAP instance itself, are asymmetric. Therefore, we cannot apply theoretical results which apply only to the symmetric QAP.  

	\subsection{Instance Subset}\label{sec:instances}
    
    Our instance subset contains benchmark instances from the QAPLIB library \citep{Burkard1997QAPLIB} available at \citep{QAPLIBsite}, benchmark instances from other sources, instances constructed using generators proposed in the literature, and instances constructed using generators proposed in this work.  The metaheurstics we are comparing in this work are often capable of finding a single reasonable-quality solution of large instances in a reasonable amount of time, but in order to obtain a dependable measure of algorithm performance on each instance we need to apply each metaheuristic several times. Therefore, we limit our analysis to instances of size $n \leq 128$. We also omit the existing benchmark instances of size $n < 25$, since the QAPLIB benchmark instances of this size were easily solved to optimality by all algorithms under consideration.

    The instance subset of 443 QAP instances used to produce the initial instance space is summarised in Table \ref{tab:instances}. 
    In this table and the rest of the paper, we categorise the instances according to the classification scheme proposed by \citet{Taillard1995classify}. The parameters used for the instance generators are listed in Appendix A. 
    In our later analysis some of these instance classes which unusual structure will be particularly important, so we briefly discuss their properties here.


    \begin{table}[t]
	    \scriptsize
	    \centering
	    \caption{Instance classes of the 443 instances used in the initial Instance Space Analysis.}
	    \rowcolors{2}{gray!15}{white} 
	    \begin{tabular}{m {3.7cm} m {3.9cm} c m{5cm} }
	      \toprule
	      Category & Instance Class and Source & \# of Inst. & Description \\ \midrule
       
	      \textbf{Real-life data} & Benchmark, QAPLIB & 23 & Instance names: bur*, ste*, esc*, kra* \\ \midrule

            \textbf{`Real-like' instances} & Benchmark, QAPLIB & 8 & Instance names: tai*b \\
            & Generator,\newline \citet{Stutzle2004Instances} & 60 & Flows may be random, structured, or structured plus noise \\ 
            & Generator, tai*b \newline \citep{Taillard1995classify} & 40 & Including variant with asymmetric distances; see Appendix A.2 \\\midrule

            \textbf{Grid-based distances} & Benchmark, QAPLIB & 21 & Instance names: nug*, sko*, tho*, wil* \\
             & Generator, \newline \citet{Stutzle2004Instances} & 60 & Flows may be random, structured, or structured plus noise \\ \midrule

            \textbf{Uniform random data} & Benchmark, QAPLIB  & 8 & Instance names: tai*a  \\
            & Generator, Random & 40 & Includes symmetric and asymmetric instances \\ \midrule

            \textbf{Other} & Benchmark, QAPLIB & 16 & Instance names: chr*, lipa*, tai*c \\
            & Benchmark, drexx\newline Section 2.1 of \newline\citep{DreznerTaillard2005} & 7 & Similar to the Grid-based distance instances with flows only linking adjacent facilities \\
            & Generator, Terminal\newline Section 2.2 of \newline \citep{DreznerTaillard2005} & 40 & Distances based on tree with uniform-size clusters of leaves, flows link facility clusters \\
            & Generator, Hypercube \newline (Section~\ref{sec:hypercube}) & 40 & Distances based on vertices of one or more hypercubes, flows match hypercube edges \\
            & Generator,\newline\citet{Palubeckis2000generate} & 40 & Grid-based distances, construction of flow matrix guarantees known optimal solution \\
            & Generator, QAPSAT\newline\citep{Verel2024} & 40 & Based on the structure of a 3-SAT problem, Contains `easy' and `hard' subcategories \\
            \bottomrule
	    \end{tabular}
     \label{tab:instances}
    \end{table}

    \subsubsection{Instances with few high-quality solutions}


    The Terminal-type instances were proposed by \citet{DreznerTaillard2005} as the Taixxeyy instance class. (We rename them here to avoid conflict with the Taixxeyy instances available at \citep{TaillardSite}, which appear to have a different structure.) The distance matrix represents a branching tree structure (like an airport terminal) in which every location is identified with a leaf of the tree. The flow matrix has a similar but inverted structure, so that solvers are strongly incentivised to place strongly linked facilities together in a single branch of the tree. 
    
    The Drexx instances were also proposed by \citet{DreznerTaillard2005}; these instances are based on a rectangular grid, and there are four optimal solutions corresponding to reflections of the rectangle. In our instance subset we include the benchmark instances available from \cite{TaillardSite}.

    The instance generator proposed by \citet{Palubeckis2000generate} randomly selects points on a grid, then takes the Manhattan distance between them to generate a distance matrix. A flow matrix with special structure is then produced to guarantee that the resulting instance has a single known optimal solution.

    \subsubsection{Hypercube} \label{sec:hypercube}
    We propose in this work an additional class of instances with few high-quality solutions, which we call the Hypercube instances. These are similar to the Drexx instances described in the previous section, generalised to a grid of more than two dimensions and with a modified distance metric. Our generator has the following integer-valued input parameters: dimension $k > 0$, side length $\ell > 1$, base distance $d_0 > 0$ and base flow $f_0 > 0$.
    The instance will have $\ell^k$ facilities and locations. 
    
    A Hypercube instance is constructed by the following procedure. Let $V = \left\{0,\dots,\ell-1 \right\}^k$. Observe that $V$ represents the vertices of a hypercube in $k$ dimensions if $\ell=2$, or the vertices of a grid of hypercubes if $\ell > 2$, and $|V| = \ell^k$. Denote the $i$th element of $V$ as $v_i$ (the ordering here determines the order of facilities and locations in the instance). Then the problem data matrices $D = [d_{ij}]$ and $F = [f_{ij}]$ are defined as follows, where $\disunif{a}{b}$ represents an integer in $\left\{ a, \dots, b \right\}$ independently selected for each $(i,j)$ using a uniform random distribution:
    $$
        d_{ij} = 
        \begin{cases}
            0 & \left|v_i-v_j\right|_1 = 0 \\
            \disunif{1}{d_0} & \left|v_i-v_j\right|_1 = 1 \\
            \disunif{2 d_0 + 1}{2 d_0 + \frac{d_0}{2}} & \left|v_i-v_j\right|_1 = 2 \\
            \disunif{4 d_0 + 1}{4 d_0 + \frac{d_0}{2}} & \left|v_i-v_j\right|_1 > 2
        \end{cases}
    $$
    $$
        f_{ij} = 
        \begin{cases}
            \disunif{f_0 + 1}{f_0 + \frac{f_0}{2}} & \left|v_i-v_j\right|_1 = 1 \\
            0 & \mathrm{otherwise} \\
        \end{cases}
    $$

    The Hypercube instances describe a problem in which we are mapping the vertices of a hypercube, or grid of hypercubes, to the vertices of another hypercube or grid with the same dimensions, in such a way that adjacency of all pairs of vertices is preserved. If we select a corner facility, place it at one of the $2^k$ corner locations, and then place its $k$ adjacent facilities at the $k$ adjacent locations, this uniquely determines a potential orientation of the hypercube and hence a high-quality solution of the problem. Therefore there are $2^k k!$ high-quality solutions, any one of which may be optimal based on the random element of the problem data. Any solution in which the arrangement of the facilities does not match the hypercube structure defined by the distance matrix will incur a large penalty. Similarly to the drexx instances, and in contrast to the Terminal instances, every pair of high-quality solutions of a Hypercube instance has few assignments in common.

	\subsection{Algorithm Space}\label{sec:algorithms}


        In this study we compare two metaheuristics for the QAP, which are described in this section.

        The Breakout Memetic Algorithm (BMA) \citep{BenlicBMA} applies the general framework of a genetic algorithm, with the additional step of applying Breakout Local Search \citep{BenlicBLS} to the initial population and to each new individual produced by crossover or mutation to find a local minimum solution. The BLS algorithm can be applied to solving QAP instances on its own, but we observed that for our instance subset it is typically outperformed by BMA. We use the implementation available at \citep{BMAsite} and the default parameters.

        The Max-Min Ant System (MMAS) algorithm \citep{Lopez2018MMAS, StutzleHoosMMAS} is an Ant Colony Optimisation (ACO) algorithm. At each iteration of an ACO algorithm, several agents known as ants construct new solutions to the QAP instance, guided by the pheromones deposited by successful ants in previous iterations. MMAS differs from other ACO algorithms in several important ways relating to the handling of pheromone values. The pheromones associated with each assignment begin at the maximum level, they decay to a lower but non-zero level if not reinforced, and may be reset to the original values if a failure to progress is detected. Only a single ant is chosen to deposit new pheromones at each iteration; depending on the state of the algorithm this may be the ant with the best solution in the current iteration, the best solution since the last pheromone reset, or the best solution since the beginning of the algorithm.  Similarly to the BMA algorithm, a local search algorithm is applied to the solutions constructed by the ants to find a local minimum solution.  We use the general ACOTSPQAP implementation available at \cite{MMASweb} with the default parameters, which set the main algorithm to MMAS and the local search to best 2-opt.

        Our selection of these metaheuristics is driven by three considerations. First, the algorithms are competitive with other state-of-the-art algorithms and each other \citep{Dantas2018label, Dantas2020landscape}.
        Second, they have significantly different designs and strengths, so we may draw interesting conclusions from our analysis. Third, efficient and well-tuned implementations of these algorithms are available for academic use, which supports reproducibility of results.     

	\subsection{Performance Space}\label{sec:performance}


    To compare the performance of the algorithms we ran them on the University of Melbourne's Spartan HPC system, using one compute node with one Intel(R) Xeon(R) Gold 6254 CPU @ 3.10GHz and 16GB RAM per test.
 
    Since the algorithms which we compare in this work are heuristics rather than exact algorithms, the performance of each algorithm is determined by two quantities: the quality of the solution produced by the algorithm, and the resources (for our purposes, time) required to produce that solution. To address this we apply a similar method to \citet{Dantas2020landscape}, giving each algorithm a fixed maximum runtime $t(n)$ based on instance size $n$ and comparing the best solution found by each algorithm in that time. To determine a reasonable definition for $t(n)$, we ran the BMA and MMAS algorithms $10$ times on each of the instances from the QAPLIB instances and the instances produced by all generators for the `real-like', grid-based and uniform random instance classes, with a maximum runtime of $2400$ seconds.
    We then plotted, for each instance and both algorithms, the average runtime (in seconds) required for the algorithms to obtain their best solution against instance size, then determined the parameters $a$ and $b$ of the exponential function $a e^{bn}$ which best fit this data. Finally, we defined the maximum runtime function as $t(n) := \max\left\{ae^{bn}, 5\right\}$ seconds. The best-fit parameters were $a = 56.45$ and $b =0.0193$.

    A natural approach to assessing the quality of a solution produced by a heuristic is to compare it to the optimal value. Unfortunately, in our instance subset the optimal values are not known for several instances, including some of the larger library problems and most of the generated instances.
    Rather than compare our heuristics' results with best-known solutions of potentially dubious quality, we instead choose to directly compare the results of the algorithms. Therefore, our results focus on the question of ``what are the relative strengths and weaknesses of these algorithms?'' rather than ``in an absolute sense, which of these instances are difficult to solve?''. We note that the ISA methodology would be capable of addressing the latter question using an alternative performance metric if the optimal solutions for the instances in our subset were known.

    The performance of each algorithm on each QAP instance is defined as follows. For a given instance of size $n$ we run each algorithm $50$ times for $t(n)$ seconds, and record the value of the best solution found in each run. We then record the average of the solutions found by each algorithm (denoted by $\mu_{BMA}$ and $\mu_{MMAS}$), and the standard deviation of the combined set of $100$ solutions (denoted by $\sigma$). Finally, the performance measure of algorithm $A$ is defined as
    \begin{equation} \label{eq:perf} \textrm{perf}(A) := \frac{\mu_A  - \min \left\{ \mu_{BMA}, \mu_{MMAS} \right\}}{\sigma}\end{equation} if $\sigma > 0$, and $\textrm{perf}(A) = 0$ otherwise. This means that the best-performing algorithm on an instance has performance measure $0$, and algorithms which perform worse have a performance measure greater than $0$. We designate an algorithm with performance less than or equal to $1$ as having good performance on the instance in question, and an algorithm with performance greater than $1$ as having bad performance.
    
    When comparing the two algorithms directly we will also refer to an overall performance criterion 
    \begin{equation} \label{eq:alg-crit}\textrm{algperf} := \textrm{perf}(MMAS) - \textrm{perf}(BMA). \end{equation} If $\textrm{algperf} > 1$ BMA is clearly the better algorithm, and if $\textrm{algperf} < -1$ MMAS is clearly better for the instance in question.

	\subsection{Theoretical Foundations for QAP Feature Design}\label{sec:all-theory}
    The final component required to apply ISA is a set of features to be measured for each instance. In this section we present some theoretical properties of the QAP. Many of these results are already well known in the literature; we reproduce them here to clarify their role in motivating and justifying our feature definitions. 
    In this and following sections, $I_n$ is the $n \times n$ identity matrix and $\ones_n$ is the $n \times n$ matrix with every entry equal to $1$. We define $\Nneq^2$ as the subset of $N^2$ containing no 2-tuples with repeating elements:
 \begin{equation*}\label{eq:define-M}
 \Nneq^2 := \left\{ (i_1,i_2) \in N^2 \mid i_1 \neq i_2 \right\}
 \end{equation*}
 
	
	\subsubsection{Similar and Equivalent Instances} \label{sec:qap-theory}
	

        In the following, $\sigma_{i,j}$ with $i,j \in N$ and $i \neq j$ denotes the permutation of $N$ which transposes $i$ and $j$: 
        \begin{equation*}
        \sigma_{i,j}(k) =
        \begin{cases}
            j & \text{if } k=i\\
            i & \text{if } k=j\\
            k & \text{otherwise}
        \end{cases}
        \end{equation*}
        
        $d_C(\varphi, \pi)$ denotes the Cayley distance between two permutations $\varphi$ and $\pi$ i.e. the minimum number of transpositions required to turn $\varphi$ into $\pi$ or vice versa. For example, if $\varphi \neq \pi$ and $\varphi \neq \sigma_{i,j} \circ \pi$ for all $i \neq j$ in $N$, but $\varphi = \sigma_{i,j} \circ \sigma_{k,l} \circ \pi$ for some $i \neq j$ and $k \neq l$ in $N$, then $d_C(\varphi, \pi) = 2$. 

        \begin{definition} [Automorphism]
            A function $f : S_n \to S_n$ is an \textbf{automorphism} of $S_n$ if it is bijective and satisfies the condition 
            $$d_C(\varphi,\pi) = 1 \iff d_C(f(\varphi),f(\pi)) = 1 \qquad \forall \varphi, \pi \in S_n.$$
        \end{definition}
        Note that the inverse of an automorphism, and a composition of automorphisms, are also automorphisms. 

        \begin{definition} [Similar and equivalent QAP instances]
        We describe a pair of QAP instances $QAP(A_1,B_1)$ and $QAP(A_2,B_2)$ of equal size $n$ as being \textbf{similar} if there exists a constant $c_1 > 0$, a constant $c_2 \in \mathbb{R}$, and an automorphism $f$ of $S_n$, such that 
        \begin{equation*}Q_{A_1,B_1}(\varphi) = c_1 Q_{A_2,B_2}(f(\varphi)) + c_2 \qquad \forall \varphi \in S_n.
        \end{equation*}
        If $f$ is the identity function i.e. $f(\varphi) = \varphi$ then we also describe the pair of instances as \textbf{equivalent}.
	\end{definition}

        Note that similarity and equivalence are both equivalence relations.

        The QAP algorithms which we compare in this work both have the same basic idea of exploring the decision space, which is naturally represented as a network of permutations linked by the potential transpositions of assignments. Therefore, it is unsurprising that when we replace an instance with a similar instance the performance of the algorithms is unaffected. We present computational results justifying this assertion in the supplemental material (Appendix B).
        
        In this section we consider several conditions which are sufficient to guarantee similarity of two QAP instances, and comment on how they will affect our choice of features. 

        The results presented in the rest of this section are straightforward and well known in the literature, so in the interests of brevity we provide proofs in the supplemental material (Appendix C). 
        In the following propositions we assume that the matrices under consideration are arbitrary real-valued $n \times n$ matrices with no particular properties beyond those explicitly stated.

    \begin{restatable}{proposition}{propreorder} \label{prop:reorder}
        Let $\tau, \theta \in S_n$. Let $A = [a_{ij}]$ and $B = [b_{ij}]$ be arbitrary real-valued $n \times n$ matrices and let $X = [x_{ij}]$ and $Y = [y_{ij}]$ be defined by $x_{ij} = a_{\tau(i)\tau(j)}$  and $y_{ij} = b_{\theta(i)\theta(j)}$ for all $i,j$ in $N$. Then $QAP(A,B)$ and $QAP(X,Y)$ are similar QAP instances. 
    \end{restatable}

    This result means that reordering the facilities and/or locations in the problem data, obtaining a pair of isomorphic QAP instances in the sense defined by \citet{Abreu2002iso}, also results in a pair of similar QAP instances. Happily the natural definitions of our features are invariant with respect to such a reordering.
    
    \begin{restatable}{proposition}{propinverse} \label{prop:inverse-problem}
        $QAP(A,B)$ and $QAP(B,A)$ are similar QAP instances.
    \end{restatable}

    Therefore, we should define features for QAP such that the feature vectors corresponding to $QAP(A,B)$ and $QAP(B,A)$ are identical. We discuss this further in Section~\ref{sec:identify-dist}.

    \begin{restatable}{proposition}{propsymm} \label{prop:either-symmetric}
	    Let $B$ be a symmetric matrix (i.e. $b_{ij} = b_{ji}$ for all $i,j$ in $N^2$). Let $C = \left[c_{ij}\right]$ be a skew-symmetric $n \times n$ matrix (i.e. $c_{ij} = - c_{ji}$ and $c_{ii} = 0$ for all $(i,j)$ in $M^2$). Then for all $\varphi \in S_n$, $$\costf_{A,B} (\varphi) = \costf_{(A+C),B}(\varphi)$$ and hence $QAP(A,B)$ and $QAP(A+C,B)$ are equivalent QAP instances.
    \end{restatable}
	In particular, if $C = \frac{1}{2}(A^\Tr - A)$, then $A+C =  \frac{1}{2}(A + A^\Tr)$ is a symmetric matrix. A similar result applies to $\costf_{A,(B+C)}$ if $A$ is symmetric. Therefore, any QAP instance in which one matrix is symmetric is equivalent to an instance in which both matrices are symmetric, so we should only consider an instance to have meaningful asymmetry if both problem data matrices are asymmetric.
	
	\begin{restatable}{proposition}{propmult}\label{prop:multiply-constants}
	    Let $c_1$ and $c_2$ be arbitrary positive constants in $\mathbb{R}$. Then for all $\varphi \in S_n$
	    $$\costf_{c_1 A,c_2 B}(\varphi) = c_1 c_2 \costf_{A,B}(\varphi)$$
	    and hence $QAP(A,B)$ and $QAP(c_1 A, c_2 B)$ are equivalent QAP instances.
	\end{restatable}
	\begin{restatable}{proposition}{propaddI}\label{prop:add-constants-diagonal}
	    Let $c_1$ and $c_2$ be arbitrary constants in $\mathbb{R}$. Then 
	    for all $\varphi \in S_n$
	    $$\costf_{(A+c_1 I_n),(B+c_2 I_n)}(\varphi) = \costf_{A,B}(\varphi) + c_2 \trace{A} + c_1 \trace{B} + n c_1 c_2$$
	    and hence $QAP(A,B)$ and $QAP(A+c_1 I_n, B+c_2 I_n)$ are equivalent QAP instances.
	\end{restatable}
	\begin{restatable}{proposition}{propaddJ}\label{prop:add-constants-elsewhere}
	    Let $c_1$ and $c_2$ be arbitrary constants in $\mathbb{R}$. Then 
	    for all $\varphi \in S_n$
	    $$\costf_{A+c_1 (\ones_n - I_n),B+c_2 (\ones_n - I_n)}(\varphi) = \costf_{A,B}(\varphi) + \sum_{(i,j) \in \Nneq^2} a_{ij} c_2 + \sum_{(i,j) \in \Nneq^2} b_{ij} c_1 + n (n-1) c_1 c_2$$
	    and hence $QAP(A,B)$ and $QAP(A+c_1 (\ones_n - I_n), B+c_2 (\ones_n - I_n))$ are equivalent QAP instances.
	\end{restatable}
	
	Proposition \ref{prop:multiply-constants} tells us that the absolute magnitude of the entries in $A$ and $B$ does not affect the difficulty of the instance. Furthermore, Propositions \ref{prop:add-constants-diagonal} and \ref{prop:add-constants-elsewhere} indicate that modifying $A$ and/or $B$ by adding or subtracting a constant from the elements on the main diagonal of these matrices, and/or from the elements not on the main diagonal, also results in an equivalent QAP instance. The instance normalisation procedure discussed in Section~\ref{sec:standard-form} addresses these results.
 

    \begin{restatable}{proposition}{propneg} \label{prop:negative-problem}
        For all $\varphi \in S_n$ $$\costf_{-A,-B}(\varphi) = \costf_{A,B}(\varphi)$$ and hence $QAP(A,B)$ and $QAP(-A,-B)$ are equivalent QAP instances.
    \end{restatable}

    Proposition~\ref{prop:negative-problem}, in combination with the preceding propositions, illustrates the notion that assigning a small flow to a short distance (indirectly) harms solution quality, since this assignment wastes the capability of the short distance to mitigate the cost of a large flow (or vice versa).

    \subsubsection{Average and Standard Deviation of Costs} \label{sec:qap-stats-1}

    It is well known that the average solution cost for a QAP instance, $E \left[\costf_{A,B}(\varphi)\right]$, can be calculated in $O(n^2)$ time:
    \begin{equation}
	    E \left[\costf_{A,B}(\varphi)\right] = \frac{1}{n(n-1)} \left( \sum_{(i,j) \in \Nneq^2} a_{ij} \right) \left( \sum_{(r,s) \in \Nneq^2} b_{rs} \right) + \frac{1}{n} \left( \sum_{i\in N} a_{ii} \right) \left( \sum_{r \in N} b_{rr} \right). \label{eq:qap-mean}
    \end{equation}
    A proof can be found in e.g. Proposition 5 of \citet{Chicano2012698}. Expressions for calculating the second moment $E\left[\left(\costf_{A,B}(\varphi)\right)^2\right]$ in polynomial time were found for symmetric QAP by \citet{Angel2001landscape} and \citet{Abreu2002iso}.  \citet{Chicano2012698} generalises this result to the asymmetric QAP in their Proposition 6, and moreover observed that it too can be calculated in $O(n^2)$ time.

    Combining these results, we can compute the standard deviation of all solution costs for $QAP(A,B)$ as the square root of the variance:
	\begin{equation} \label{eq:qap-stdev}
	    \sigma = \sqrt{\textrm{Var}\left(\costf_{A,B}(\varphi)\right)} = \sqrt{E\left[\left(\costf_{A,B}(\varphi)\right)^2\right] - E\left[\costf_{A,B}(\varphi)\right]^2}
	\end{equation}
    As discussed in Section~\ref{sec:standard-form} we will in practice apply this result to instances where $E\left[\costf_{A,B}(\varphi)\right] = 0$, which has the advantage here of removing a potential source of numerical error due to catastrophic cancellation.

    \subsubsection{Average Cost after Fixing One Assignment} \label{sec:qap-stats-3}

    Theorem 2 of \citet{ChmielCondExp2019} gives an expression for the average cost of a QAP instance under the condition that one or more assignments are predefined and cannot be altered. For a single specified list of assignments to be fixed, this expression can be evaluated in $O(n^2)$ time.
    In this section we restrict our attention to the cases where only one assignment is fixed. We demonstrate a more efficient method for calculating all $n^2$ of these conditional expected values in $O(n^2)$ time, rather than the $O(n^4)$ time which would be required to naively apply Chmiel's result to each of the $n^2$ cases.

    Formally, for all $x$ and $y \in N$ we define $S_n^{xy}$ as the set of all permutations $\varphi \in S_n$ satisfying the condition $\varphi(x)=y$. For all $x\in N$ we define $N_x = N \setminus \left\{x\right\}$, and define $\Nneq_x^2$ as the subset of $N_x^2$ which contains no pairs with repeated values (or, equivalently, the subset of $M^2$ containing no pairs which contain $x$).

        \begin{proposition} \label{prop:fix-one-mean}
            For all $x \in N$ and $y \in N$, the average cost of a solution $\varphi$ of $QAP(A,B)$, where $\varphi$ is selected randomly from $S_n^{xy}$ (in effect fixing $\varphi(x) = y$), can be expressed as
            { \small
            \begin{multline} \label{eq:fix-one-mean}
                E \left[\costf_{A,B}(\varphi) \mid \varphi \in S_n^{xy}\right] = \frac{1}{(n-1)(n-2)} \left( \sum_{(i,j) \in M_x^2} a_{ij} \right) \left( \sum_{(r,s) \in M_y^2} b_{rs} \right) + \frac{1}{n-1} \left( \sum_{i \in N_x} a_{ii} \right) \left( \sum_{r \in N_y} b_{rr} \right) + \\ \frac{1}{n-1} \left( \sum_{i \in N_x} a_{ix} \right) \left( \sum_{r \in N_y} b_{ry} \right)
                + \frac{1}{n-1} \left( \sum_{j \in N_x} a_{xj} \right) \left( \sum_{s \in N_y} b_{ys} \right) + a_{xx} b_{yy}
            \end{multline}
            }
        \end{proposition}
        \begin{proof}
            {\small 
            \begin{eqnarray*}
    	    E\left[\left(\costf_{A,B}(\varphi) \mid \varphi \in S_n^{xy}\right)\right] &=& \frac{1}{\left|S_n^{xy}\right|}  \sum_{\varphi \in S_n^{xy}} \left(\costf_{A,B}(\varphi)\right) \nonumber \\
                &=& \frac{1}{(n-1)!} \sum_{\varphi \in S_n^{xy}} \left( \sum_{(i,j) \in M_x^2} a_{ij} b_{\varphi(i)\varphi(j)} + \sum_{i \in N_x} a_{ii} b_{\varphi(i)\varphi(i)} + \sum_{j \in N_x} a_{xj} b_{y\varphi(j)} + \right. \nonumber \\ & & \left. \sum_{i \in N_x} a_{ix} b_{\varphi(i)y} + a_{xx} b_{yy} \right)  \nonumber \
            \end{eqnarray*}
            }
            
            Here we note that for fixed $(i,j) \in M_x^2$ and $(r,s) \in M_x^2$ in  there are exactly $(n-3)!$ permutations $\varphi \in S^{xy}_n$ satisfying the conditions $\varphi(i) = r$ and $\varphi(j) = s$. Using this and similar observations, we rearrange the sums then replace the sums over $S^{xy}_n$ with sums indexed by $r$ and $s$ to obtain the desired result:
            
            {
            \small
            
            \begin{eqnarray*}
    	    E\left[\left(\costf_{A,B}(\varphi) \mid \varphi \in S_n^{xy}\right)\right] &=& \frac{1}{(n-1)!}  \left( \sum_{(i,j) \in M_x^2} \sum_{\varphi \in S_n^{xy}} a_{ij} b_{\varphi(i)\varphi(j)} + \sum_{i \in N_x} \sum_{\varphi \in S_n^{xy}} a_{ii} b_{\varphi(i)\varphi(i)} + \right. \nonumber \\ & & \left. \sum_{j \in N_x} \sum_{\varphi \in S_n^{xy}} a_{xj} b_{y\varphi(j)} + \sum_{i \in N_x} \sum_{\varphi \in S_n^{xy}} a_{ix} b_{\varphi(i)y} + \sum_{\varphi \in S_n^{xy}} a_{xx} b_{yy} \right) \nonumber \\
            &=& \frac{1}{(n-1)(n-2)} \left( \sum_{(i,j) \in M_x^2} a_{ij} \right) \left( \sum_{(r,s) \in M_y^2} b_{rs} \right) + \frac{1}{n-1} \left( \sum_{i \in N_x} a_{ii} \right) \left( \sum_{r \in N_y} b_{rr} \right) + \nonumber \\ & & \frac{1}{n-1} \left( \sum_{i \in N_x} a_{ix} \right) \left( \sum_{r \in N_y} b_{ry} \right)
                + \frac{1}{n-1} \left( \sum_{j \in N_x} a_{xj} \right) \left( \sum_{s \in N_y} b_{ys} \right) + a_{xx} b_{yy} \nonumber
            \end{eqnarray*}
            }
        \end{proof}

    We give particular attention to the first term of \eqref{eq:fix-one-mean}:
    $$\frac{1}{(n-1)(n-2)} \left( \sum_{(i,j) \in M_x^2} a_{ij} \right) \left( \sum_{(r,s) \in M_y^2} b_{rs} \right)$$
    For a fixed $x \in N$ the sum over $M_x^2$ requires $O(n^2)$ time to evaluate. A naive attempt to evaluate this sum for all $x \in N$ would require $O(n^3)$ time, but we can instead decompose the sum as follows:
    $$\sum_{(i,j) \in M_x^2} a_{ij} = \sum_{(i,j) \in M^2} a_{ij} - \sum_{j \in N_x} a_{xj} - \sum_{i \in N_x} a_{ix}$$
    In this expression the first sum is over a set of size $O(n^2)$ but is not parameterised by $x$; the remaining sums are parameterised by $x$ but are over sets of size $O(n)$. Therefore we can evaluate all of them, for all $x \in N$, in $O(n^2)$ time. A similar method can be applied to decompose the sums over $M_y^2$. All of the other sums in Proposition~\ref{prop:fix-one-mean} are over sets of size $O(n)$, and evaluating them for all $x,y \in N$ requires at most $O(n^2)$ time. Therefore, we can evaluate $E \left[\costf_{A,B}(\varphi) \mid \varphi \in S_n^{xy}\right]$ for all $x,y \in N$ simultaneously in $O(n^2)$ time.

    \subsection{Feature Space}\label{sec:features}

    The features described in this section are defined with three considerations in mind. 
    
    First, they should (singly or in combination) be predictive of instance difficulty for the algorithms being studied. In particular we desire to avoid a situation in which two similar or equivalent instances have widely disparate feature values and hence conflicting performance predictions. As we have already observed, actual algorithm performance remains consistent when comparing similar or equivalent instances.
    
    Second, it should be possible to evaluate the features of an instance in considerably less time than we would expect the algorithms to run for, to facilitate the construction of an algorithm selector with practical utility. Third, the value of the features should be consistently measurable. The second and third points come into tension when we consider fitness landscape features based on random sampling of the solution space; a larger number of samples increases feature reliability but also increases the computational cost.

    \subsubsection{Identifying the Problem Matrices} \label{sec:identify-dist}

    As per the comments following Proposition \ref{prop:inverse-problem}, exchanging the order of the problem data matrices results in a similar QAP instance. Therefore, to obtain consistent feature definitions we must identify the matrices in a way which is invariant with respect to their order. The literature typically describes a QAP instance as having a `distance' matrix and a `flow' matrix, with the identification being made by inspection or based on the statements of the instance's author. Neither of these approaches is suitable in general when we may encounter a number of instances too large to classify by hand, or esoteric instance classes where neither matrix fits our intuitive understanding of `distances' or `flows'. Therefore, in this section we propose an algorithm to efficiently identify the matrices of a QAP instance.

    The structure of a flow matrix may vary significantly based on the characteristics of the problem being modelled by a particular QAP instance. Instead, we attempt to identify the distance matrix based on whether it appears to respect the triangle inequality. This procedure is complicated by Proposition~\ref{prop:add-constants-elsewhere}, since for any QAP instance we can trivially create an equivalent instance which respects a naive application of the triangle inequality by adding a large constant to all matrix entries. We would also prefer our measure to be robust against a small number of outlier terms in the matrix. 
    
    To address these issues we propose the TRiangle Inequality Property Of Distances (TRIPOD) score, for which we provide algorithm pseudocode in the supplemental material (Appendix D.1). 
    This algorithm measures the degree to which a single input matrix resembles a distance matrix, and has two parameters, $\alpha$ and $\beta$. We begin by rescaling the problem data so that matrix entries which were $\alpha$ standard deviations or more below the mean in the original matrix are zero in the new matrix. We then examine each triplet of locations, applying a penalty to the TRIPOD score if the corresponding distances in the new matrix do not respect the triangle inequality. The maximum penalty which can be inflicted by a single triplet is determined by the parameter $\beta$. We chose the parameter values $\alpha=3$ and $\beta=10$. With these settings the TRIPOD score correctly identifies (i.e. is larger for) the distance matrix for all QAPLIB instances for which the authors state one matrix should be interpreted as defining distances, which establishes face validity for the procedure.

    \subsubsection{Standard Form for Problem Matrices} \label{sec:standard-form}

    In this section we define a standard form for the problem data matrices $A$ and $B$ with the property that any pair of QAP instances which can be shown to be equivalent by applying Propositions~\ref{prop:either-symmetric}, \ref{prop:multiply-constants}, \ref{prop:add-constants-diagonal} and \ref{prop:add-constants-elsewhere} will have identical data matrices in standard form. The standard form has the following properties:
    \begin{itemize}
        \item The matrices are either both symmetric or both asymmetric.
        \item The mean value of all entries in both matrices, and the mean cost of all possible assignments, is 0.
        \item The standard deviation of the costs of all possible assignments is 1.
        \item The standard deviation of the entries in the two matrices is identical.
        \item One of the two matrices is identified as a distance matrix using the TRIPOD score.
    \end{itemize}

    \begin{definition} [Zero-mean form]
	We define the zero-mean form $\zeromean(A)$ of a QAP data matrix $A$ as
	$$ \zeromean(A) := A - \left[ \frac{1}{n} \sum_{i \in N} a_{ii} \right] I_n - \left[ \frac{1}{n(n-1)} \sum_{i=1}^n\sum_{\substack{j=1\\j \neq i}}^n a_{ij} \right] (\ones_n - I_n) $$
    \end{definition}
	
	The zero-mean form represents the QAP instance obtained from $QAP(A,B)$ by applying Propositions \ref{prop:add-constants-diagonal} and \ref{prop:add-constants-elsewhere} such that the mean of the values in both matrices is zero. Furthermore, the trace of both matrices, and the sum of all values outside the main diagonal of both matrices, are also zero. By inspection of \eqref{eq:qap-mean} we have $E[Q_{\zeromean(A),\zeromean(B)}(\varphi)] = 0$.


    \begin{definition} [Standard form] \label{def:standard-form}
    Given a pair of QAP problem data matrices $A$ and $B$, we define their standard form as the outputs $D$ and $F$ of the following procedure:
    \begin{enumerate}
        \item If $A$ is symmetric and $B$ is not, replace $B$ with $\frac{1}{2} \left( B + B^\Tr \right)$. Similarly, if $B$ is symmetric and $A$ is not, replace $A$ with $\frac{1}{2} \left( A + A^\Tr \right)$.
        \item Calculate the zero-mean forms $\zeromean(A)$ and $\zeromean(B)$ of the individual data matrices.
        \item Calculate the standard deviation $\sigma_Q$ for the problem $QAP(\zeromean(A),\zeromean(B))$ as in \eqref{eq:qap-stdev}.
        \item Calculate the standard deviations $\sigma_A$ and $\sigma_B$ of the entries in $\frac{1}{\sqrt{\sigma_Q}} \zeromean(A)$ and $\frac{1}{\sqrt{\sigma_Q}} \zeromean(B)$ respectively.
        \item Calculate the TRIPOD score of $\sqrt{\frac{\sigma_B}{\sigma_Q \sigma_A}} \zeromean(A)$ and $\sqrt{\frac{\sigma_A}{\sigma_Q \sigma_B}} \zeromean(B)$. Assign the matrix with the highest value for this measure as the distance matrix $D$, and the other matrix as the flow matrix $F$.
    \end{enumerate}
    \end{definition}

    For some features it is more natural to work with problem data matrices where the minimum value, rather than the mean value, of the matrices is zero.
 
    \begin{definition} [Zero-minimum form] \label{def:reduced-form}
	Given a QAP problem data matrix $A$, we define its zero-minimum form as
	$$ \reduced(A) := A - \left[ \min_{i \in N} a_{ii} \right] I_n - \left[ \min_{(i,j) \in N^2, i \neq j} a_{ij} \right] (\ones_n - I_n) $$
    \end{definition}

    We now present the list of features divided into three categories; features of the individual distance and flow matrices, features of the matrices combined, and fitness landscape features. For all features we assume that the problem data matrices have been placed in standard form, and then placed in zero-minimum form if indicated by the feature definition.

    In cases where the natural interpretation of a feature is unbounded, or the typical range of variation of the feature within our instance subset is considerably smaller than the theoretical range, we apply the following normalisation process:
	\begin{equation}\label{eq:feature-scaling}
	    \textrm{normalised feature} \gets \frac{\tan^{-1}(\textrm{raw feature value}/{\gamma})}{\pi/2}
	\end{equation}
	The value of $\gamma$ is specified independently for each feature where this normalisation is applied, with the goal that the mapping should be reasonably close to linear for typical raw feature values. 

	\subsubsection{Features applied to individual matrices} \label{sec:one-matrix-features}
	
	\begin{table}[t] 
	    \scriptsize
	    \centering
	    \caption{Description of features which we apply to each of the problem data matrices individually. Given initial problem data $QAP(A,B)$ we calculate the standard form $D$ and $F$. The features below are stated in terms of the distance matrix $D$; replace $D$ with $F$ for the flow matrix features which we also measure, for a total of 24 features in this category. Starred features are set to zero if the matrix contains only zeroes.}
	    \label{tab:features-one-matrix}
	    \rowcolors{2}{gray!15}{white} 
	    \begin{tabular}{m{3.0cm} c m{11cm} }
	      \toprule
	      Feature Name & Bound & Description \\ \midrule
	      
	      1. Normalised Mean$^*$ & $[0,1]$ & Mean of all elements of $\reduced(D)$, divided by $\max(\reduced(D))$. \\ 
		  
		  2. Trace Proportion$^*$ & $[0,1]$ & Sum of elements on the main diagonal of $\reduced(D)$, divided by the sum of all elements in $\reduced(D)$. \\
		  
		  3. Sparsity & $[0,1]$ & Proportion of elements in $\reduced(D)$ which are equal to zero. \\

		  
		  4. Dominance & [0,1] & Coefficient of variation of $\reduced(D)$, compared to that of reference matrices. See Equation~\eqref{eq:dominance}. \\
		  
		  5. TRIPOD score & $[0,1]$ & As described in Section~\ref{sec:identify-dist} and Appendix D.1. \\ 

            6. Diversity & [0,1] & Measure of diversity in elements of $D$. See Appendix D.2.\\ 

            7. Outliers & [0,1] & Proportion of entries in $D$ which are outside three standard deviations from the mean. \\

            8. Skewness & [0,1] & Skewness of elements of $D$. Rescaled using \eqref{eq:feature-scaling} with $\gamma=2$. \\

            9. Kurtosis & [0,1] & Kurtosis of elements of $D$. Rescaled using \eqref{eq:feature-scaling} with $\gamma=10$. \\

            10. Betafit Alpha & [0,1] & Linearly rescale all elements of $D$ to the interval [0,1]; then this feature is the alpha parameter of the best fit beta function to that data. Rescaled using \eqref{eq:feature-scaling} with $\gamma=1$. \\

            11. Betafit Beta & [0,1] & Linearly rescale all elements of $D$ to the interval [0,1]; then this feature is the beta parameter of the best fit beta function to that data. Rescaled using \eqref{eq:feature-scaling} with $\gamma=1$. \\

            12. Near Similarity & [0,1] & Treating $D$ as a distance matrix, measures whether or not two locations with a short distance between them tend to have similar properties. See Appendix D.3.\\
		  
		   \bottomrule
	    \end{tabular}
            \label{table:feat-single-matrix}
	\end{table}
	
	The features we apply to the distance and flow matrices individually are summarised in Table~\ref{table:feat-single-matrix}. 

    The \textbf{Dominance} feature is derived from the coefficient of variation of the values in the input matrix:
    \begin{equation*} \label{eq:coeff-variation}
        \textrm{CV}(A) = \frac{\textrm{standard deviation of values in A}}{\textrm{mean of values in A}}
    \end{equation*}
    Note that we apply this feature to the zero-minimum form of the matrix. To scale this feature we compute the coefficient of variation of two reference matrices of the same size as $A$. The high-dominance reference matrix $H_n$ has zeros everywhere except for one entry equal to $1$, while the low-dominance reference matrix $L_n$ is equal to $J_n - I_n$.
    The Dominance feature is then defined as
    \begin{equation} \label{eq:dominance}
        \text{Dominance(A)} := \frac{\textrm{CV}(H_n) - \textrm{CV}(\reduced(A))}{\textrm{CV}(H_n) - \textrm{CV}(L_n)}
    \end{equation}
    For all instances in our instance subset the value of this feature is between $0$ and $1$.
	
	\subsubsection{Exactly computed features of full QAP}
	
	\begin{table}[t] 
	    \scriptsize
	    \centering
	    \caption{Description of features applied to the two standard form QAP problem data matrices in combination. Given initial problem data $QAP(A,B)$ we calculate the standard form $D$ and $F$, and then apply these features to $QAP(D,F)$.}
	    \label{tab:features-exact}
	    \rowcolors{2}{gray!15}{white} 
	    \begin{tabular}{m{3.0cm} c m{11cm}  }
	      \toprule
	      Feature Name & Bound & Description \\ \midrule
	      
	      1. Instance size  & $[0,1]$ & Number of facilities/locations. Rescaled using \eqref{eq:feature-scaling} with $\gamma=200$. \\ 
		  
		  2. Maximum Symmetry & $[0,1]$ & Maximum of the symmetry measure calculated for the distance and flow matrices. \\

            3. Gilmore Lawler Bound & [0,1] & Lower bound proposed by \citet{Gilmore1962LB,Lawler1963LB}. Rescaled using \eqref{eq:feature-scaling} with $\gamma=-50$. \\

            4. Distribution Similarity & [0,1] & Measure of the similarity in the number of short distances and large flows, etc. See Equation~\eqref{eq:cumulative-integral}. \\

            5. Ruggedness Coefficient & [0,100] & Measure of the ruggedness of the solution space as defined by \citet{Angel2001landscape} and generalised to asymmetric QAP by \citet{Chicano2012698}. \\

            6. Least Dominance & [0,1] & The minimum of the distance dominance and flow dominance. \\

            7. Most Dominance & [0,1] & The maximum of the distance dominance and flow dominance. \\


            8. Improvement of $\mu$ from One Assignment & [0,0] & Measures the proportion of individual assignments which meaningfully improve average solution quality, using result from Section~\ref{sec:qap-stats-3}. \\

		   \bottomrule
	    \end{tabular}
            \label{table:feat-full-matrix}
	\end{table}

    The features we deterministically calculate using both matrices are listed in Table~\ref{table:feat-full-matrix}.

    The \textbf{Maximum Symmetry} feature is motivated by Proposition \ref{prop:either-symmetric}, which suggests that measuring the overall degree of symmetry of a QAP instance requires us to consider both matrices. The symmetry measure applied to each of the matrices is
    \begin{equation} \label{eq:symmetry}
        \textrm{Symmetry}(A) = \frac{\left| \left\{ (i,j) \in M^2 \mid a_{ij} = a_{ji} \right\} \right|}{n(n-1)}.
    \end{equation}

    The \textbf{Distribution Similarity} feature is intended to measure the `fit' between the distances and flows; an instance with a number of large flows about equal to the number of small distances is qualitatively different to an instance where these quantities are substantially different. Define the `cumulative distribution' function $C_A : \mathbb{R} \to [0,1]$ as follows: $$C_A(x) := \frac{\left| \left\{ (i,j) \in M^2 \mid a_{ij} \leq x \right\} \right|}{n(n-1)}$$
    Then the definition of the feature is
    \begin{equation} \label{eq:cumulative-integral}
    \textrm{DistributionSimilarity}(A,B) = \int_{-\infty}^{\infty} \left|C_A(x) - C_{-B}(x) \right| dx
    \end{equation}
 
    \subsubsection{Fitness landscape features using random sampling}

	\begin{table}[t]
	    \scriptsize
	    \centering
	    \caption{Description of fitness landscape features using random sampling, applied to the standard form matrices $D$ and $F$.}
	    \rowcolors{2}{gray!15}{white} 
	    \begin{tabular}{m{3.0cm} c m{11cm}  }
	      \toprule
	      Feature Name & Bound & Description \\ \midrule

            1. Average Distance to Optima  & [0,1] & Average distance to closest pseudo-global optimum, divided by instance size. \\ 
            
	      2. Fitness Distance Correlation  & [-1,1] & Correlation between distance to closest pseudo-global optimum and solution cost. \\ 

            3. Accumulated Escape Probability  & [0,1] & Average chance to improve solution with a single assignment swap. \\ 

            4. Base Dispersion Metric  & [-1,1] & From initial random sample, the dispersion of the best 5 percent of solutions minus the dispersion of the first 5 percent. \\ 

            5. Optima Dispersion Metric  & [-1,1] & From local optima found by steepest descent, the dispersion of the best 5 percent of solutions minus the dispersion of the first 5 percent. \\ 

            6. Average Descent  & [0,2] & Average number of steps taken by steepest descent before finding local optimum, divided by instance size. \\ 

            7. Optima Fitness Coefficient & [-1,0] & Standard deviation of the fitness of all local optima found, divided by their mean. \\

            8. Entropy Difference & [-1, 1] & Entropy of the initial sample, minus the entropy of the local optima; see \citet{Taillard1995classify} \\
       
		   \bottomrule
	    \end{tabular}
            \label{table:feat-fitness-landscape}
	\end{table}

	Our feature space also includes the fitness landscape features previously applied to the QAP by \citet{Dantas2020landscape}. To apply these features we begin by randomly generating a sample of $1000$ solutions, 
    and then apply a local steepest-descent search starting from each of these solutions to find a set of local minima. The lowest-cost solution or solutions among these local minima are referred to as `pseudo-global optima'. We then measure features of the original sample and pseudo-global optima, as summarised in Table~\ref{table:feat-fitness-landscape}.

    In addition to the fitness landscape features measured by \citet{Dantas2020landscape} and described more fully in that paper, we include an additional feature suggested by \citet{Taillard1995classify}, the \textbf{Entropy Difference}. For this feature we measure the entropy of the initial sample and the local optima, as defined by \citet{FleurentEntropy}, and take the differece between these entropy values.

We implemented the code required to measure these features using a MATLAB interface which calls a C++ library for the fitness landscape analysis, and ran this code on the same computational infrastructure as the algorithms (see Section~\ref{sec:performance} for details). We observe that the feature measurement including fitness landscape analysis carries a meaningful but not unreasonable computational cost. For example, for a Terminal instance of size $n=125$ it requires about $500$ seconds to measure the features. By comparison, we allocated about $630$ seconds to each algorithm per attempt to solve this instance. However, due to the high variance in the solution quality achieved in these tests, much more time would be required in practice to confidently achieve a high quality solution. In this context the computational cost of calculating our proposed features appears reasonable.  

In total we record 40 features: $12$ distance features, $12$ flow features, $8$ features measured deterministically using both matrices, and $8$ fitness landscape features measured using random sampling.
	\section{Constructing an Initial Instance Space}\label{sec:initial}


In this section we present and discuss the results of applying the ISA toolkit to the instance subset, performance data and feature data described in Section~\ref{sec:isa}. As discussed in more detail by \cite{Smith-Miles2023MATILDA}, the projection to the $2$D instance space is performed by four methods applied in sequence. The data is initially scaled, bounded and otherwise prepared by the PRELIM method. The SIFTED method then selects a subset of the input features which appear most relevant to explaining the performance of the algorithms on the instance subset. Following this, the PILOT method produces a projection from the features selected by PRELIM to the $2$D instance space, with the goal of producing visible trends in algorithm performance. Finally, the CLOISTER method infers a boundary of the instance space from experimental or theoretical bounds on each feature and the observed correlation between the features. Once the projection is in place, the PYTHIA method quantifies the predictive power of the projection by training a Support Vector Machine (SVM) for each algorithm to predict where in the instance space each algorithm will perform well. 

The parameters we applied to each of the ISA components can be found at \url{https://matilda.unimelb.edu.au/matilda/problems/opt/qap}. In particular, our choice to have SIFTED select 6 of the 40 features for this projection was based on visual inspection of the plots produced by varying this parameter, then choosing the projection appearing to best separate the instances for which one algorithm was clearly superior.

\subsection{Properties of the Instance Space} \label{sec:ispace-properties}

The ISA toolkit produces the following projection to construct the $2$D instance space:

\begin{equation}
		\begin{bmatrix} 
		Z_1 \\ 
		Z_2
		\end{bmatrix}
		= 
		\begin{bmatrix}[r]
			0.4430 &   -0.6851 \\
            -0.4695 &   0.2124 \\
            0.0509 &    0.3762 \\
            -0.7521 &   0.1012 \\
            -0.3867 &   -0.3110 \\
            0.1125 &   0.3864 \\
		\end{bmatrix}^\Tr
		\begin{bmatrix}
		\text{Distance Sparsity} \\
		\text{Distance TRIPOD Score} \\
		\text{Distance Betafit Alpha} \\
		\text{Distance Near Similarity} \\ %
		\text{Distribution Similarity} \\
		\text{Average Distance to Optima} \\
		\end{bmatrix}
		\label{eq:initialproj}
\end{equation}

\begin{figure}[t]
\centering
\begin{minipage}{.45\textwidth}
  \centering
	    \includegraphics[trim=0.0cm 0.0cm 0.0cm 0.0cm, clip=true, width=0.8\textwidth]{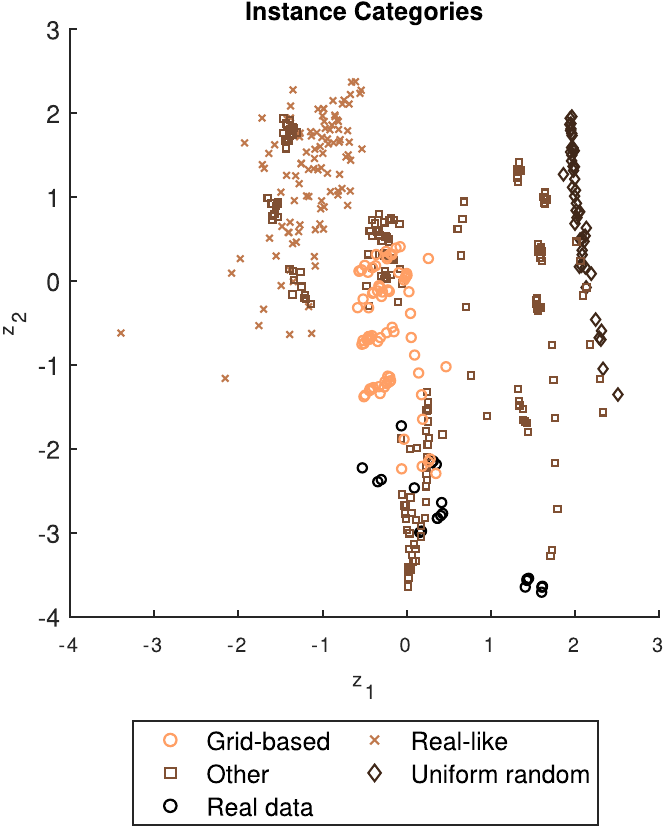}
	    \caption{Categories of instances in the initial QAP instance space, with axes as defined by~\eqref{eq:initialproj}.}
        \label{fig:initial-classes-macro}
\end{minipage}%
\hspace{0.3cm}
\begin{minipage}{.45\textwidth}
  \centering
	    \includegraphics[trim=0.0cm 0.0cm 0.0cm 0.0cm, clip=true, width=0.8\textwidth]{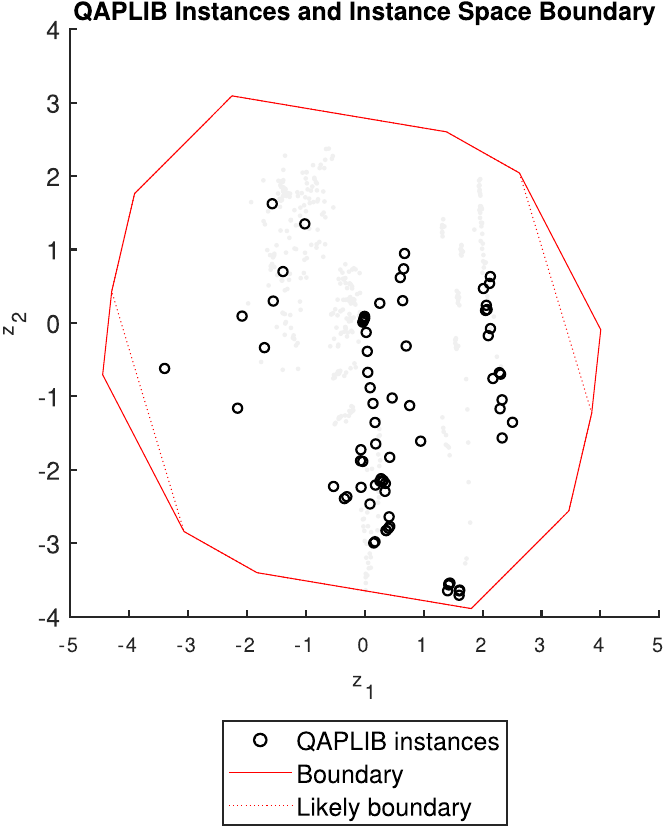}
	    \caption{Location of QAPLIB instances, and the boundary of the space inferred from observed variation in features.}
        \label{fig:initial-qaplib}
\end{minipage}
\end{figure}

\begin{figure}[t]
    \centering
    \subfloat{\includegraphics[trim=0.0cm -0.2cm 0.0cm 0.0cm, clip=true, width=0.24\textwidth]{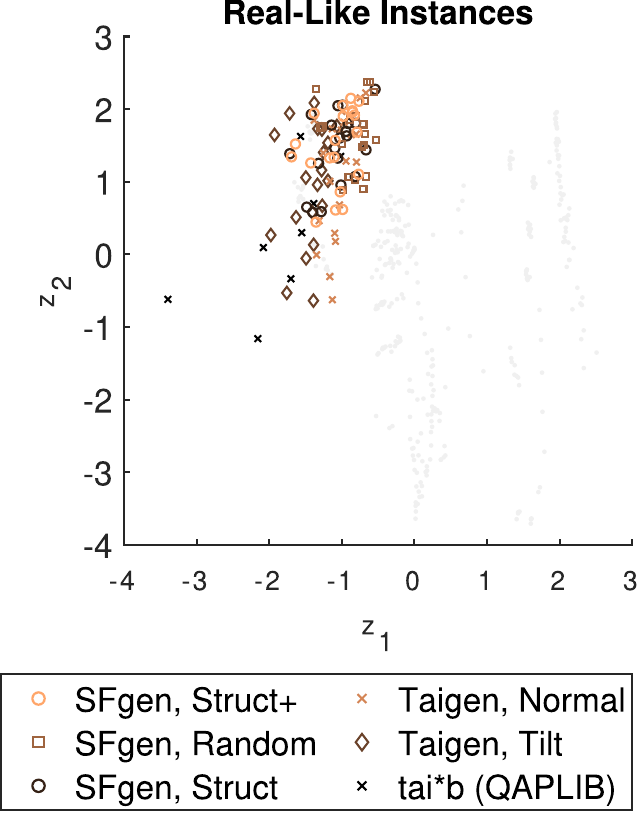}}%
    \subfloat{\includegraphics[trim=0.0cm 0.0cm 0.0cm 0.0cm, clip=true, width=0.24\textwidth]{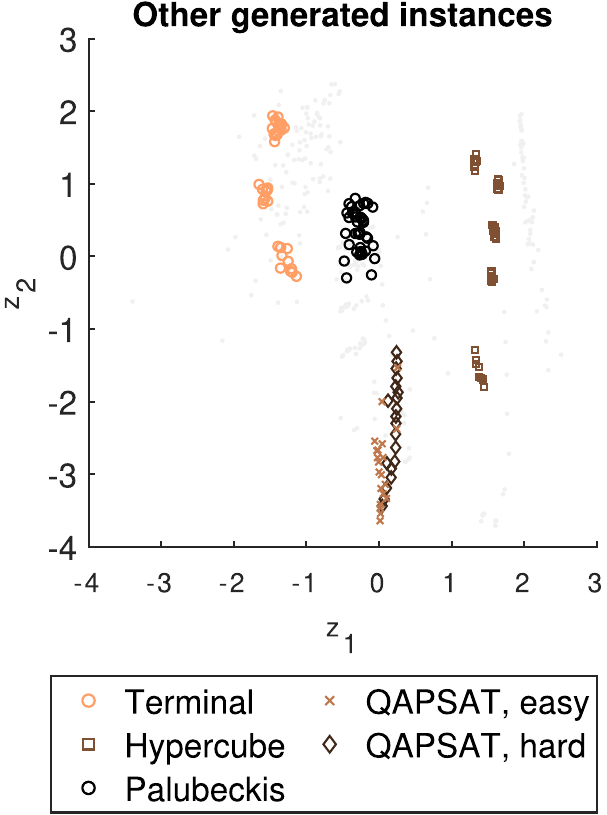}}%
    \subfloat{\includegraphics[trim=0.0cm -0.8cm 0.0cm 0.0cm, clip=true, width=0.24\textwidth]{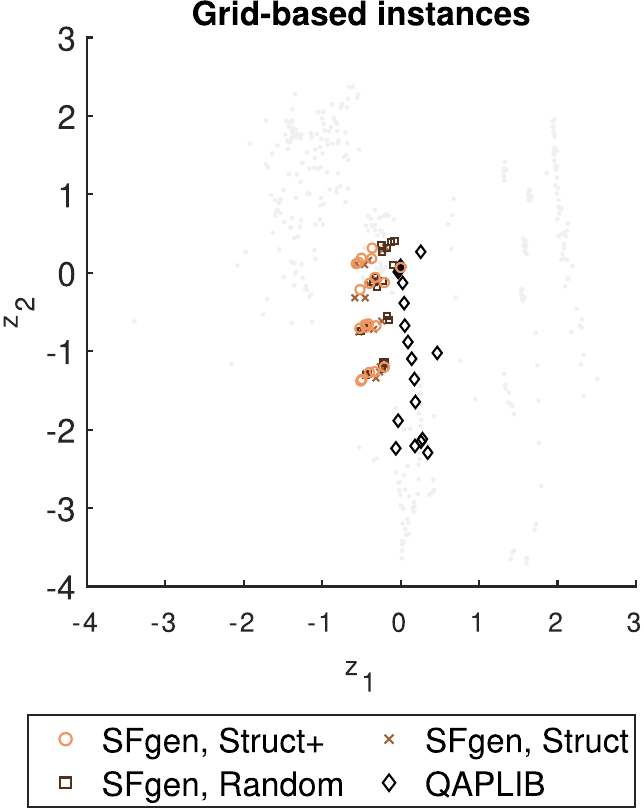}}%
    \subfloat{\includegraphics[trim=0.0cm -0.8cm 0.0cm 0.0cm, clip=true, width=0.24\textwidth]{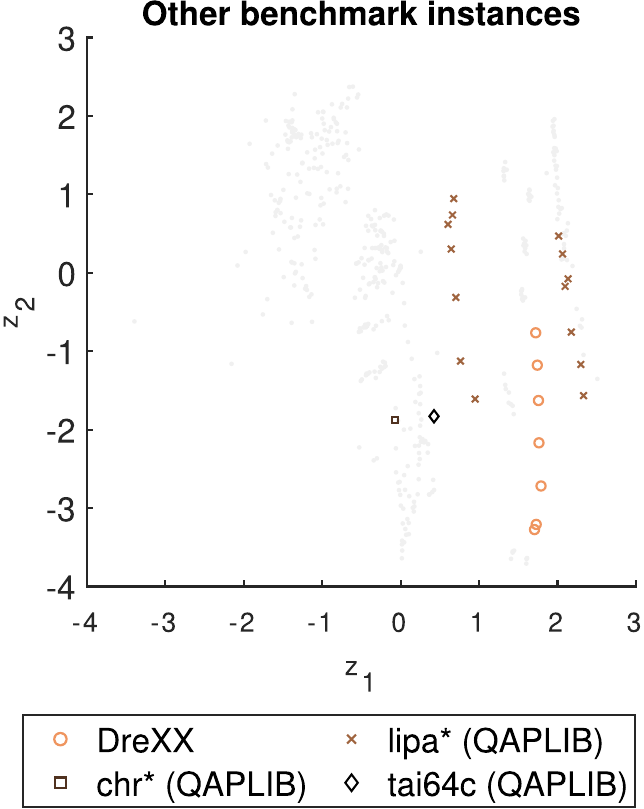}}%
    \caption{Position of instance classes in the initial QAP instance space. Taigen and SFgen refer to the instance generators proposed by \cite{DreznerTaillard2005} and \cite{Stutzle2004Instances} respectively, as described in Section~\ref{sec:instances}.}
    \label{fig:initial-classes-detail}
\end{figure}

\begin{figure}[t]
    \centering
    \includegraphics[trim=0.0cm -0.2cm 0.0cm 0.0cm, clip=true, width=0.5\textwidth]{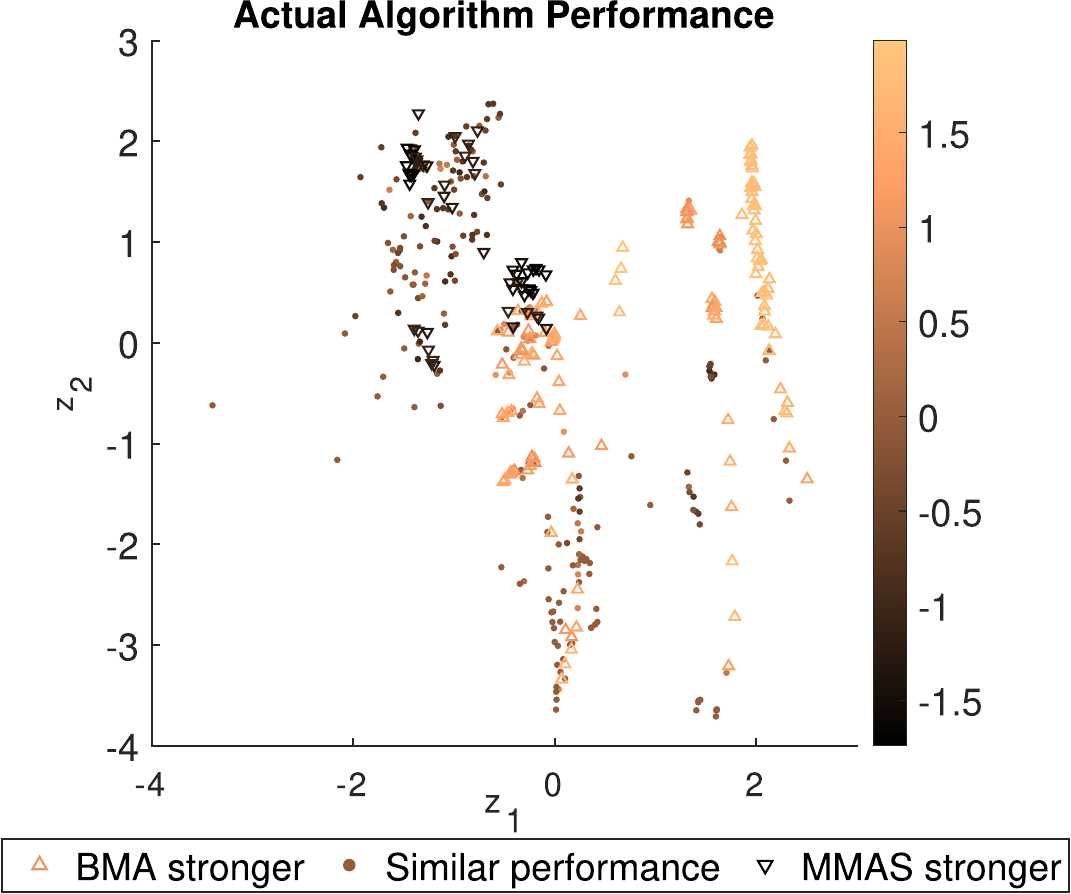}
    \caption{Comparison in performance of the algorithms as defined in \eqref{eq:alg-crit}. Larger values and lighter colours indicates that BMA is better compared with MMAS. Symbols indicate where each algorithm is clearly superior.}
    \label{fig:initial-algorithm}
\end{figure}

\begin{figure}[t]
    \centering
    \subfloat{\includegraphics[trim=0.0cm 0.0cm 0.0cm 0.0cm, clip=true, width=0.24\textwidth]{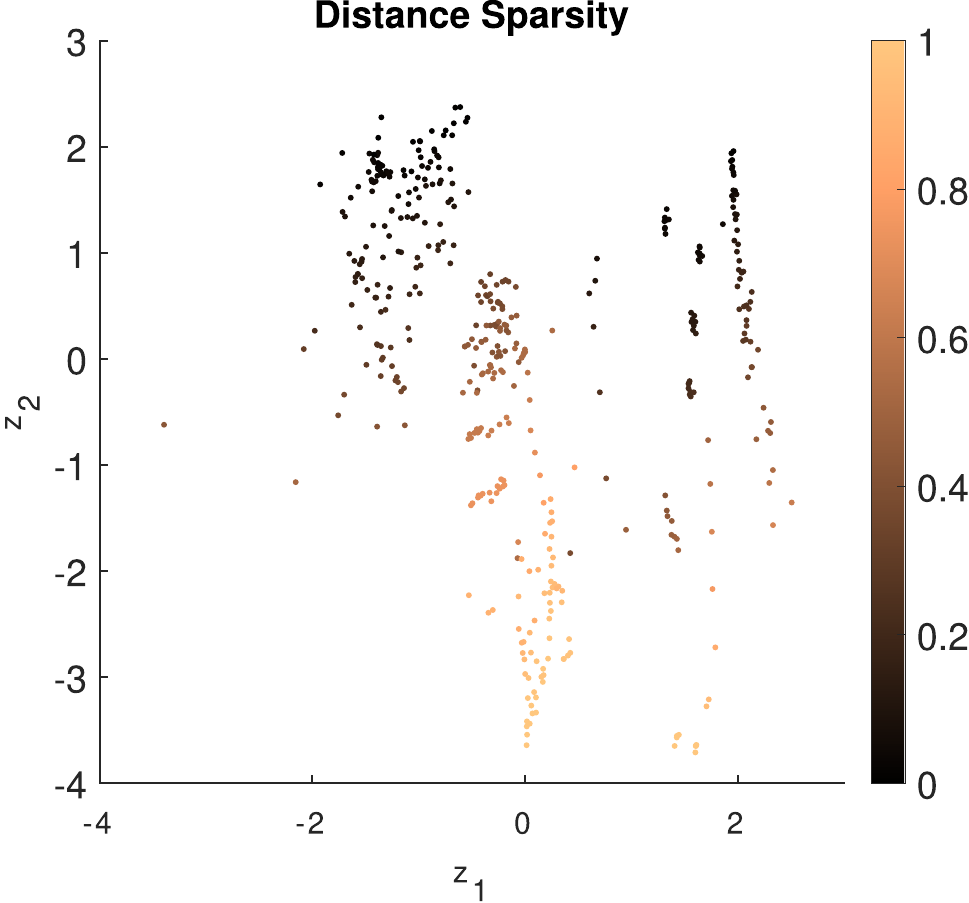}}%
    \subfloat{\includegraphics[trim=0.0cm 0.0cm 0.0cm 0.0cm, clip=true, width=0.24\textwidth]{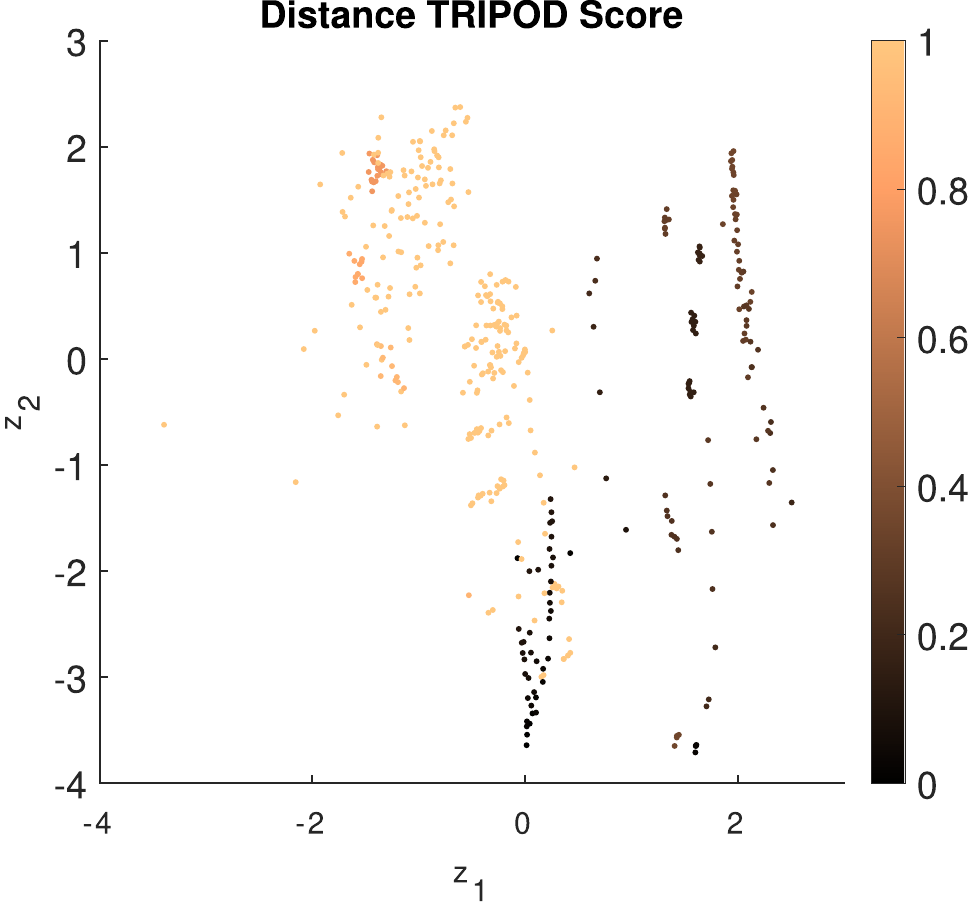}}%
    \subfloat{\includegraphics[trim=0.0cm 0.0cm 0.0cm 0.0cm, clip=true, width=0.24\textwidth]{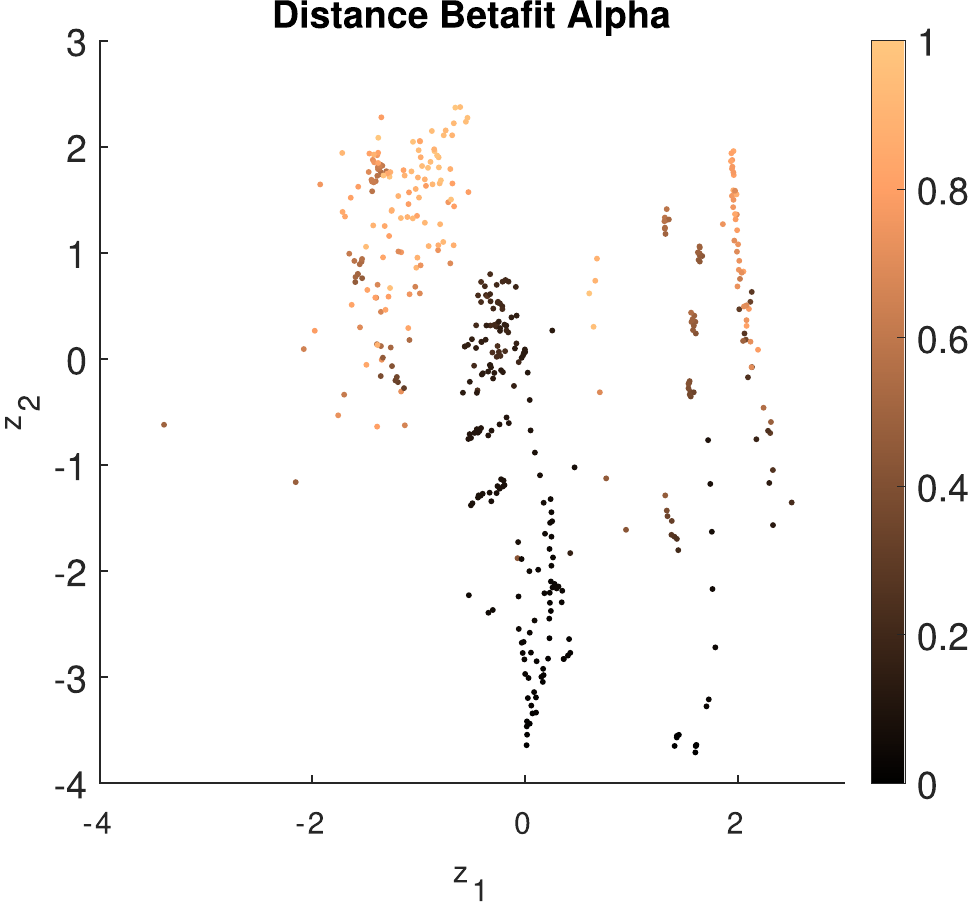}}\\%
    \subfloat{\includegraphics[trim=0.0cm 0.0cm 0.0cm 0.0cm, clip=true, width=0.24\textwidth]{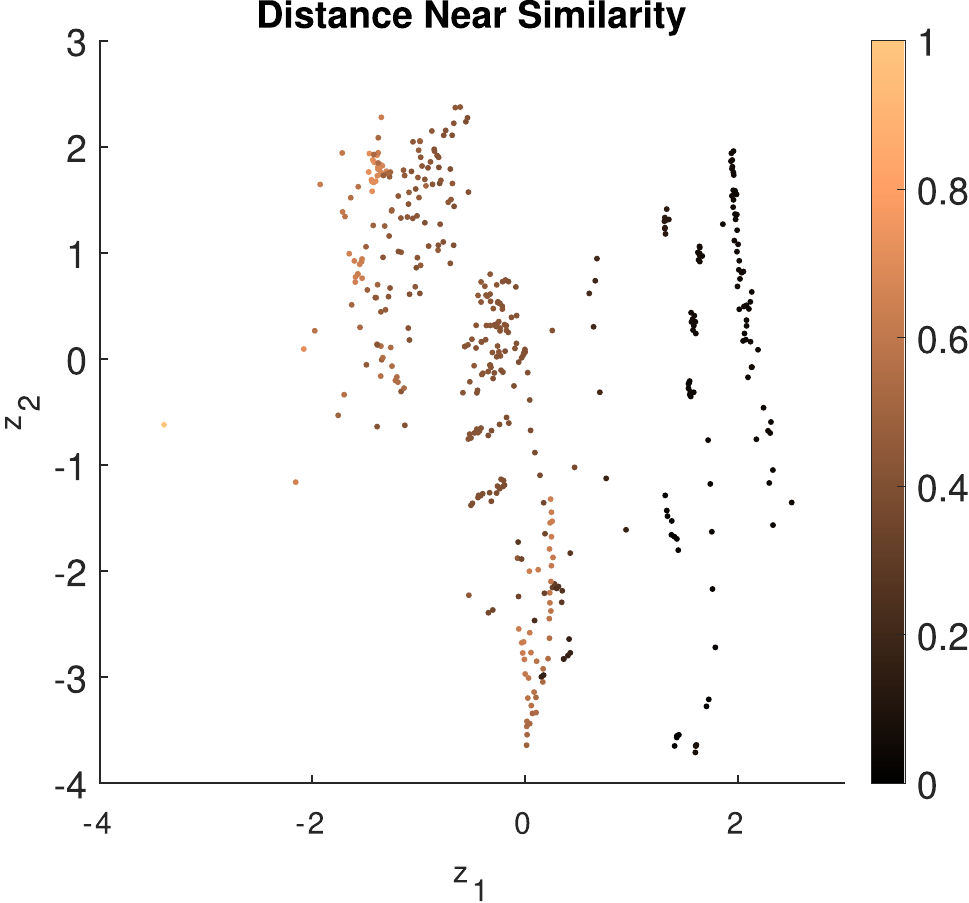}}%
    \subfloat{\includegraphics[trim=0.0cm 0.0cm 0.0cm 0.0cm, clip=true, width=0.24\textwidth]{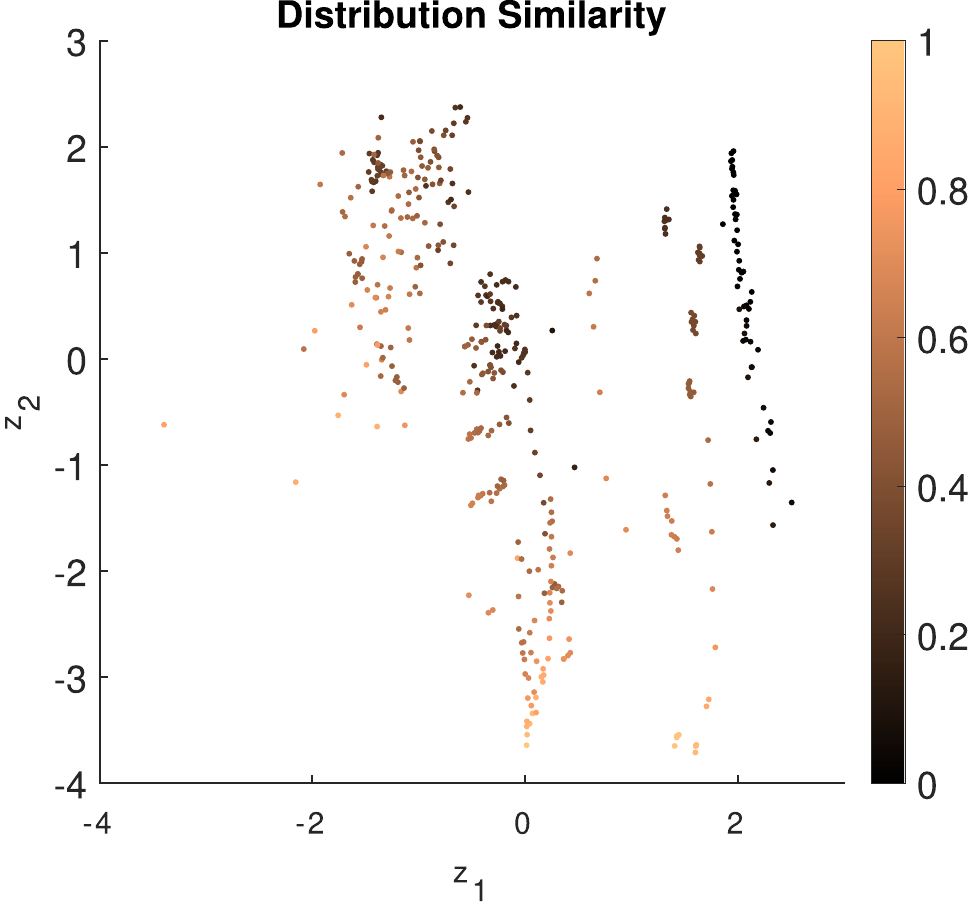}}
    \subfloat{\includegraphics[trim=0.0cm 0.0cm 0.0cm 0.0cm, clip=true, width=0.24\textwidth]{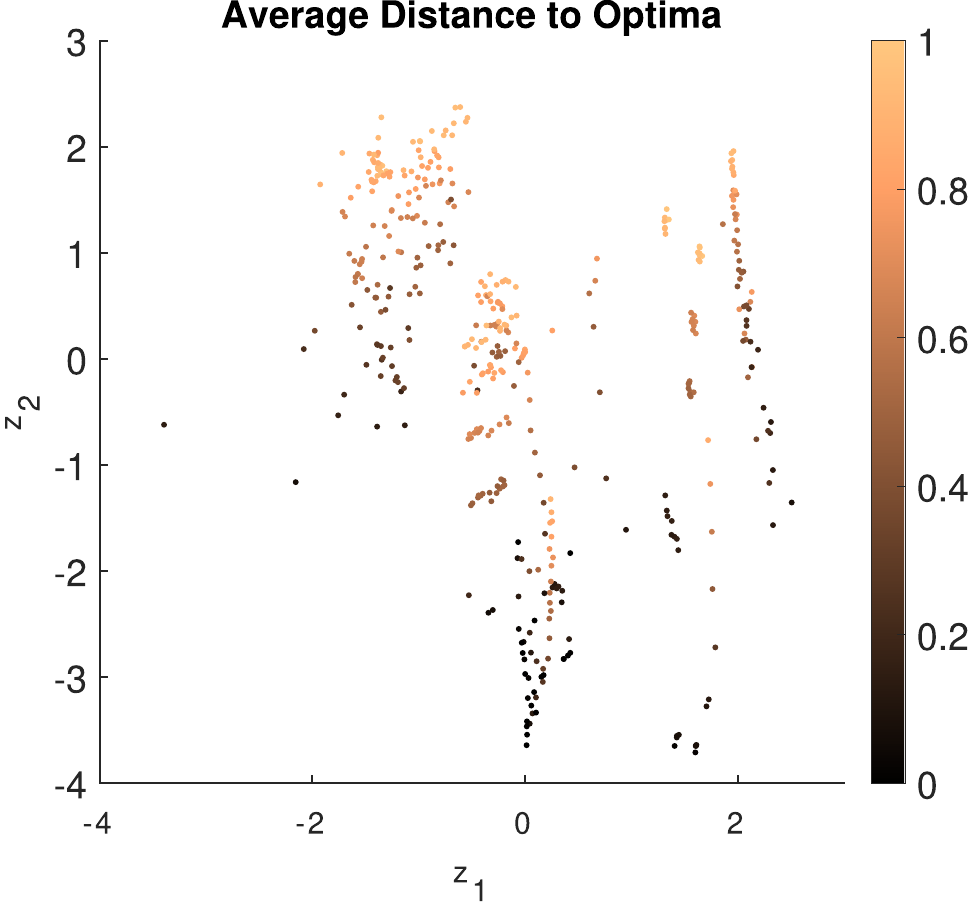}}%
    \caption{Feature distributions.}
    \label{fig:initial-feature}
\end{figure}

The results of this projection are presented in Figures~\ref{fig:initial-classes-macro}-\ref{fig:initial-feature}. Figures \ref{fig:initial-classes-macro} and \ref{fig:initial-classes-detail} indicate the class of each instance. Figure~\ref{fig:initial-qaplib} shows the location of the QAPLIB instances, and the instance space boundaries inferred by CLOISTER. Figure~\ref{fig:initial-algorithm} gives the performance of the algorithms for each instance. Finally, Figure~\ref{fig:initial-feature} shows the normalised feature values for each instance.

In these plots we see that the MMAS algorithm tends to perform best in the upper-left corner of the instance space, corresponding to instances whose distance matrix has a high TRIPOD score and low sparsity. These are characteristics we would expect of distance matrices based on an ordinary distance metric, and so it is no surprise to find the real-like instances based on Euclidean distances in this region. The Terminal and Palubeckis instance classes are also found in this area, and the MMAS algorithm performs particularly well on these instances. 

The grid-based instances with distances based on Manhattan distance on a grid are found towards the centre of the space. Recall that the Distance Sparsity feature is applied to the zero-minimum form of the distance matrix (Definition~\ref{def:reduced-form}), so the Sparsity feature is counting all distances between adjacent locations in a grid-based problem as zero, resulting in a larger value for this feature compared to the Euclidean-distance instances. For these grid-based instances, the BMA algorithm tends to outperform MMAS. The BMA algorithm is also generally superior when applied to the Uniform Random, Hypercube and DreXX instances on the right-hand side of the space. 

On the bottom part of the space the BMA algorithm has an edge on a few instances, but for most instances here the algorithms perform similarly. Most of the QAPLIB instances based on real-life data, as well as the QAPSAT instances of the easy and hard varieties, are found in this region. Since our performance measure is relative, we cannot reliably distinguish between instances where both algorithms trivially discover the optimal solution, and instances where both algorithms struggle to a similar degree.

Notably, the real-like instances based on the generators of \citet{Taillard1995classify} and \citet{Stutzle2004Instances} fall in a quite different part of the instance space to the real-life instances found in QAPLIB. The reason for this is plain enough upon an inspection of the real-life instances; in general they do not conform to the assumption made in the real-like instances of a distance matrix based on Euclidean distances between points in a low-dimensional space, with distances instead being based on such problems as typewriter layout \citep{Burkard1997QAPLIB} and circuit design \citep{EscQAPLIB}. We note that all of the real-life instances in our subset are relatively small, which tends to make them easy for both algorithms.

Comparing Figures~\ref{fig:initial-qaplib} and \ref{fig:initial-algorithm} we see that the region of the instance space where MMAS often outperforms BMA has very little representation in the QAPLIB benchmark library, the tai*b instances being the only exception. This suggests that a compararison between algorithms based on QAPLIB instances alone is not guaranteed to offer a fair assessment of an algorithm's qualities.

This 2-D projection is reasonably successful at separating the existing instances according to the strength of the algorithms at solving them; we assess the predictive power of the projection in more detail in Section~\ref{sec:ispace-selection}. However, there are hints that the instance subset is not representative of the entire potential QAP instance space, and therefore our algorithm predictor may fail when applied to instances with properties not represented in the instance subset. In particular, a notable property of the features selected for the projection in \eqref{eq:initialproj} is the emphasis on the distance matrix. Four of the six features are solely dependent on the distance matrix, with one feature based on both matrices and one fitness landscape feature. None of the features focus on the flow matrix alone. Furthermore, a comparison of Figures~\ref{fig:initial-algorithm} and \ref{fig:initial-feature} shows that the distance-based features do most of the work to separate the regions where each algorithm is particularly strong.


Since the role of the distance and flow matrices in defining a QAP instances are symmetrical, it seems unlikely that instance difficulty is inherently always driven by whichever matrix happens to be identified as the distance matrix by the procedure defined in Section~\ref{sec:identify-dist}. An alternative explanation is that the instance classes composing our instance subset have been designed to include a large variety of distance matrix structures, with the design of the flow matrices being subsidiary. This explanation is indirectly supported by the QAP instance classification scheme proposed by \citet{Taillard1995classify} and used by us in Section~\ref{sec:instances}, which primarily refers to the characteristics of the distance matrix to distinguish between instances. In Section~\ref{sec:more-instances} we will expand our instance subset with a view to exploring the role which the flow matrix may play in determining the properties and difficulty of a QAP instance.

\subsection{Learning Algorithm Footprints and Automated Algorithm Selection} \label{sec:ispace-selection}


\begin{figure}
\centering
\begin{minipage}{.45\textwidth}
  \centering
	    \centering
        \subfloat{\includegraphics[trim=0.0cm 0.0cm 0.0cm 0.0cm, clip=true, width=0.48\textwidth]{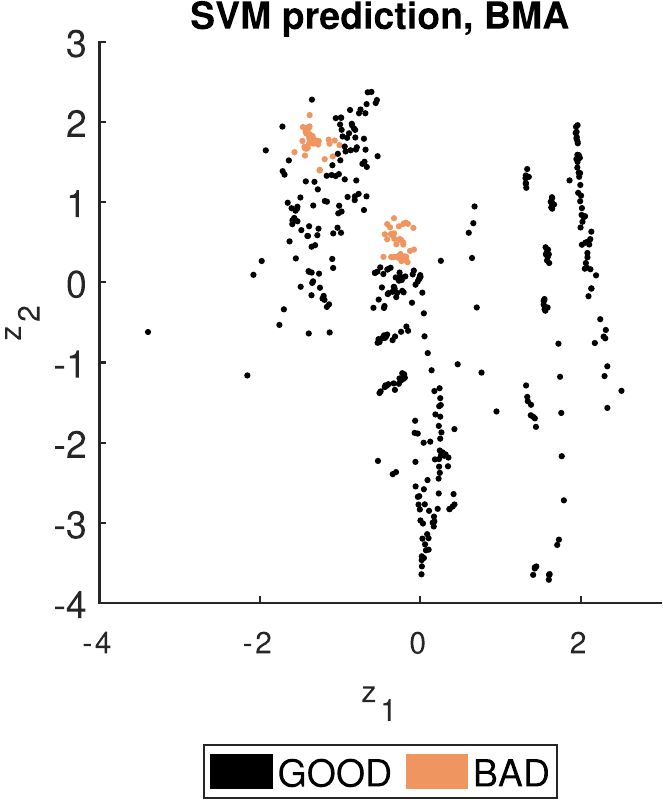}}%
        \subfloat{\includegraphics[trim=0.0cm 0.0cm 0.0cm 0.0cm, clip=true, width=0.48\textwidth]{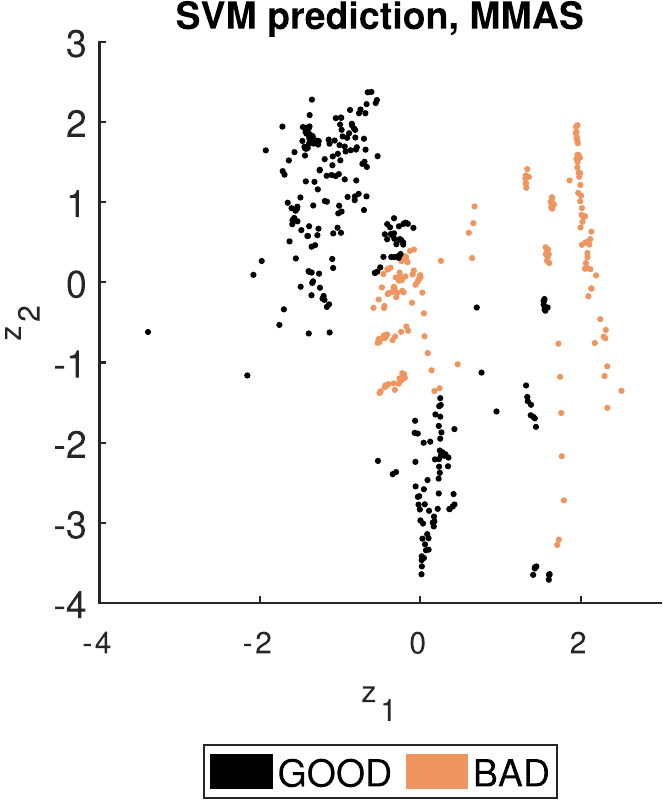}}\\
        \caption{Individual SVM predictions.}
        \label{fig:initial-svm}
\end{minipage}%
\hspace{0.3cm}
\begin{minipage}{.45\textwidth}
  \centering
    \subfloat{\includegraphics[trim=0.0cm 0.0cm 0.0cm 0.0cm, clip=true, width=0.48\textwidth]{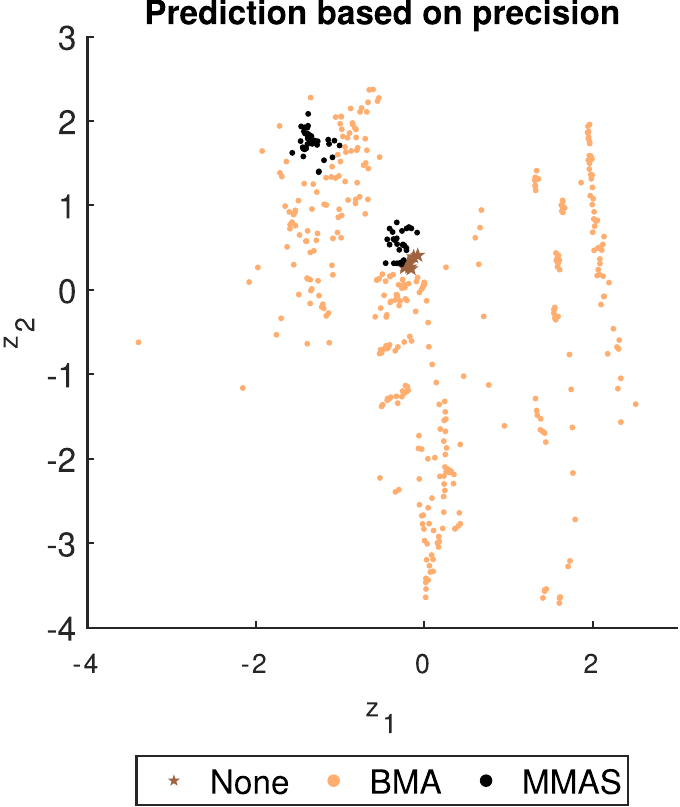}}%
    \subfloat{\includegraphics[trim=0.0cm 0.0cm 0.0cm 0.0cm, clip=true, width=0.48\textwidth]{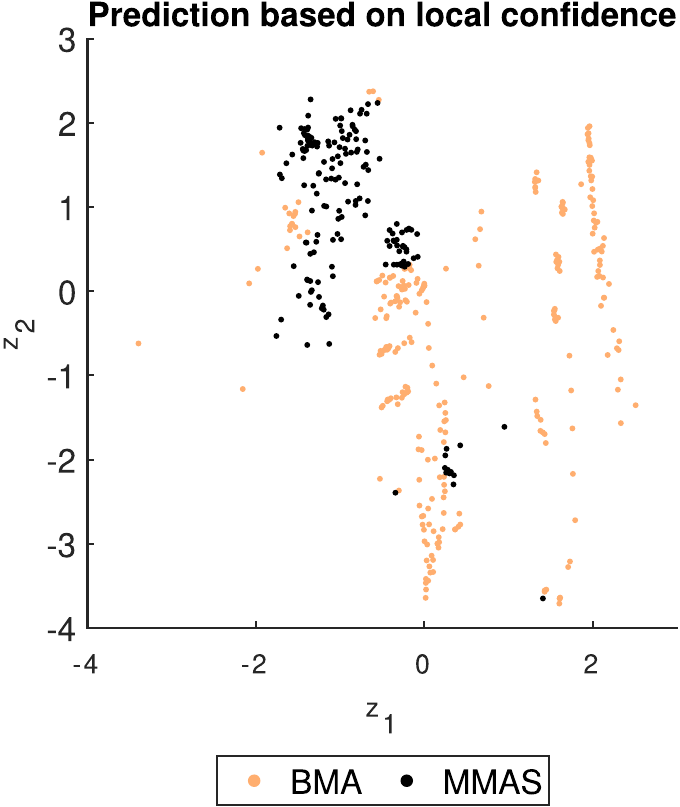}}\\
    \caption{Combined SVM predictions.}
    \label{fig:initial-svm-combined}
\end{minipage}
\end{figure}

The results of the SVMs trained by PYTHIA for the two algorithms are given in Figure~\ref{fig:initial-svm}. The SVM predicting where MMAS will perform well gives results matching our previous intuition about the space in Section~\ref{sec:ispace-properties}; the MMAS algorithm performs particularly well compared to BMA in the upper-left part of the space, and in the lower part of the space both algorithms perform similarly. However, the SVM predicting BMA's performance expects it to perform well over most of the space, including a large part of the upper-left region. This is understandable in the context that the SVM is only interested in distinguishing instances where BMA is actively bad from instances where it has at minimum similar performance to MMAS; in the upper-left region of the space there are many instances (especially in the real-like instances category) where MMAS is moderately stronger than BMA, but not by enough to make BMA a bad choice according to our performance criterion. 

Figure~\ref{fig:initial-svm-combined} shows the results of two algorithm selectors based on the individual algorithm SVMs. The first selector combines the SVMs using the method described by \citet{Smith-Miles2023MATILDA}, using the precision of the SVMs to break ties if both predict their algorithm will perform well. Since the SVM predicting BMA performance has higher precision, this selector chooses BMA for most of the instance space, including much of the upper-left region where MMAS is in fact superior. The precision of an SVM is a global measure of its quality, which allows the BMA SVM to `run up the score' on the large regions in the bottom and right of the space where the BMA algorithm never performs poorly.

The second selector, rather than using the yes/no output of the individual algorithm SVMs, makes its predictions based on the distance of each point from the separating hyperplane produced by the SVMs. The motivation for this idea is that this distance reflects a local measure of confidence which the SVM has in its prediction for a particular instance. These distances are not directly comparable from one SVM to another, so for each SVM we plot the actual performance of the algorithm against the distance from the hyperplane, then find the line of best fit in a least-squares sense. Finally, we use this line to predict the performance of each algorithm on each instance based on its distance from the hyperplane in the corresponding SVM, and select the algorithm which is predicted to perform better. We compare the performance of the two selectors in Table~\ref{table:initial-selector}. The second selector clearly performs better both in terms of confidently selecting an algorithm that is not bad, and actually choosing the stronger of the two algorithms.

\begin{table}[t]
	    \scriptsize
	    \centering
	    \caption{Proportion of instances for which each algorithm selector (based on the SVMs) selects a good algorithm, or the best algorithm (including equal-best).}
	    \label{tab:Instances}
	    \begin{tabular}{m{3.0cm} m{3.5cm} m{3.0cm} }
        \toprule
        Selector & Chooses good algorithm & Chooses best algorithm \\ \hline
        Precision-based & 0.9368 & 0.7065 \\
        Local confidence-based & 0.9819 & 0.8397 \\
        \bottomrule
        \end{tabular} \label{table:initial-selector}
\end{table}

\subsection{Further Analysis of Algorithm Performance} \label{sec:ispace-analysis}

For further insight into the differences in algorithm performance over the instance space, we focus on the Hypercube, Drexx, Terminal and Palubeckis instance classes. All four of these instance classes are constructed such that the small distances and large flows in their respective matrices have a particular structure; in order for a solution to be of high quality it must match the structures so that all of the large flows are assigned to small distances. The Palubeckis and Drexx instances have a small number of structure-matching solutions, whereas the Terminal and Hypercube instances have a large number of structure-matching solutions. However, this does not seem to be predictive of algorithm performance; MMAS typically outperforms BMA on the Palubeckis and Terminal instances, but is outperformed by BMA on the Drexx and larger Hypercube instances. 

In search of a more predictive criterion, we begin with a comparison of the Hypercube and Terminal instances, which both have a large number of structure-matching solutions and penalties of similar magnitude for a solution which is not structure-matching. If we find a small number of adjacent locations and fix the assignment of structure-matching facilities to those locations, we will be left with many structure-matching assignments to choose from for a Terminal instance, but few or only one possible solution for a Hypercube instances. 
If the early stages of the MMAS algorithm result in high pheromone values being attached to these few assignments, then the algorithm has more scope for finding a particularly good solution in the case of a Terminal instance.

The Drexx instances have four structure-matching solutions of equal cost. However, there are many solutions which are nearly structure-matching (for example, moving all facilities one location to the right, wrapping around the right-most facilities) which depending on the value of the random penalties for non-adjacent facilities, may still be a moderately high quality solution. Encountering a component of one such solution and assigning a high pheromone value to those assignments may lead the MMAS algorithm astray. Finally, the Palubeckis instances are designed to have only one strucutre-matching and optimal solution, so if the MMAS algorithm encounters a component of this solution then this will lead it in a desirable direction.


The above discussion leads us to the following conjecture:
\begin{conjecture} \label{hyp:1}
    For a given QAP instance, consider the question `if we find a few assignments which are part of a high-quality solution, would fixing these assignments permanently then solving the rest of the problem be a relatively effective strategy for a solver?' If we can answer this question in the affirmative, we predict that MMAS will perform well relative to BMA. 
\end{conjecture}



Having formed Conjecture~\ref{hyp:1} through analysis of instances with special structure and strong relationships between their distance and flow matrices, it is reasonable to inquire whether the conjecture makes sense in the context of instances where these conditions do not apply. As we have previously noted, BMA generally significantly outperforms MMAS on the grid-based instances. The structure of the DreXX instances' distance matrices is fairly similar to a grid-based instance, and it is reasonable to expect MMAS struggles for much the same reasons on these instances. 

On the other hand, when applied to the real-like instances based on a Euclidean distance metric, the MMAS algorithm tends to have a slight edge but only rarely achieves a substantial performance improvement over BMA. This also has a plausible explanation in terms of our conjecture. Both of the real-like instance generators define their locations in clusters. If an instance happens to contain clusters and groups of facilities that match up especially well together, and the MMAS algorithm can identify these assignments, then we can expect it to perform well. However, since for these instances the distance and flow matrix generation procedures are independent, not all real-like instances will contain such a `jackpot' outcome.

Another reasonable question is whether we can measure the property which makes an instance well- or ill-suited for MMAS with an automated feature, for example with some sort of fitness landscape measurement. Unfortunately, this seems difficult to accomplish for two reasons. First, the obvious method for determining how the pheromones of the MMAS algorithm will develop is to directly apply MMAS to the instance in question; in the context of the Algorithm Selection Problem, where we desire to cheaply determine the most appropriate algorithm, this is an unpleasantly expensive prospect. Second, even if we can cheaply guess which assignments the MMAS algorithm will assign high pheromone values to, it is not obvious how to determine whether these pheromones will lead to a good outcome without \emph{a priori} knowledge about the high-quality solutions of the instance. A thorough fitness landscape analysis might uncover the required information, but again such a procedure is likely to take too long to be practical in an algorithm selection context.

Even though we do not have a practical automated method for algorithm selection based on Conjecture~\ref{hyp:1}, we can still test the conjecture by designing new instances where we predict one algorithm or the other will perform strongly, and then running the algorithms to confirm or refute our prediction. Constructing novel instances with these ideas in mind, and using them to expand our instance subset, can also improve our understanding of the QAP in general. We pursue this idea in Section~\ref{sec:flow-mmas}.

\section{Expanding the Instance Space}\label{sec:more-instances}

\subsection{Hybrid Instance Classes} \label{sec:hybrid}

The first method we consider for expanding the instance space is combining distance matrices produced by one existing generator with flow matrices produced by a different generator. The distance and flow matrix generators we include in this method are listed in Tables~\ref{table:hybrid-dist-gen} and \ref{table:hybrid-flow-gen} respectively. In these tables SFgen refers to the instance generator proposed by \citet{Stutzle2004Instances}. All generator parameters not stated are as in the previously generated instances.

Table \ref{table:recombinator} indicates the combinations of distance and flow generators we include in this method, and lists the maximum and minimum algorithm performance criterion. We omit pairs of generators which generate matrices of incompatible sizes, or which would result in duplication or near-duplication of an existing instance class. 

\begin{table}[t]
\centering
\begin{minipage}{.45\textwidth}
  \scriptsize
    \centering
    \caption{Distance generators used in hybrid instances. }
    \rowcolors{2}{gray!15}{white} 
    \begin{tabular}{l c l }
      \toprule
      Generator & Size & Parameters  \\ \midrule
        Uniform Random & $81$ & -- \\
        Hypercube & $81$ & $d_0 = 20, \ell=3, k=4$ \\
        Palubeckis & $81$ & -- \\
        Terminal & $81$ & $4$ levels, $3$ branches each \\
        SFgen, Real-like & $80$ & -- \\
        SFgen, Manhattan & $80$ & $8 \times 10$ grid \\
        Drexx & $80$ & $8 \times 10$ grid \\
       \bottomrule
    \end{tabular}
        \label{table:hybrid-dist-gen}
\end{minipage}%
\hspace{0.3cm}
\begin{minipage}{.45\textwidth}
  \scriptsize
    \centering
    \caption{Flow generators used in hybrid instances.}
    \rowcolors{2}{gray!15}{white} 
    \begin{tabular}{l c l}
      \toprule
      Generator & Size & Parameters  \\ \midrule
        Uniform Random & $80$ or $81$ & -- \\
        Hypercube & $81$ & $f_0 = 20, \ell=3, k=4$ \\
        Palubeckis & $80$ or $81$ & -- \\
        Terminal & $81$ & $4$ levels, $3$ branches each \\
        SFgen, Structured Plus & $80$ or $81$ \\
       \bottomrule
    \end{tabular}
        \label{table:hybrid-flow-gen}
\end{minipage}
\end{table}

\begin{table}[t] \footnotesize
        \centering
        \caption{Summary of hybrid instance classes. $\checkmark$ indicates a combination of flow and distance matrix used to produce new hybrid instances. D indicates a combination which would duplicate an existing instance class. I indicates a combination which would be incompatible in terms of instance size.} 
        \rowcolors{2}{gray!15}{white} 
		\begin{tabular}{ c c c c c c }
			\toprule
			  $\downarrow$ Distances $\backslash$ Flows $\rightarrow$ & Uniform Random & Hypercube & Palubeckis & Terminal & Structured Plus \\
                \midrule
                Uniform Random & D & $\checkmark$ & $\checkmark$ & $\checkmark$ & $\checkmark$ \\ 
                Hypercube & $\checkmark$ & D & $\checkmark$ & $\checkmark$ & $\checkmark$ \\ 
                Palubeckis & $\checkmark$ & $\checkmark$ & D & $\checkmark$ & $\checkmark$ \\ 
                Terminal & $\checkmark$ & $\checkmark$ & $\checkmark$ & D & $\checkmark$ \\ 
                Real-like & D & $\checkmark$ & $\checkmark$ & $\checkmark$ & D \\ 
                Manhattan & D & I & $\checkmark$ & I & D \\ 
                Drexx & $\checkmark$ & I & $\checkmark$ & I & $\checkmark$  \\ \bottomrule
			\end{tabular}
		\label{table:recombinator}
\end{table}

\begin{figure}[t]
    \centering
    \subfloat{\includegraphics[trim=0.0cm 0.0cm 0.0cm 0.0cm, clip=true, width=0.28\textwidth]{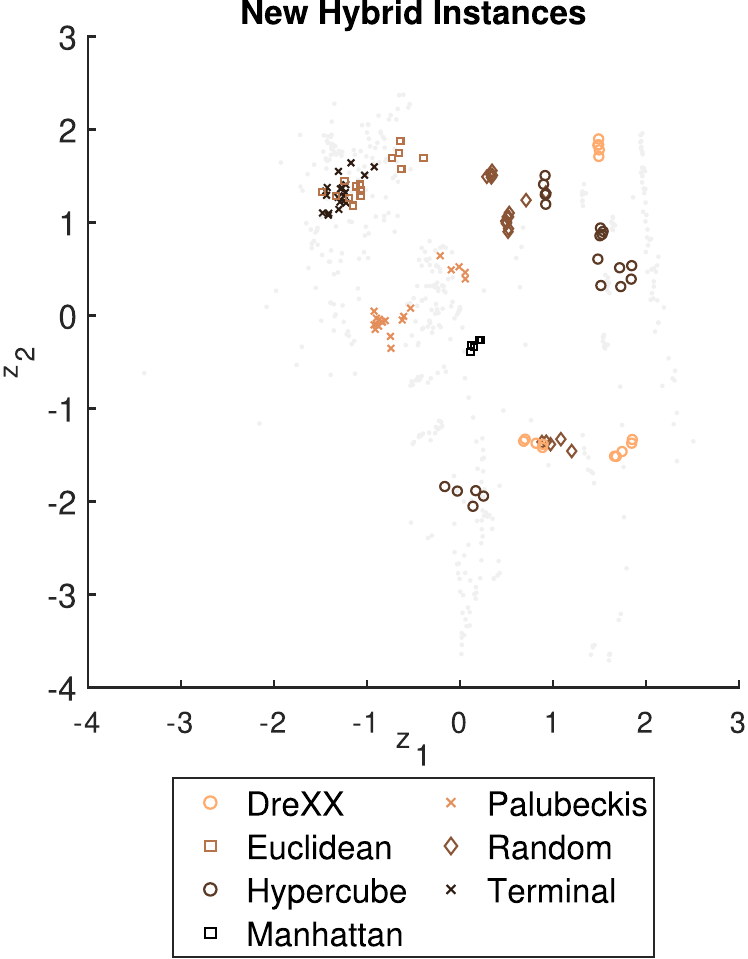}}%
    \hspace{5mm}
    \subfloat{\includegraphics[trim=0.0cm -2.2cm 0.0cm 0.0cm, clip=true, width=0.38\textwidth]{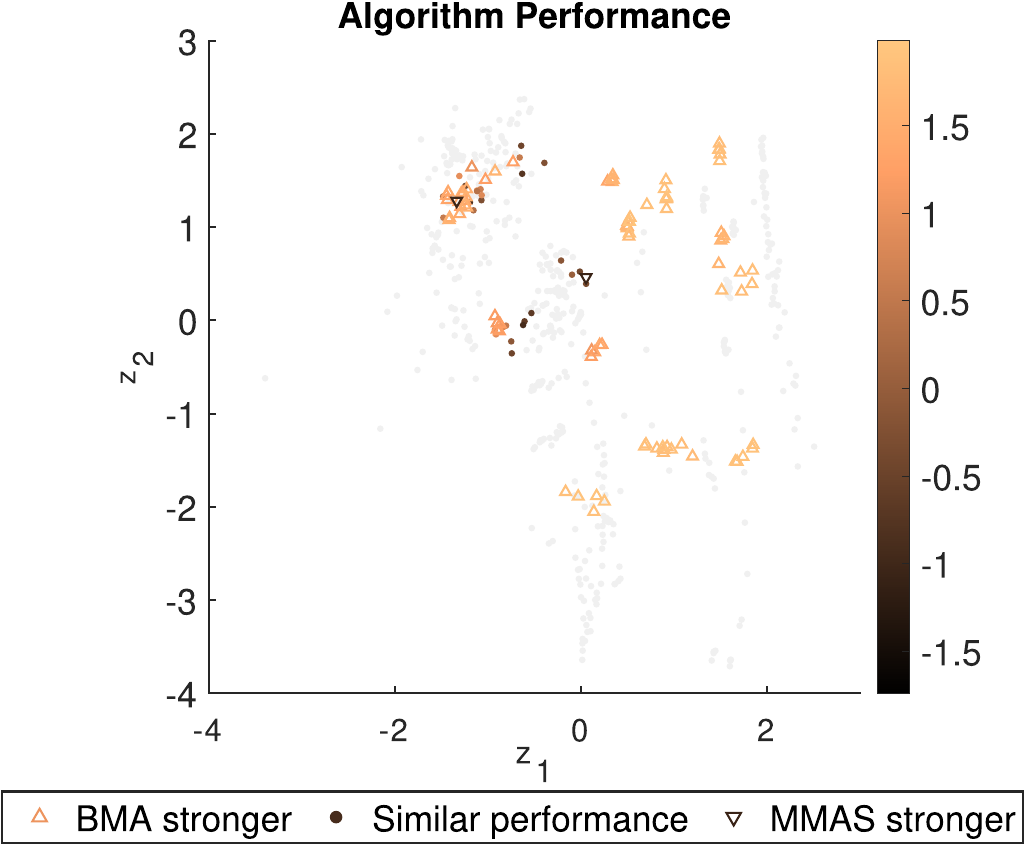}}
    \caption{Position of hybrid instances in the initial instance space, and comparative algorithm performance on these new instances. Instances are labelled according to the distance matrix generator used to produce them.}
    \label{fig:hybrid-initial-instance}
\end{figure}

The location of the new hybrid instances in the instance space is shown in Figure~\ref{fig:hybrid-initial-instance}. Using distance generators similar to those already in the instance space, this method has produced some instances in the upper-middle region of the space, which was previously poorly explored. Note that the distance matrix selecting procedure defined in Section~\ref{sec:identify-dist} is applied to these instances as usual, so for any given instance the matrix produced by the `distance matrix generator' is not necessarily identified as the distance matrix for the purpose of feature measurement.

Figure~\ref{fig:hybrid-initial-instance} also displays the performance of the algorithms on these new instances; we see that BMA is the stronger algorithm for almost all of these instances, including most instances in the top-left region where, in the initial instance subset, MMAS was stronger. Therefore, contrary to the conclusion we might draw from the features chosen for the projection which largely ignored the flow matrix, we see that the flow matrix can signficantly affect problem difficulty. However, it seems that creating new instances in the mannner proposed in this section generally favours BMA over MMAS.

\subsection{Designing Flow Generators to Favour MMAS}\label{sec:flow-mmas}

In the context of Conjecture~\ref{hyp:1} the conclusion of the previous section perhaps ought not to surprise us. An instance without some kind of relationship between the structure of the distance and flow matrices is unlikely to contain the especially good assignments which we conjectured to be desirable for MMAS. In this section we take distance matrices from instance classes where BMA is favoured, and apply Conjecture~\ref{hyp:1} to design flow matrix generators, with the goal that MMAS should be the superior algorithm for the resulting instances.



We consider six instance subclasses in this section, combining two existing distance generators with three novel flow generators. The distance matrices are based on either the distances of the Drexx instance class (using an $8 
\times 10$ grid, so $n=80$) or the Hypercube instance class (with $k=4$, $\ell = 3$, $d_0 = 20$, so $n=81$). Both of these distance generators are based on a uniform grid of locations. The original instances of these classes had flow matrices which linked the facilities in a single monolithic group, which had to exactly match the structure of the distance matrix in a high-quality solution. In contrast, the flow generators proposed in this section divide the facilities into smaller clusters with large flows only within each cluster; hence we denote the class as `flow-cluster' instances. The motivating idea is to emulate the structure of the Terminal instance class, in which the assignment of a single cluster of locations has a moderate, but not overwhelming, influence on the rest of the solution.

As earlier, $\disunif{a}{b}$ represents an integer in $\left\{ a, \dots, b \right\}$ selected using a uniform random distribution, generated independently for each flow it is used to define. The three flow generators are defined as follows:

\paragraph*{Triangle flows}
This generator has two parameters, base flow $f_0 = 20$ and number of triangles $t$ which we vary. Begin by creating a grid of facilities which matches the structure of the distance matrix i.e. a $8 \times 10$ rectangle or a $3 \times 3 \times 3 \times 3$ hypercube. Also create an initial flow matrix of all zeros. Then for each of the $t$ desired triangles, we randomly select a facility $i \in N$, randomly select two of its neighbours $j, k \in N$, and assign each of the flows $f_{ij}$, $f_{jk}$ and $f_{ki}$ a random value from $\disunif{f_0+1}{ f_0 + \frac{f_0}{2}}$. Note that if a particular flow between two facilities is selected more than once, its value is overwritten not increased.

\paragraph*{Tree flows}
This generator has two parameters, base flow $f_0 = 20$ and number of groups $s$ which we vary. The generator is based on the idea of $s$ groups of facilities $S_\sigma$ for $\sigma = 1,\dots,s$, which may be larger or smaller depending on their weighting $w_\sigma$. The facilities within each group are connected by links with no cycles, hence the `tree' description; all facilities in a group have some flow between them, but directly linked facilities share stronger flows.

As in the triangle flow generator, we create a grid of facilities which matches the structure of the distance matrix and an initial flow matrix of all zeros. We initialise each $S_\sigma$ with a randomly selected facility as its root and a weighting $w_\sigma$ from $\disunif{1}{5}$. We then randomly select a group $S_\sigma$ with probability proportional to its weighting, randomly select a facility already in $S_\sigma$, and link that facility to a randomly selected neighbour which is added to $S_\sigma$. This process is repeated until all facilities are in a group. Facilities which are directly linked receive flows chosen from $\disunif{f_0+1}{ f_0 + \frac{f_0}{2}}$ in both directions between them; facilities which are in the same group, but not directly linked, instead receive flows chosen from $\disunif{1}{\frac{f_0}{2}}$. All other flows are zero.

\paragraph*{Square flows}
This generator has one parameter, the number of squares $s$. We independently generate $s$ $4 \times 4$ matrices $F_i$ with the form 
\begin{equation*}
F_i = \begin{pmatrix} 0 & 10^4 & \star & 10^4 \\
10^4 & 0 & 10^4 & \star \\
\star & 10^4 & 0 & 10^4 \\
10^4 & \star & 10^4 & 0 \\
\end{pmatrix}
\end{equation*}
where each $\star$ is replaced by an independently chosen random value from $\disunif{1}{100}$. The flow matrix is then created by merging the $F_i$ matrices in a block-diagonal form, with any additional rows and columns required to form an $n \times n$ matrix filled with zeros.



\begin{figure}[t]
    \centering

    \subfloat{\includegraphics[trim=0.0cm 0.0cm 0.0cm 0.0cm, clip=true, width=0.33\textwidth]{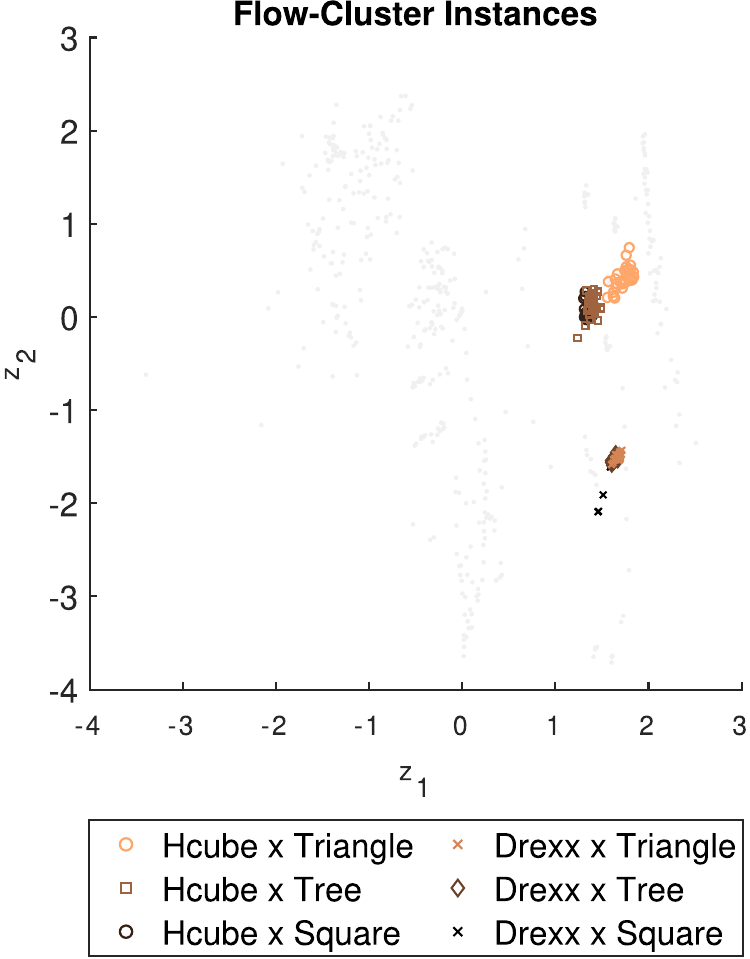}}%
    \caption{Position of flow-cluster instances in the initial instance space.}
    \label{fig:ispace-flowcluster}
\end{figure}

\begin{figure}[t]
    \centering
    \subfloat{\includegraphics[trim=0.0cm 0.0cm 0.0cm 0.0cm, clip=true, width=0.33\textwidth]{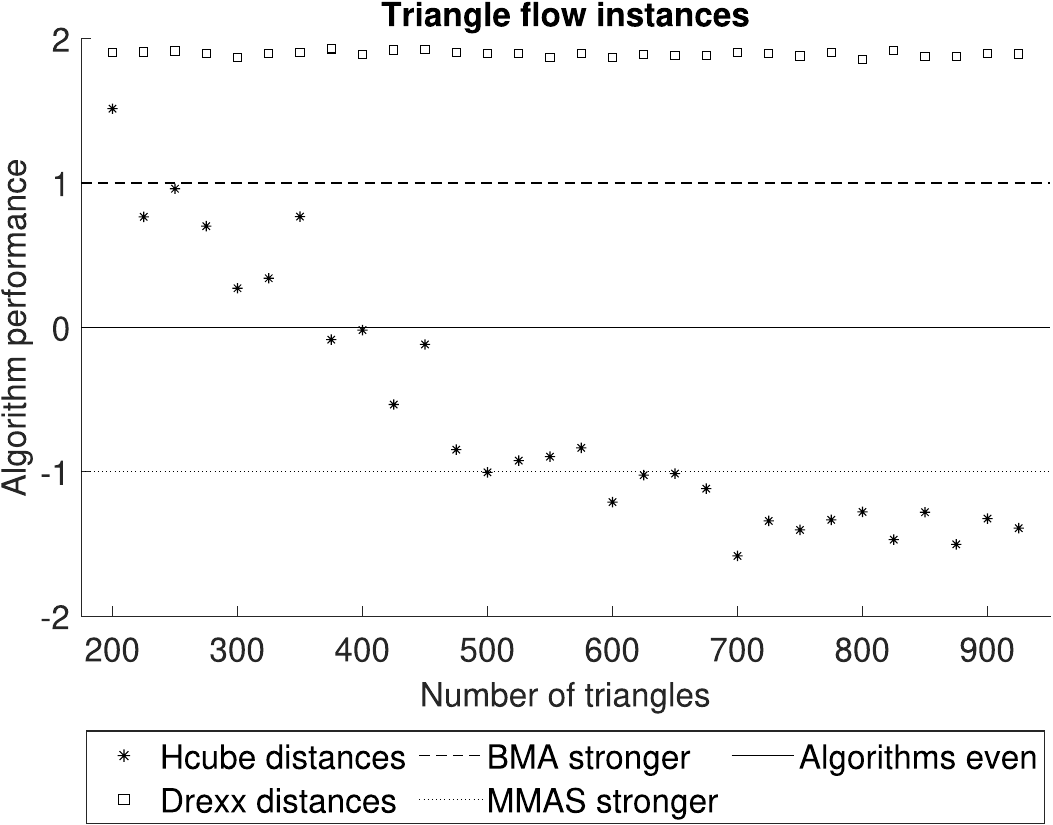}}
    \subfloat{\includegraphics[trim=0.0cm 0.0cm 0.0cm 0.0cm, clip=true, width=0.33\textwidth]{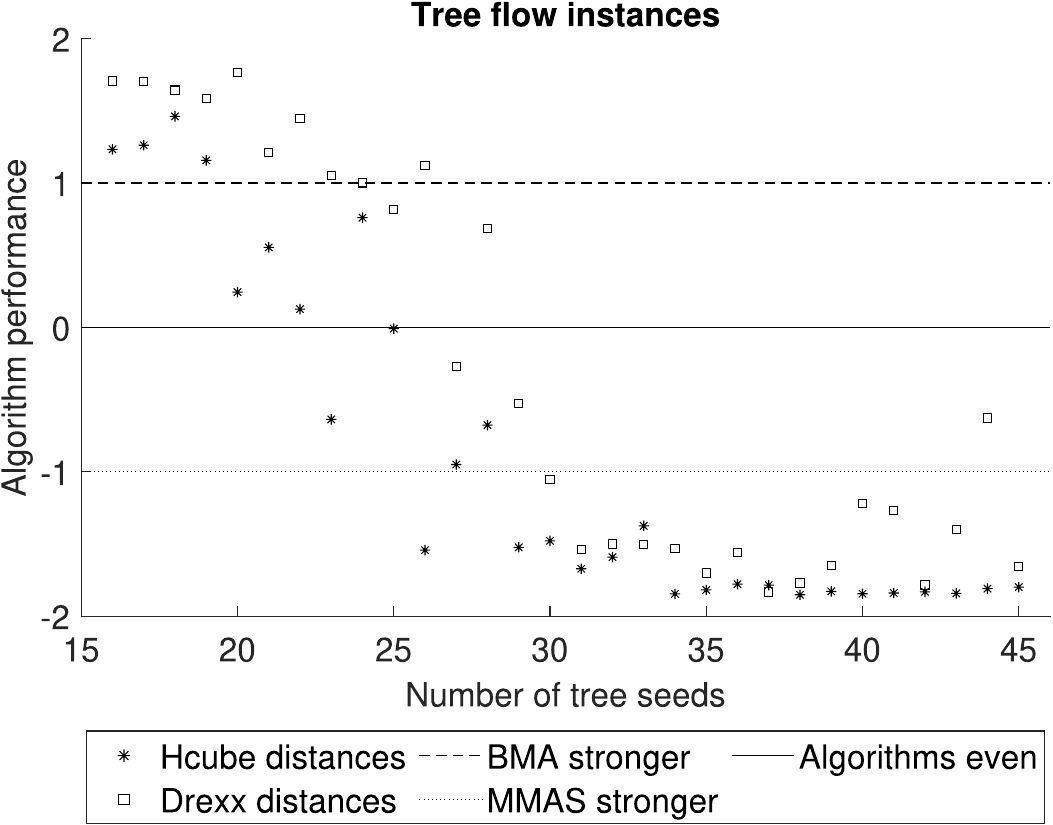}}
    \subfloat{\includegraphics[trim=0.0cm 0.0cm 0.0cm 0.0cm, clip=true, width=0.33\textwidth]{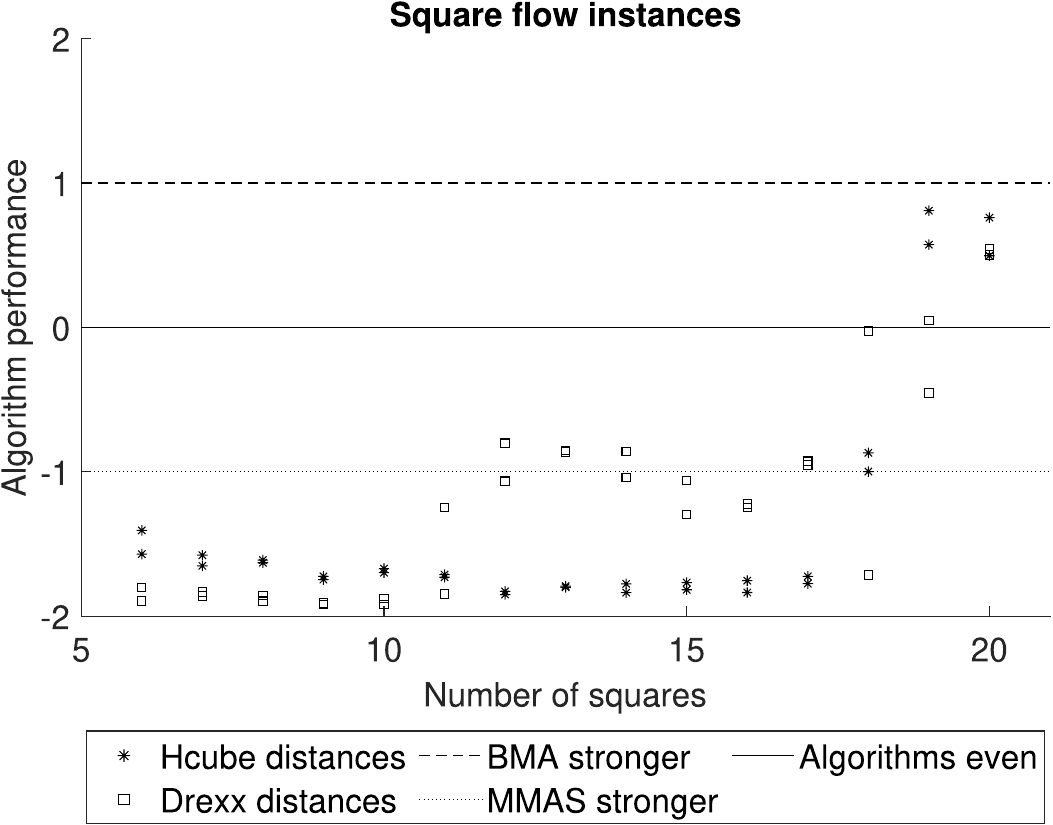}}
    \caption{Comparative algorithm performance for flow-cluster instances, plotted against each generator's parameter.}
    \label{fig:ispace-fcperf}
\end{figure}

The position of the flow-cluster instances in the instance space is shown in Figure~\ref{fig:ispace-flowcluster}. As expected, these instances occupy a similar position to the original Hypercube and Drexx instances, on the right-hand side of the space. In Figure~\ref{fig:ispace-fcperf} we plot the algorithm performance against the parameter setting for each flow matrix generator. 

With suitable settings of the generator parameters, five of the six instance subclasses considered in this section reliably produce instances which favour MMAS over BMA. The exception is the instance subclass with DreXX distances and triangle flows, which always strongly favours BMA. The triangle flow generator puts a particular emphasis on locations separated by two edges. In the DreXX distance generator all distances between non-adjacent locations are generated from the same distribution; however, the Hypercube distance generator assigns a shorter distance to locations separated by exactly two edges in the grid. 

The other five instance subclasses follow a pattern where one extreme of the flow generator parameter produces instances which favour BMA, and moving towards the other extreme produces instances which favour MMAS. For the tree and square flow instances the trend is consistent for both the DreXX and Hypercube instances, with the DreXX instances often being a little more favourable for BMA. The easiest dependence on generator parameter to understand is for the square-based flow instances: when the number of squares becomes large so that all or nearly all locations must be occupied by a facility belonging to one of the squares, if MMAS fixes a square `out of alignment' then there will not be enough space to fit in the remaining squares while keeping all of the large flows only between adjacent locations. 

The Triangle and Tree flow generators have less rigid structure, so it is more difficult to confidently determine the cause of the relationship between these generators' respective parameters and the relative instance difficulty. For the Triangle flow instances with Hypercube distances, it seems plausible that when the number of cycles is low that the resulting problem has so much flexibility in finding potentially-good solutions that fixing a few of the cycles in one place is too restrictive. On the other hand, when the number of tree seeds in the tree-based flow generator is too low, the resulting clusters of linked facilities may be too large to give MMAS an advantage.

\subsection{New Instance Space}\label{sec:new-isa}

We now perform a second instance space analysis on the 443 original instances plus the 115 hybrid and 180 flow-cluster instances defined above. For this analysis we use the same parameters as the first, except that we now set the number of selected features is now set to $7$. Based on visual inspection, the projections produced using a greater or fewer number of features were less effective at distinguishing instances where MMAS outperformed BMA. The projection which produces the new $2$D instance space is:

\begin{equation}
		\begin{bmatrix} 
		Z_1 \\ 
		Z_2
		\end{bmatrix}
		= 
		\begin{bmatrix}[r]
			-0.0234 &   -0.2319 \\
            -0.2466 &   -0.0596 \\
            -0.1477 &   0.2326 \\
            0.3352 &    0.1186 \\
            -0.4022 &   0.2386 \\
            -0.2025 &   -0.1974 \\
            0.1967 &   0.1066 \\
		\end{bmatrix}^\Tr
		\begin{bmatrix}
		\text{Distance Normalised Mean} \\
		\text{Flow Normalised Mean} \\
		\text{Distance Sparsity} \\
		\text{Flow Dominance} \\ %
		\text{Distance Skewness} \\
		\text{Gilmore Lawler Bound} \\
        \text{Escape Probability} \\
		\end{bmatrix}
		\label{eq:laterproj}
\end{equation}

\begin{figure}[t]
    \centering
    \subfloat{\includegraphics[trim=0.0cm 0.0cm 0.0cm 0.0cm, clip=true, width=0.33\textwidth]{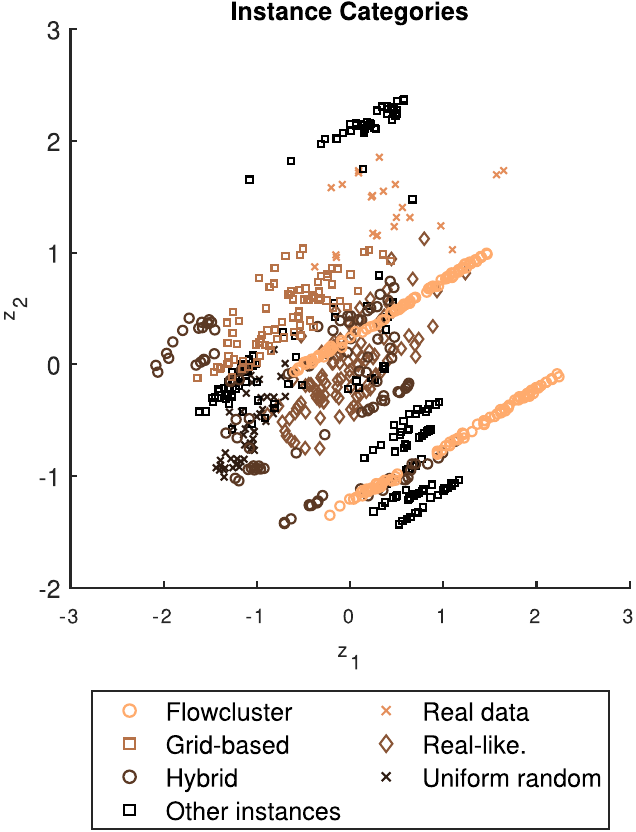}}%
    \subfloat{\includegraphics[trim=0.0cm 0.0cm 0.0cm 0.0cm, clip=true, width=0.33\textwidth]{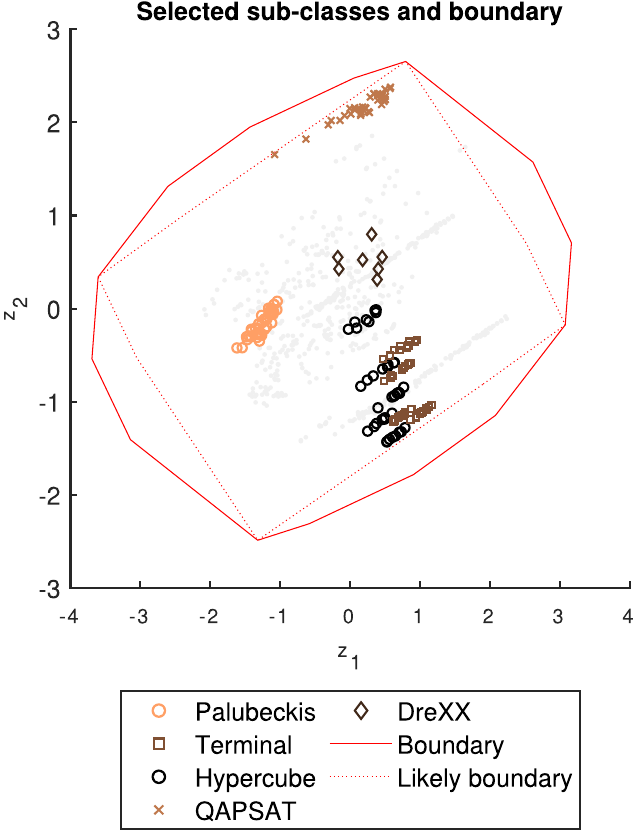}}%
    \subfloat{\includegraphics[trim=0.0cm -0.7cm 0.0cm 0.0cm, clip=true, width=0.33\textwidth]{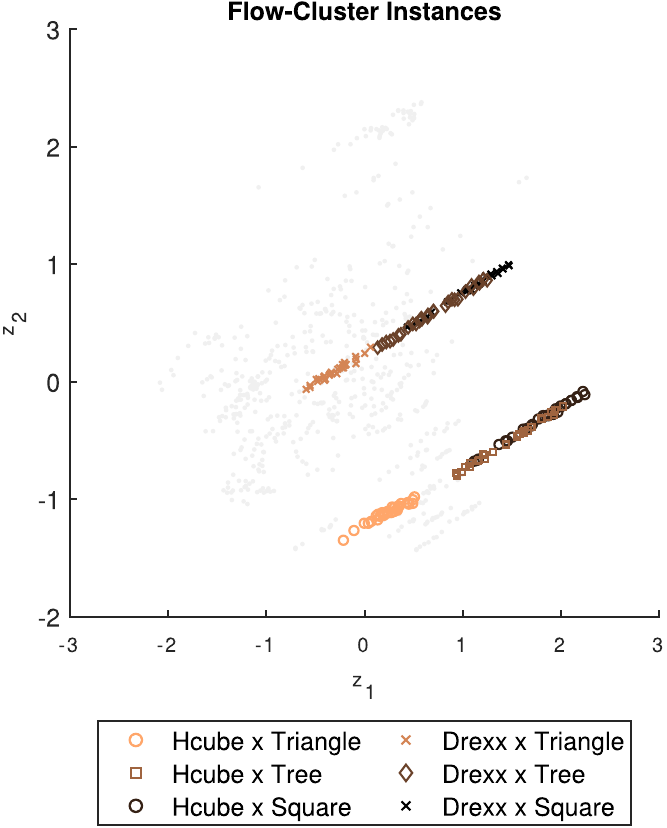}}%
    \caption{Position of instance classes in the second QAP instance space.}
    \label{fig:ext-classes}
\end{figure}

\begin{figure}[t]
    \centering
    \includegraphics[trim=0.0cm 0.0cm 0.0cm 0.0cm, clip=true, width=0.4\textwidth]{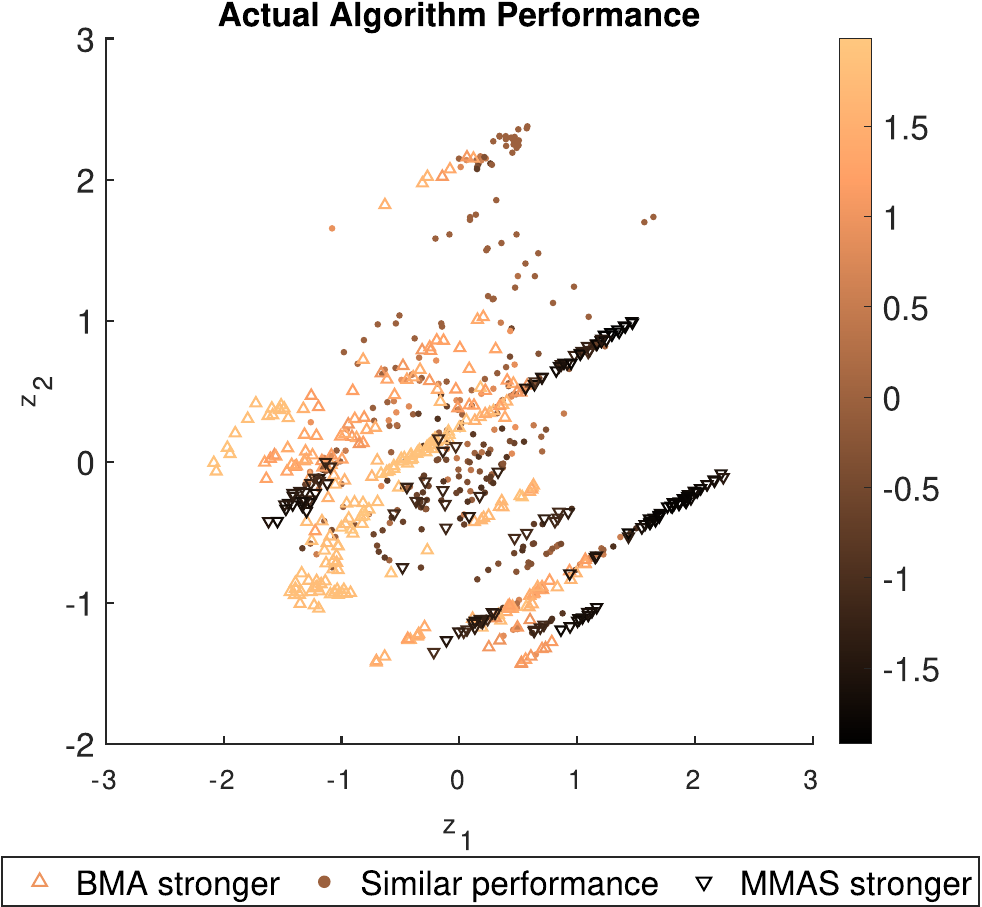}
    \caption{Comparison in performance of the algorithms as defined in \eqref{eq:alg-crit}. Larger values and lighter colours indicates that BMA is better compared with MMAS. Symbols indicate where each algorithm is clearly superior.}
    \label{fig:ext-algorithm}
\end{figure}

\begin{figure}[t]
    \centering
    \includegraphics[trim=0.0cm 0.0cm 0.0cm 0.0cm, clip=true, width=0.24\textwidth]{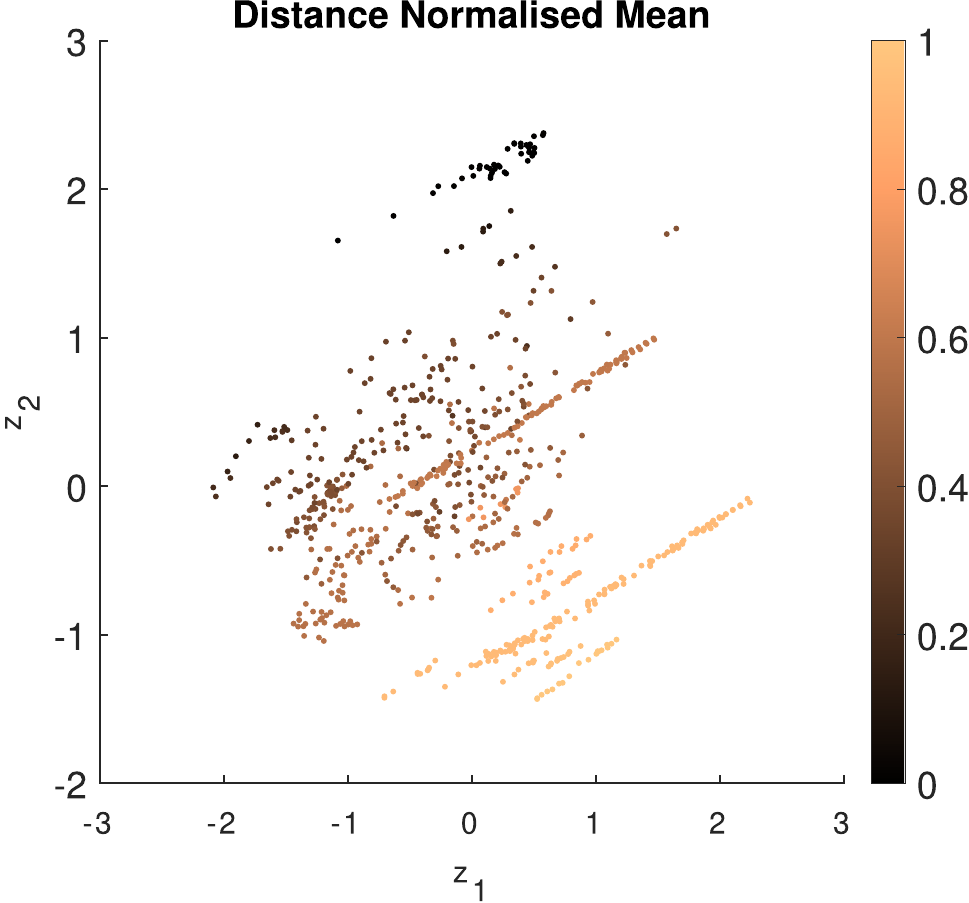}
    \includegraphics[trim=0.0cm 0.0cm 0.0cm 0.0cm, clip=true, width=0.24\textwidth]{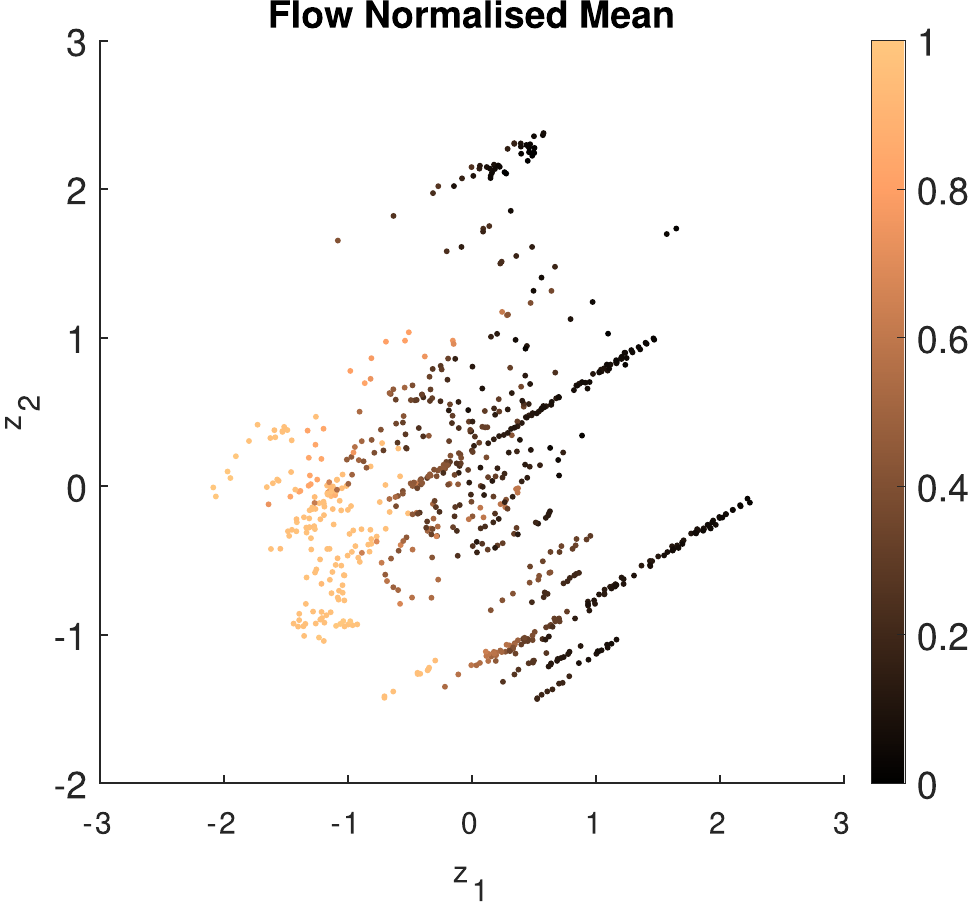}
    \includegraphics[trim=0.0cm 0.0cm 0.0cm 0.0cm, clip=true, width=0.24\textwidth]{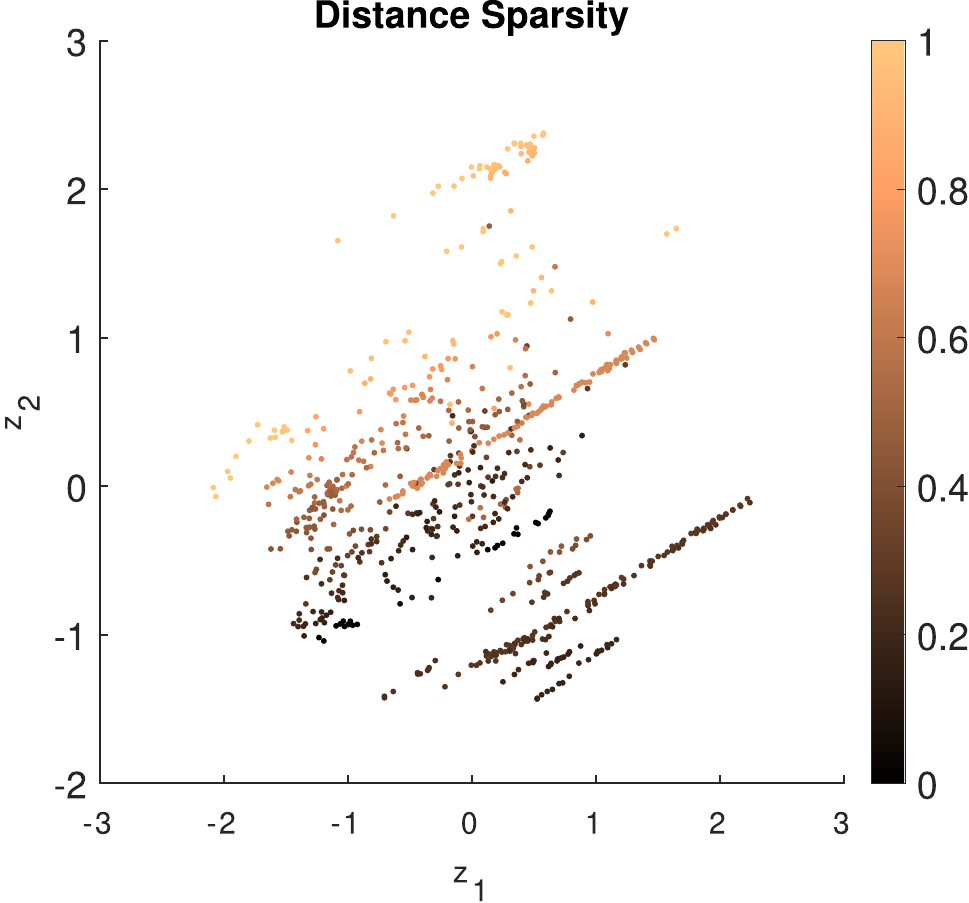}
    \includegraphics[trim=0.0cm 0.0cm 0.0cm 0.0cm, clip=true, width=0.24\textwidth]{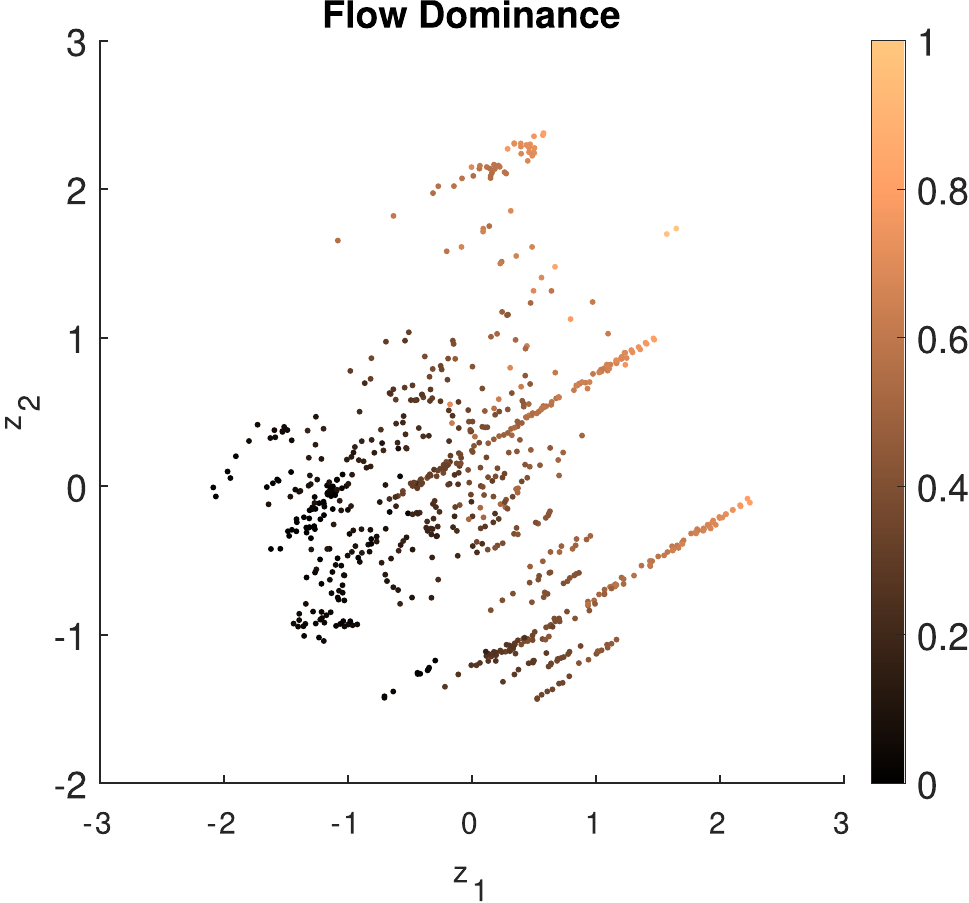} \\
    \includegraphics[trim=0.0cm 0.0cm 0.0cm 0.0cm, clip=true, width=0.24\textwidth]{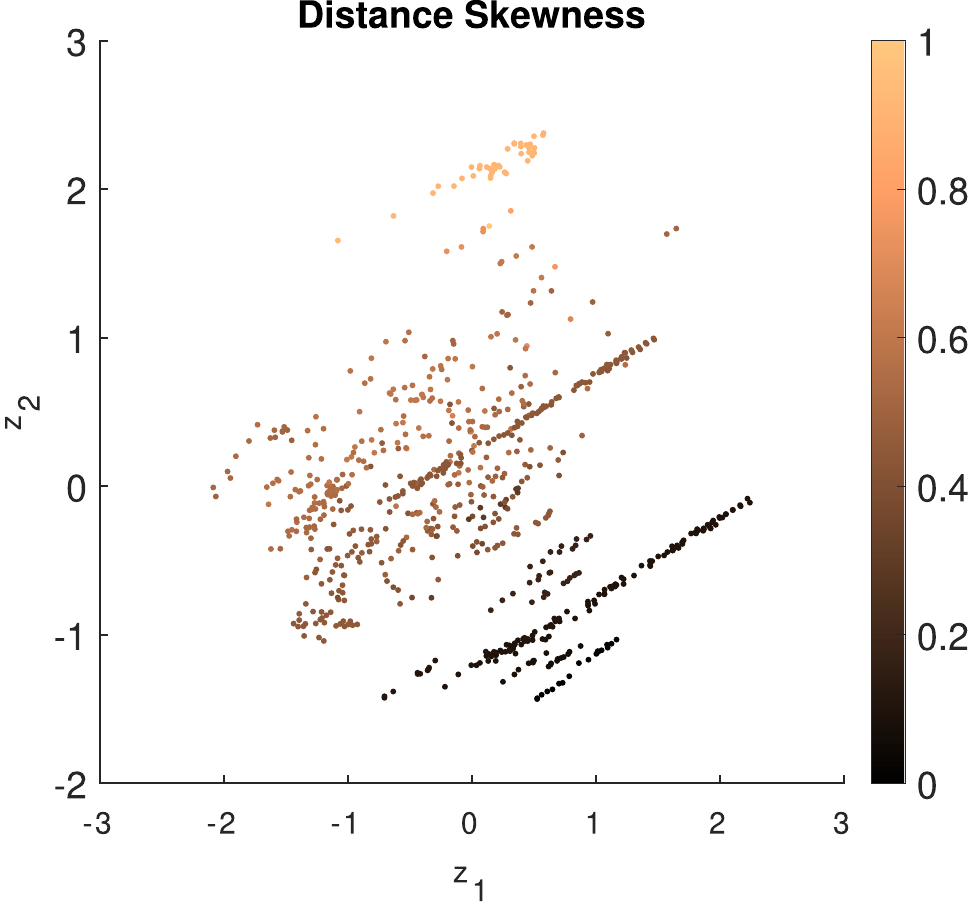}
    \includegraphics[trim=0.0cm 0.0cm 0.0cm 0.0cm, clip=true, width=0.24\textwidth]{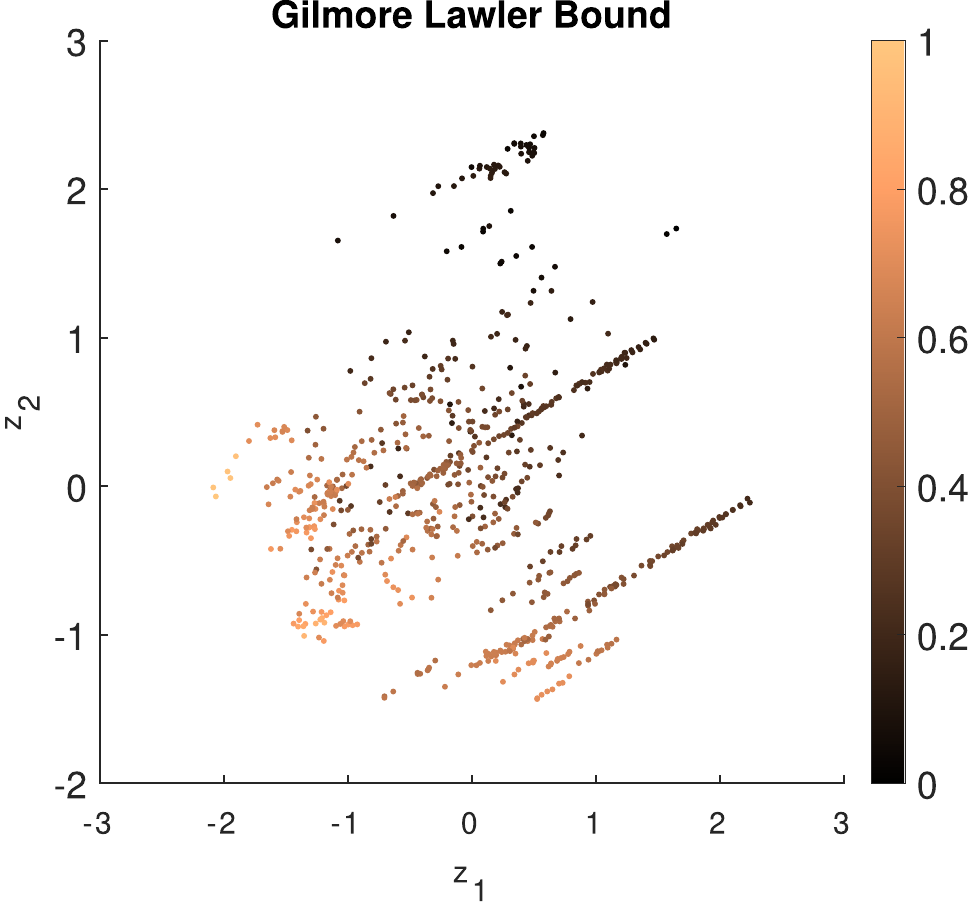}
    \includegraphics[trim=0.0cm 0.0cm 0.0cm 0.0cm, clip=true, width=0.24\textwidth]{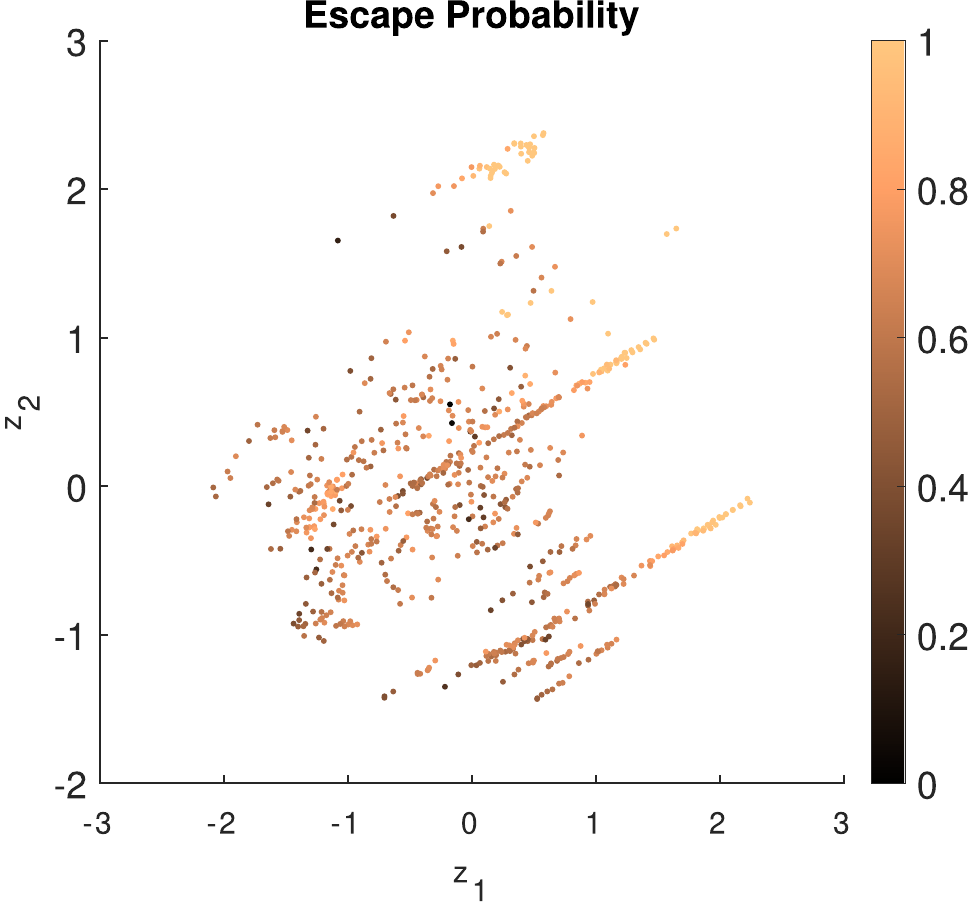}
    \caption{Feature distribution in second space.}
    \label{fig:ext-features}
\end{figure}

The results of this projection are presented in Figures~\ref{fig:ext-classes}-\ref{fig:ext-features}. Figure~\ref{fig:ext-classes} show the location of instance classes of particular interest in the space, Figure~\ref{fig:ext-algorithm} displays the relative performance of the algorithms, and Figure~\ref{fig:ext-features} shows the normalised feature values for each instance. The Distance Sparsity feature is the only feature in common between the two projections.

An immediate observation is that unlike the initial projection \eqref{eq:initialproj} the new projection includes two features derived solely from the flow matrix, the Normalised Mean and Dominance. 
Figure~\ref{fig:ext-features} shows that the variation in both of these features, as well as the Gilmore-Lawler bound, tends to have a consistent gradient from the lower-left to the upper-right of the space. The Escape Probability feature shows a similar though less consistent pattern. On the other hand, the three distance-based features have fairly consistent gradients from the upper-left to the lower-right. 

The instances favourable to MMAS are found in several clusters. The Terminal instances, and most of the flow-cluster instances where MMAS performs strongly, are found along the lower-right edge of the space; these instances are characterised by large values of the Flow Dominance feature, and low values of the Flow Normalised Mean feature. The flow-cluster instances with Hypercube distances and Triangle flows which favour MMAS are on the lower edge of the space, while the Palubeckis instances are on the left side of the space; both of these groups are surrounded by instances which instead favour BMA. Unlike the first instance space, the instances where MMAS performs well are spread around the space and not clearly separated from the instances where BMA is preferred. This suggests that the features being measured are not sufficient to fully explain the variation in algorithm performance within the expanded instance subset, and hence that the ideas proposed in Section~\ref{sec:ispace-analysis} which informed and motivated the new instance classes are not captured by any combination of our existing feature set. 

\subsection{Evolving New Instances}\label{sec:evolving}

To further expand the instance subset and investigate the relationship between the flow matrix and problem difficulty, we apply an genetic algorithm to produce instances near target points in the instance space. The approach of directly representing the problem data of a single QAP instance in the genotype of an individual did not perform well in initial experiments, typically producing degenerate instances (for example with very high or total sparsity in the matrices). Instead, we choose to evolve the parameters of an instance generator to produce instances in the desired region of the instance space, similar to the structured instances evolved in our previous work on the 0-1 Knapsack Problem \citep{Smith-Miles2021KP}.

As we have seen, the properties and difficulty of a QAP instance depends on the complicated relationship between the distance and flow matrices. Our initial experiments indicated that modifying both matrices at the same time has unpredictable effects on the resulting instance. Therefore, the instance generator we apply here uses a fixed distance matrix taken from existing instances, and we evolve parameters which determine the properties of the flow matrix. 

Pseudocode for the flow generator is provided in the supplementary material (Appendix D.4). We have equipped the generator with parameters with the goal that no combination of the parameters will produce a degenerate QAP instance. We evolve five `raw' parameters taking values in $[0,1]$: $p_1$ determining the number of clusters, $p_2$ determining the flow density within a cluster, $p_3$ determining the variation in flows within a cluster, $p_4$ determining the number of links between clusters, and $p_5 $ determining the strength of flows outside a cluster. The actual generator configuration is defined as follows: number of clusters $C = \left(\frac{n}{2}\right)^{p_1}$ rounded to the nearest integer, minimum density within cluster $C_{dens} = p_2$, maximum flow within cluster $C_{max} = \lfloor 100 + 10^{1 + 2p_3} \rfloor$, frequency of noise $N_{freq} = \frac{p_4}{10}$, and maximum noise flow $N_{max} = \lfloor 1 + 99p_5 \rfloor$. 

The objective of the evolutionary algorithm is to produce a set of parameters for this generator such that, when the resulting instances are projected into the $2$D instance space, they tend to be close (in a Euclidean distance sense) to a given target point in the space. We note that there is a random element in both the generation procedure for each instance, and the measurement of the Escape Probability fitness landscape feature which contributes to the $2$D projection. Therefore, a reasonable sample of potential outcomes from each individual (i.e. vector of generator parameters) must be taken to ensure a reliable assessment of their fitness.

In order to expand the existing instance space we selected eight target points on the outside of the space, with less emphasis on the top side since the two algorithms tended to perform similarly on the existing instances in this region. For each combination of distance matrix and target point we run the genetic algorithm with 20 individuals. Each individual parameter vector is used to generate 10 QAP instances. We calculate the distance of each instance from the target point in the $2$D instance space, then order the instances in ascending order of distances. and take the average of the second through fifth distances as the fitness of the individual. Since our goal is to find generator parameters which consistently produce instances with the desired characteristics, rather than parameters which may occasionally produce a lucky result, we take the average of the second through fifth distances in this list as the fitness of the individual parameter vector.

We apply MATLAB's default genetic algorithm implementation with a maximum of 20 generations. Once this is complete, we use the best set of parameters found to generate 50 instances, and keep the closest of these instances to the target point. We generally found that the best few out of the 50 instances had similar distance to the target point. 

Performing the fitness landscape analysis for each instance produced by each individual set of parameters at each generation of the genetic algorithm is quite time-consuming. Therefore, we choose distance matrices of fairly modest size, with variety in terms of structure and position in the instance space: the sko72 instance from QAPLIB with grid-based distances, dre72 from the DreXX test set, three instances (stf60es2, stf60er1 and stf100ep3) from the real-like instance generator of \citet{Stutzle2004Instances}, a uniform random instance xran70A1, a Terminal instance term75_4, and a Hypercube instance hyp64_3. 

\begin{figure}[t]
    \centering
    \subfloat{\includegraphics[trim=0.0cm 0.0cm 0.0cm 0.0cm, clip=true, width=0.38\textwidth]{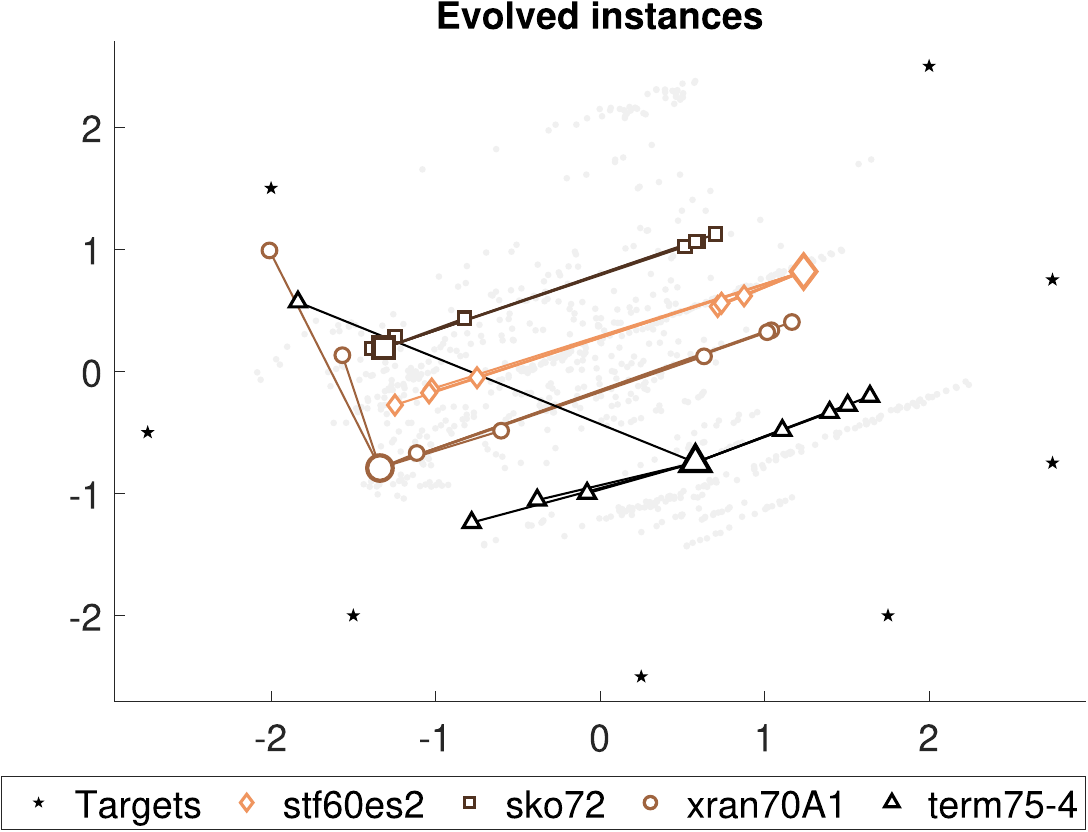}}
    \hspace{2mm}
    \subfloat{\includegraphics[trim=0.0cm 0.0cm 0.0cm 0.0cm, clip=true, width=0.38\textwidth]{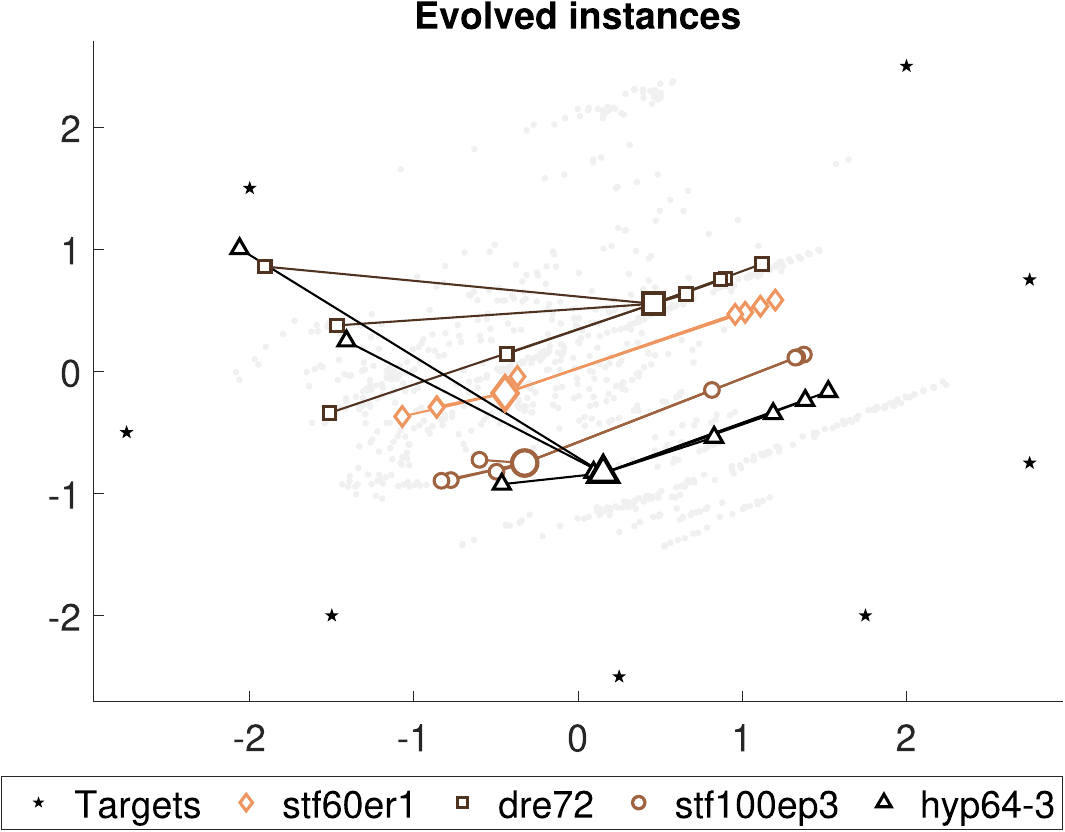}}
    \caption{Positions of evolved instances in the instance space. Large symbols indicate the original instance. Small symbols indicate the instances with evolved flows, with lines connecting them to the original instance. Stars indicate target points.}
    \label{fig:spider}
\end{figure}

\begin{figure}[t]
    \centering
    \subfloat{\includegraphics[trim=0.0cm 0.0cm 0.0cm 0.0cm, clip=true, width=0.38\textwidth]{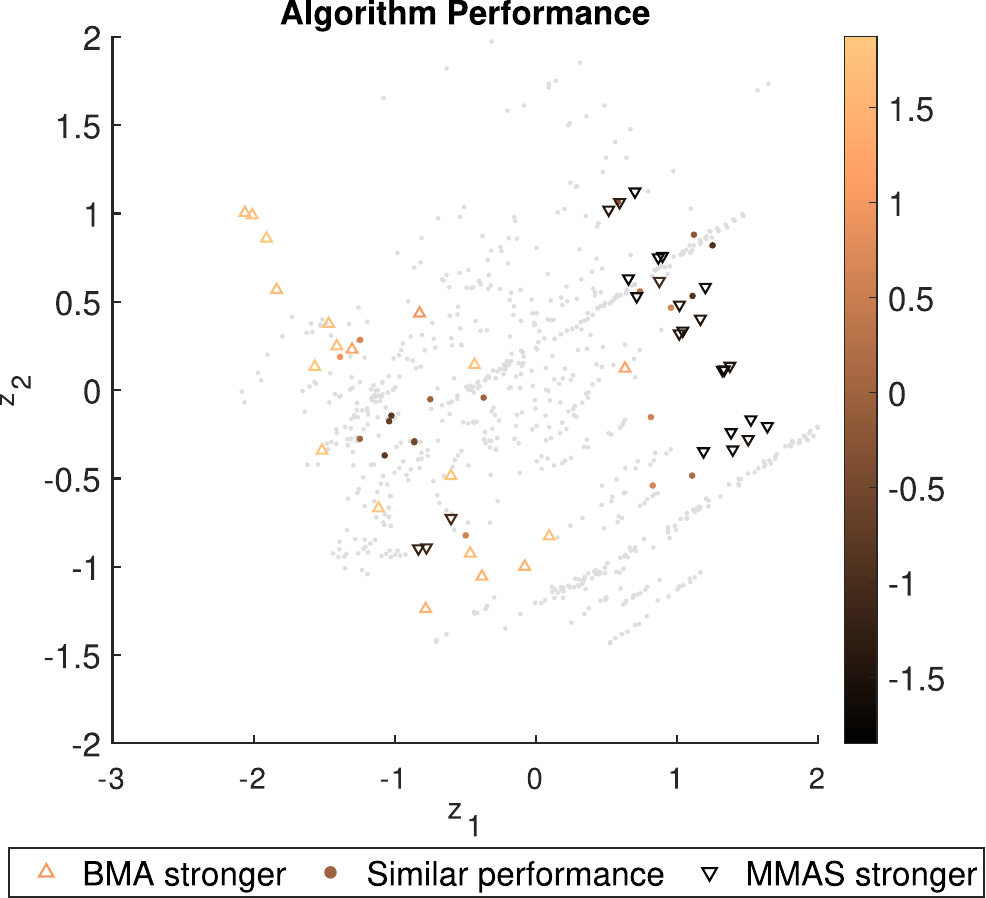}}
    \caption{Comparison in performance (as defined by \eqref{eq:alg-crit}) of the algorithms on the evolved instances. Larger values and lighter colours indicates that BMA is better compared with MMAS. Symbols indicate where each algorithm is clearly superior.}
    \label{fig:evoperf}
\end{figure}

The position of the evolved instances in the instance space is shown in Figure~\ref{fig:spider}. It is immediately apparent that, in this instance space, most of each group of instances starting with the same distance matrix forms a line across the space, in a similar direction as the line formed by the flow-cluster instances. While there is some variation in the values of the Gilmore Lawler Bound and Escape Probability features these do not result in much deviation from the line. A few of the instances with the term75_4, xran70A1, hyp64_3 and dre72 distance matrices are instead found in the far upper-left part of the space. For these unusual instances, the `flow' matrix produced using the generator in this section has actually been identified as the distance matrix by the procedure in Section~\ref{sec:identify-dist}, resulting in very different values for the first five features in the projection.

The performance of the algorithms when applied to these evolved instances is summarised in Figure~\ref{fig:evoperf}. Many of the evolved instances are strongly favourable for one of the two algorithms, with most of the instances favouring MMAS found on the right side of the instance space, in a similar position to the flow-cluster instances with square and tree flows which also favour MMAS. The parameter vectors which produce these instances consistently have high values for the $p_1$ and $p_3$ parameters and low values for $p_4$ and $p_5$, so the resulting flow matrices contain many pairs of facilities with high flow between them and minimal flow between distinct pairs.

Including the evolved instances in the instance subset and generating a new $2$D instance space results in an instance space with very similar characteristics and six of the seven features in common; the Flow Dominance feature is replaced by the Most Dominance feature. Therefore, in the interest of brevity we do not include plots for the new instance space here; they may be found at \url{https://matilda.unimelb.edu.au/matilda/problems/opt/qap}. 
 
\section{Conclusion and Future Work} \label{sec:conclusion}

In this paper we have applied Instance Space Analysis to obtain new insights into the problem space of the Quadratic Assignment Problem, and into the relative strengths and weaknesses of the BMA and MMAS algorithms for solving QAP instances. Our contribution begins with an examination of the theory surrounding the QAP, starting with a review of results which are well-known in the literature. Our discussion of these results emphasises some practical consequences which have not been fully considered in existing predictors of QAP instance difficulty. We have also presented a method for calculation of conditional expectation values which improves on previous work by \citet{ChmielCondExp2019} in terms of computational efficiency. While our application of this result is limited to defining a new feature, the method may also be applicable to QAP algorithm design.

Based on our review of QAP theory, we defined a standardisation process which reliably identifies one of the two problem data matrices of a QAP instance as the distance matrix, then rescales the problem data. This procedure allows us to measure QAP instance features which behave consistently across instances which may superficially appear different, but in fact have similar or identical properties from the point of view of a solver. Taking advantage of this standard form, we have proposed several new features for measuring instance difficulty, such as the Normalised Mean of the values in each matrix and the Distribution Similarity between the distance and flow matrices. Furthermore, existing feature ideas such as the flow and distance dominance behave more consistently once our instance standardisation procedure has been applied.


This study has exposed several weaknesses in the currently available QAP test instances which compromise their ability to provide a fair assessment of algorithm quality. Our initial instance space revealed that the instances contained in QAPLIB have very poor coverage of a large region of the potential instance space. Even when augmented by instances from other sources in the literature and the Hypercube instances proposed in this work, we found that the resulting instance test set of 443 instances over-emphasised the role of the distance matrix in defining the characteristics of the instances. Motivated by this observation and Conjecture~\ref{hyp:1}, we diversified our instance subset with three more sources of instances: 115 hybrid, 180 flow-cluster and 64 evolved instances, each with corresponding subclasses, combining to form an expanded test set of 802 instances. These classes take an existing distance matrix or generator and pair it with a new flow generator, which in many cases produces instances with a significantly different algorithm performance profile. 

Our initial instance space analysis suggested that the evolution-based BMA algorithm tended to outperform or at least perform competitively with the ant colony-based MMAS algorithm on most of the instance space. A superficial examination of the space would suggest that the exceptional instances where MMAS outperforms BMA are characterised by a distance matrix which, after normalisation, appears similar to a matrix which might be produced by distances between points using an ordinary distance metric. However, the new instances we constructed in Section~\ref{sec:more-instances} contained several instances which met this characterisation and yet favoured BMA, or did not meet it and still favoured MMAS. Our Conjecture~\ref{hyp:1}, which instead characterises instances expected to favour MMAS in terms of how successful we would expect a strategy of fixing a few initially good-looking assignments to be, seems to better fit the observed results. Turning this conjecture into a feature which can be calculated in a reasonable amount of time, without prior knowledge about the structure of a QAP instance, is a potential topic for future study.

We note an interesting parallel between the results in this work and a previous study \citep{Smith-Miles2014GCP} projecting Graph Colouring Problem (GCP) instances to a $2$D space using Principal Component Analysis; in the resulting instance space, an evolutionary algorithm is superior on a substantial proportion of the space, but an ant colony algorithm is better in a smaller region corresponding to high algebraic connectivity, high energy, and moderate graph density. Further consideration of the similarities between these GCP instances and the QAP instances on which MMAS outperforms BMA may lead to a greater understanding of the relative strengths of evolutionary and ant-colony metaheuristics for optimisation problems in general.

There is more work to be done towards a library of QAP instances which is fully representative of the potential instance space. 
Just as the comparison between BMA and MMAS has inspired the instances proposed in this work, comparing other algorithms may lead to further ideas for broadening test instance libraries. 
There is potential to expand on our method for evolving new instances by defining new distance and flow generators with suitable parameters, then alternately evolving the parameters for one matrix while holding the other constant.
We also encourage the publication of large instances derived from real-life problems, to ensure that analysis of the QAP remains grounded in practical applications. As computational capabilities continue to increase, the existing smaller real-life instances contained in QAPLIB will become less useful for this purpose.

The generality of the problem space considered in this paper naturally imposes limits on the features which can be measured and the analysis which can be conducted. Therefore, a potential line of future inquiry is to apply ISA or similar techniques to a specific category of QAP instances. An analysis focusing on symmetric QAP instances could apply the more powerful theoretical results which apply to this subclass of the QAP to define new features. Focusing on a set of instances where the locations and distances are related to points in a specific metric space would permit an even broader range of features. Limiting the analysis to instances with a known optimal solution would enable the use of a non-relative performance metric. Considering only instances of the same size would make it considerably easier to both define features which will behave consistently over the entire problem space, and evolve new instances to fill holes in the instance space. 

Finally, a practical algorithm selector for the QAP would ideally not only select the algorithm to apply to a particular instance, but would also tune the parameters of the chosen algorithm to achieve the best possible results. The current test instance libraries and general understanding of the problem space still have significant shortcomings, making the construction of such a selector with reliable strong performance very challenging. If these shortcomings can be addressed, we look forward to the potential for future developments in algorithm selection procedures for the QAP.

\section*{Acknowledgements}

This research was supported by the Australian Research Council through the ARC Training Centre in Optimisation Technologies, Integrated Methodologies and Applications (OPTIMA), grant number IC200100009. This research was also supported by the University of Melbourne’s Research Computing Services and the Petascale Campus Initiative.

    \bibliographystyle{apalike}
    \bibliography{references}

\appendix
\pagebreak
\section{Parameters used in instance generators}

 In the following descriptions we use 
$\disunif{a}{b}$
to represent an integer in $\left\{a,\dots,b\right\}$ chosen using a discrete uniform distribution between $a$ and $b$, and 
$\conunif{a}{b}$
to represent a real number in $[a,b]$ chosen using a continuous uniform distribution. Where generator parameters are specified using a random distribution, each parameter was chosen independently for each instance.

\subsection{St\"utzle and Fernandes (2004) generators}

We generated $20$ instances with each of the $6$ total combinations of distance and flow generator described in Tables~\ref{tab-app:sdist} and \ref{tab-app:sflow}. We included an equal number of instances of sizes $n=60, 80, 100$ and $120$.

\begin{table}[ht]
    \scriptsize
    \centering
    \caption{Parameters used for St\"utzle and Fernandes (2004) distance generators}
    \rowcolors{2}{gray!15}{white} 
    \begin{tabular}{m{3.0cm} m{7.0cm}}
    \toprule
    Subclass & Parameters \\ \midrule
    Euclidean &  Maximum cluster size: $K = \lceil \frac{n}{\conunif{2}{15}} \rceil$ \newline Cluster width: $m = \conunif{20}{50}$ \\ 
    Manhattan & Grid has width $4$ or $10$, with equal probability \\
    \bottomrule
    \end{tabular}
    \label{tab-app:sdist}
\end{table}

\begin{table}[ht]
    \scriptsize
    \centering
    \caption{Parameters used for St\"utzle and Fernandes (2004) flow generators}
    \rowcolors{2}{gray!15}{white} 
    \begin{tabular}{m{3.0cm} m{7.0cm}}
    \toprule
    Subclass & Parameters \\ \midrule
    Random & Flow magnitude parameters: $b = \conunif{1}{7}$, $a = 100^{1/b}$ \newline Sparsity parameter: $sp = \conunif{0.6}{0.8}$ \\ 
    Structured & Threshold parameter: $d = \conunif{10}{50}$ \newline $a,b$ as above \\
    Structured Plus & Probability of extra flows: $p = 0.05$ \newline $a,b,d$ as above \\
    \bottomrule
    \end{tabular}
    \label{tab-app:sflow}
\end{table}

\subsection{Taillard (1995) tai*b generator}

The symmetric instances use the generator as proposed by Taillard (1995). In these instances the distance matrix is based on the Euclidean distance between points $(x_1,y_1), \dots, (x_n, y_n)$. The flow matrix is asymmetric, but as we have observed a QAP problem which is symmetric in only one matrix has an equivalent problem which is symmetric in both matrices. We created an asymmetric variant of the generator which defines the distances according to the following rule, analogous to a tilted plane where more work is required to travel uphill:

$$d_{ij} = \begin{cases} \lceil \sqrt{(1+t)(x_i-x_j)^2 + (y_i-y_j)^2} \rceil & x_i > x_j \\ 
    \lceil \sqrt{(1-t)(x_i-x_j)^2 + (y_i-y_j)^2} \rceil & x_i \leq x_j 
    \end{cases}
$$

The quantity $t$ is a parameter controlling the tilt in the $x$-direction.

\begin{table}[ht]
    \scriptsize
    \centering
    \caption{Parameters used for Taillard (1995) tai*b generator}
    \rowcolors{2}{gray!15}{white} 
    \begin{tabular}{m{3.0cm} m{2.0cm} m{7.0cm}}
    \toprule
    Subclass & \# of Inst. & Parameters \\ \midrule
    Symmetric & 20 & Two instances of each size $n = 30, 35,\dots,120,125$ \newline Flow magnitude parameters: $B = 5$, $A = \conunif{-15}{0}$ \newline Overall radius: $M = 1000$ \newline Cluster radius: $m = \conunif{10}{100}$ \newline Max. cluster size: $K = \lceil \frac{n}{\conunif{2}{15}} \rceil$ \\ 
    Asymmetric & 20 & Tilt parameter: $t = \conunif{0.1}{0.3}$ \newline Other parameters as above\\
    \bottomrule
    \end{tabular}
    \label{tab-app:tai}
\end{table}

\pagebreak

\subsection{QAPSAT generator}

The quantity $M$ used to define the number of A-clauses $m$ in the table below represents a threshold defined by Verel et al. (2024b) believed to separate easy and hard instances. It is defined as $$M := k n^{\alpha_1} m_1^{\alpha_2},$$ with $k \approx 5.23$, $\alpha_1 \approx -0.760$ and $\alpha_2 \approx 0.903$ determined empirically. Note the dependence of $M$ on $n$ and $m_1$.

\begin{table}[ht]
    \scriptsize
    \centering
    \caption{Parameters used for QAPSAT generator}
    \rowcolors{2}{gray!15}{white} 
    \begin{tabular}{m{3.0cm} m{2.0cm} m{7.0cm}}
    \toprule
    Subclass & \# of Inst. & Parameters \\ \midrule
    QAPSAT (easy) & 20 & $n = 25, 30,\dots,115, 120$, one instance of each size \newline Clause size: $k = 3$ \newline Number of B-clauses: $m_1 = 2n$ \newline Number of A-clauses: $m = \lfloor \conunif{0.5}{1.0} \times M \rfloor$\\ 
    QAPSAT (hard) & 20 & Number of A-clauses: $m = \lfloor \conunif{1.5}{2.0} \times M \rfloor$ \newline Other parameters as above\\
    \bottomrule
    \end{tabular}
    \label{tab-app:qapsat}
\end{table}

\subsection{Other generators}

\begin{table}[ht]
    \scriptsize
    \centering
    \caption{Parameters used for other instance generators}
    \rowcolors{2}{gray!15}{white} 
    \begin{tabular}{m{3.0cm} m{1.5cm} m{9.5cm}}
    \toprule
    Subclass & \# of Inst. & Parameters \\ \midrule
    Uniform Random & 40 & $n = 30, 35,\dots,120,125$, $2$ instances of each size \newline Each individual distance and flow chosen from $\disunif{0}{999}$ \\ 
    Terminal & 40 & Base distance and flow: $d_0 = f_0 = 20$ \newline Tree branches: $(s,t,u) = (5,3,3), (5,5,3), (5,5,5)$ or $(3,5,7)$, \newline 10 instances with each setting \newline (Implied instance size: $n = s \times t \times u)$\\
    Hypercube & 40 &  Base distance and flow: $d_0 = f_0 = 20$ \newline Side length and dimension: $(\ell,k) = (2,5), (2,6), (2,7), (3,4)$ or $(5,3)$, \newline 8 instances with each setting \newline (Implied instance size: $n = \ell^k$)\\
    Palubeckis (2000) & 40 & $n = 60, 80, 100, 120$, 10 instances of each size \newline Grid dimension: $n_x = n_y = 7$ \newline Number of graphs: $h=\lfloor\frac{n}{2}\rfloor$ \newline Size of graph bounds: $m_{lower}=3$, $m_{upper} = 19$ \newline Bound on edge weights: $w = 10$ \\
    \bottomrule
    \end{tabular}
    \label{tab-app:other}
\end{table}

\section{Computational comparison between equivalent QAP instances}\label{app:comp-equivalent}

In each table Column 1 gives the name of the original instance in the first row, and the transformations applied to this instance in the following rows. The first transformation consists of exchanging the two problem data matrices, which does not rescale the objective values obtained. Therefore the objective values achieved for this transformed problem are directly comparable with those for the original problem. The remaining four transformations do rescale the objective value of each potential assignment, but each transformation does this in a predictable and consistent way.
In these rows we report the solution values obtained by reversing these rescalings, so that the results in all six rows are directly comparable.

In columns 2 and 3 (corresponding to BMA and MMAS results respectively), the first line is the average solution cost obtained by each algorithm for each instance, followed by the standard deviation in brackets. The second line is the average time required to find the best solution, followed by the standard deviation in brackets. Each algorithm is run 50 times on each instance variant. 

In column 4 we record the dominance of the two matrices; we give the two values in the same order as the matrices appear in the data file. In column 5 we similarly record the sparsities of the two matrices according to the ordinary definition i.e. the proportion of entries which are equal to zero.

We observe that while some of these transformations have substantial effects on the dominance and sparsity values, the performance of the algorithms is not substantially affected.

\begin{table}[hptb] \footnotesize
    \centering
    \caption{Comparison of algorithm performance and selected features between the dre56 instance and other equivalent instances. See text at beginning of the appendix for details.} 
		\begin{tabular}{c c c c c c}
			\toprule
			Instance & BMA & MMAS & Dominances & Sparsities \\
			\midrule
            \rowcolor{white} & 1232.60 (99.14) & 1617.76 (64.28) & 0.078 & 0.938 \\ \rowcolor{white} \multirow{-2}{*}{dre56 = $QAP(A,B)$} & 63.61s (47.86s) & 80.09s (54.48s) & 0.007 & 0.018 \\ 
        \rowcolor{gray!12} & 1192.68 (116.85) & 1623.84 (65.37) & 0.007 & 0.018 \\ \rowcolor{gray!12} \multirow{-2}{*}{$QAP(B,A)$} & 63.90s (46.65s) & 68.02s (51.12s) & 0.078 & 0.938 \\ 
        \rowcolor{white} & 1214.88 (107.85) & 1624.44 (65.73) & 0.078 & 0.938 \\ \rowcolor{white} \multirow{-2}{*}{$QAP(2A, 2B)$} & 65.09s (51.84s) & 78.95s (49.50s) & 0.007 & 0.018 \\ 
        \rowcolor{gray!12} & 1182.44 (104.72) & 1653.20 (77.33) & 0.000 & 0.000 \\ \rowcolor{gray!12} \multirow{-2}{*}{$QAP(A + 10 J_n, B + 10 J_n)$} & 78.07s (50.90s) & 75.25s (53.48s) & 0.001 & 0.000 \\ 
        \rowcolor{white} & 1203.40 (109.01) & 1616.20 (94.11) & 0.078 & 0.938 \\ \rowcolor{white} \multirow{-2}{*}{$QAP(\reduced(A),\reduced(B))$} & 46.81s (38.66s) & 74.95s (51.70s) & 0.009 & 0.080 \\ 
        \rowcolor{gray!12} & 1186.92 (100.53) & 1650.56 (51.00) & 0.001 & 0.027 \\ \rowcolor{gray!12} \multirow{-2}{*}{$QAP(\reduced(-A),\reduced(-B))$} & 65.61s (52.14s) & 66.24s (51.36s) & 0.009 & 0.119 \\ \bottomrule
			\end{tabular}
		
		\label{table:results-tmp-2}
\end{table}

\begin{table}[hptb] \footnotesize
\centering
\caption{Comparison of algorithm performance and selected features between the tai60a instance and other equivalent instances. See text at beginning of the appendix for details.} 
		\begin{tabular}{c c c c c c }
			\toprule
			Instance & BMA & MMAS & Dominances & Sparsities \\
			\midrule 
            \rowcolor{white} & 7240875.76 (9661.34) & 7317982.20 (16216.61) & 0.008 & 0.024 \\ \rowcolor{white} \multirow{-2}{*}{tai60a = $QAP(A,B)$} & 82.30s (56.02s) & 67.30s (56.87s) & 0.008 & 0.027 \\ 
\rowcolor{gray!12} & 7240961.40 (9532.80) & 7319907.60 (19838.76) & 0.008 & 0.027 \\ \rowcolor{gray!12} \multirow{-2}{*}{$QAP(B,A)$} & 69.91s (51.91s) & 49.15s (52.95s) & 0.008 & 0.024 \\ 
\rowcolor{white} & 7242344.92 (8054.25) & 7314545.44 (19089.69) & 0.008 & 0.024 \\ \rowcolor{white} \multirow{-2}{*}{$QAP(2A, 2B)$} & 64.38s (53.66s) & 48.69s (49.39s) & 0.008 & 0.027 \\ 
\rowcolor{gray!12} & 7243626.00 (5774.86) & 7322690.84 (22777.99) & 0.006 & 0.000 \\ \rowcolor{gray!12} \multirow{-2}{*}{$QAP(A + 10 J_n, B + 10 J_n)$} & 74.62s (52.71s) & 49.02s (48.66s) & 0.006 & 0.000 \\ 
\rowcolor{white} & 7242561.92 (9289.77) & 7314579.48 (20306.07) & 0.008 & 0.024 \\ \rowcolor{white} \multirow{-2}{*}{$QAP(\reduced(A),\reduced(B))$} & 76.18s (47.23s) & 52.62s (51.74s) & 0.008 & 0.027 \\ 
\rowcolor{gray!12} & 7239655.04 (8745.51) & 7321162.64 (20184.21) & 0.008 & 0.026 \\ \rowcolor{gray!12} \multirow{-2}{*}{$QAP(\reduced(-A),\reduced(-B))$} & 95.71s (62.06s) & 54.07s (50.66s) & 0.008 & 0.028 \\ \bottomrule 

			\end{tabular}
		
		\label{table:results-tmp-3}
\end{table}

\begin{table}[hptb] \footnotesize
    \centering
    \caption{Comparison of algorithm performance and selected features between the term45_1 instance and other equivalent instances. See text at beginning of the appendix for details.}
		\begin{tabular}{c c c c c c}
			\toprule
			Instance & BMA & MMAS & Dominances & Sparcities \\
			\midrule
            \rowcolor{white} & 45155.24 (124.52) & 45048.00 (0.00) & 0.006 & 0.022 \\ \rowcolor{white} \multirow{-2}{*}{term45_1 = $QAP(A,B)$} & 40.35s (44.03s) & 15.66s (23.35s) & 0.063 & 0.836 \\ 
\rowcolor{gray!12} & 45179.88 (109.22) & 45048.00 (0.00) & 0.063 & 0.836 \\ \rowcolor{gray!12} \multirow{-2}{*}{$QAP(B,A)$} & 23.56s (31.65s) & 19.63s (19.87s) & 0.006 & 0.022 \\ 
\rowcolor{white} & 45176.00 (126.95) & 45048.00 (0.00) & 0.006 & 0.022 \\ \rowcolor{white} \multirow{-2}{*}{$QAP(2A, 2B)$} & 22.73s (28.34s) & 23.38s (27.91s) & 0.063 & 0.836 \\ 
\rowcolor{gray!12} & 45155.32 (113.46) & 45049.88 (13.29) & 0.004 & 0.000 \\ \rowcolor{gray!12} \multirow{-2}{*}{$QAP(A + 10 J_n, B + 10 J_n)$} & 32.44s (37.56s) & 31.00s (36.87s) & 0.004 & 0.000 \\ 
\rowcolor{white} & 45193.48 (134.38) & 45048.00 (0.00) & 0.006 & 0.024 \\ \rowcolor{white} \multirow{-2}{*}{$QAP(\reduced(A),\reduced(B))$} & 37.23s (37.98s) & 20.83s (20.89s) & 0.063 & 0.836 \\ 
\rowcolor{gray!12} & 45164.04 (103.86) & 45055.72 (26.86) & 0.024 & 0.081 \\ \rowcolor{gray!12} \multirow{-2}{*}{$QAP(\reduced(-A),\reduced(-B))$} & 33.66s (39.66s) & 37.22s (33.38s) & 0.002 & 0.026 \\ \bottomrule
			\end{tabular}
		 
		\label{table:results-tmp-4}
\end{table}

\pagebreak
\section{Proofs for propositions in Theoretical Foundations section}\label{app:sim-proofs}

Recall that we define $\Nneq^2$ as the subset of $N^2$ containing no 2-tuples with repeating elements:
\begin{equation*}
 \Nneq^2 := \left\{ (i_1,i_2) \in N^2 \mid i_1 \neq i_2 \right\}
 \end{equation*}
 We observe that
 \begin{equation*}
    (i_1, i_2) \in \Nneq^2 \iff (i_2, i_1) \in \Nneq^2
 \end{equation*}
 and
 \begin{equation*}
    (i_1, i_2) \in \Nneq^2 \iff (\varphi(i_1), \varphi(i_2)) \in \Nneq^2 \quad \forall \varphi \in S_n,
 \end{equation*}
 which allows for convenient reordering of sums over $M^2$.
    In the following results we will implicitly apply the following reformulation of the QAP cost function:
    \begin{equation} \label{eq:QAP-cost-function-2}
	    \costf_{A,B}(\varphi):=\sum_{(i,j)\in M^2} a_{ij} b_{\varphi(i)\varphi(j)} + \sum_{i\in N} a_{ii} b_{\varphi(i) \varphi(i)}
    \end{equation}
    The first sum in this decomposition contains all of the cost terms which depend on two distinct assignments $\varphi(i)$ and $\varphi(j)$, while the second sum contains the terms which only depend on a single assignment $\varphi(i)$. Note that $\costf_{A,B}(\varphi)$ contains no terms of the form $a_{ij} b_{rr}$ or $a_{ii} b_{rs}$ for any $(i,j),  (r,s) \in M^2$.
     
     \begin{lemma} \label{lemma-apx:isomorphic-distpres}
        For all $\theta, \tau \in S_n$,  $f_{\theta,\tau}(\varphi) = \theta \circ \varphi \circ \tau$ is an automorphism of $S_n$.
    \end{lemma}
    \begin{proof}
        $f^{-1}_{\theta,\tau}(\varphi) = \theta^{-1} \circ \varphi \circ \tau^{-1} = f_{\theta^{-1},\tau^{-1}}(\varphi)$, so $f_{\theta,\tau}$ is bijective. 
        
        ($\implies$) Let $\varphi$ and $\pi$ be distinct elements of $S_n$ satisfying $d_C(\varphi,\pi) = 1$. Then there exists distinct $i$ and $j$ in $N$ such that $\varphi = \sigma_{i,j} \circ \pi.$
        Now we have
        \begin{eqnarray*}
	        \theta \circ \varphi \circ \tau &=& \theta \circ \sigma_{i,j} \circ \pi \circ \tau \\ 
	        &=& \sigma_{\theta(i),\theta(j)} \circ \theta \circ \pi \circ \tau \\
	\end{eqnarray*}
        and $\theta \circ \varphi \circ \tau \neq \theta \circ \pi \circ \tau$, so $d_C(\theta \circ \varphi \circ \tau, \theta \circ \pi \circ \tau) = 1$ as desired to prove the forward implication.
        
        ($\impliedby$) We apply the forward implication using the function $f_{\theta^{-1},\tau^{-1}}$:
        $$d_C(\theta \circ \varphi \circ \tau, \theta \circ \pi \circ \tau) = 1 \implies d_C(\theta^{-1} \circ \theta \circ \varphi \circ \tau \circ \tau^{-1},\theta^{-1} \circ \theta \circ \pi \circ \tau \circ \tau^{-1}) = d_C(\varphi, \pi) = 1$$
    \end{proof}
    
    \begin{proposition}
        Let $\tau, \theta \in S_n$. Let $A = [a_{ij}]$ and $B = [b_{ij}]$ be arbitrary real-valued $n \times n$ matrices and let $X = [x_{ij}]$ and $Y = [y_{ij}]$ be defined by $x_{ij} = a_{\tau(i)\tau(j)}$  and $y_{ij} = b_{\theta(i)\theta(j)}$ for all $i,j$ in $N$. Then $QAP(A,B)$ and $QAP(X,Y)$ are similar QAP instances. 
    \end{proposition}
    \begin{proof}
    By substituting $k = \tau(i)$ and $l = \tau(j)$ and reordering the sums, we have
    \begin{eqnarray*}
        \costf_{X,Y}(\varphi) &=& \sum_{i \in N}\sum_{j \in N} x_{ij} y_{\varphi(i)\varphi(j)} \\
        &=& \sum_{i \in N}\sum_{j \in N} a_{\tau(i)\tau(j)} b_{\theta \circ \varphi(i)\theta \circ \varphi(j)} \\
        &=& \sum_{k \in N}\sum_{l \in N} a_{kl} b_{\theta \circ \varphi \circ \tau^{-1}(k)\theta \circ \varphi \circ \tau^{-1}(l)} \\
        &=& \costf_{A,B}(\theta \circ \varphi \circ \tau^{-1})
    \end{eqnarray*}
        The desired result follows from Lemma~\ref{lemma-apx:isomorphic-distpres}.
    \end{proof}

    \begin{lemma} \label{lemma-apx:inverse-distpres}
        $f(\varphi) = \varphi^{-1}$ is an automorphism of $S_n$.
    \end{lemma}
    \begin{proof}
        $f^{-1}(\varphi) = \varphi^{-1} = f(\varphi)$, so $f$ is bijective. 
        
        ($\implies$) Let $\varphi$ and $\pi$ be distinct elements of $S_n$ satisfying $d_C(\varphi,\pi) = 1$. Then there exists distinct $i$ and $j$ in $N$ such that $\varphi = \sigma_{i,j} \circ \pi.$
        Now we have
        \begin{eqnarray*}
	        \varphi^{-1} &=& \left( \sigma_{i,j} \circ \pi \right)^{-1} \\ 
	        &=& \left( \pi \circ \sigma_{\pi^{-1}(i),\pi^{-1}(j)} \right)^{-1} \\
                &=& \sigma^{-1}_{\pi^{-1}(i),\pi^{-1}(j)} \circ \pi^{-1} \\
                &=& \sigma_{\pi^{-1}(i),\pi^{-1}(j)} \circ \pi^{-1}
	\end{eqnarray*}
        and $\varphi^{-1} \neq \pi^{-1}$, so $d_C(\varphi^{-1}, \pi^{-1}) = 1$ as desired to prove the forward implication.
        
        ($\impliedby$) We apply the forward implication:
        $$d_C(\varphi^{-1}, \pi^{-1}) = 1 \implies d_C((\varphi^{-1})^{-1}, (\pi^{-1})^{-1}) = d_C(\varphi, \pi) = 1$$
    \end{proof}
    
    \begin{proposition}
        $QAP(A,B)$ and $QAP(B,A)$ are similar QAP instances.
    \end{proposition}
    \begin{proof}
        By substituting $k = \varphi(i)$ and $l = \varphi(j)$ and reordering the sums, we have $$\costf_{A,B}(\varphi) = \sum_{i \in N} \sum_{j \in N} a_{ij} b_{\varphi(i)\varphi(j)} = \sum_{k \in N} \sum_{l \in N} a_{\varphi^{-1}(k)\varphi^{-1}(l)} b_{kl} = \costf_{B,A}(\varphi^{-1}).$$
        The desired result follows from Lemma \ref{lemma-apx:inverse-distpres}.
    \end{proof}

    \begin{proposition}
	    Let $B$ be a symmetric matrix (i.e. $b_{ij} = b_{ji}$ for all $i,j$ in $N^2$). Let $C = \left[c_{ij}\right]$ be a skew-symmetric $n \times n$ matrix (i.e. $c_{ij} = - c_{ji}$ and $c_{ii} = 0$ for all $(i,j)$ in $M^2$). Then for all $\varphi \in S_n$, $$\costf_{A,B} (\varphi) = \costf_{(A+C),B}(\varphi)$$ and hence $QAP(A,B)$ and $QAP(A+C,B)$ are equivalent QAP instances.
\end{proposition}
	\begin{proof}
    Applying the formulation of the cost function in~\eqref{eq:QAP-cost-function-2}, symmetry of $B$ and skew-symmetry of $C$:
	    \begin{eqnarray*}
	        \costf_{A,B} (\varphi) &=& \sum_{(i,j)\in M^2} a_{ij} b_{\varphi(i)\varphi(j)} + \sum_{i\in N} a_{ii} b_{\varphi(i) \varphi(i)} \\ 
            &=& \sum_{(i,j)\in M^2} \frac{1}{2} (a_{ij} + a_{ji}) b_{\varphi(i)\varphi(j)}  + \sum_{i \in N} a_{ii} b_{\varphi(i)\varphi(i)} \\ 
	        &=& \sum_{(i,j)\in M^2} \frac{1}{2} (a_{ij}+c_{ij}+a_{ji}+c_{ji})b_{\varphi(i)\varphi(j)} + \sum_{i=1}^n (a_{ii} + c_{ii}) b_{\varphi(i)\varphi(i)} \\
            &=& \sum_{(i,j)\in M^2} (a_{ij} + c_{ij}) b_{\varphi(i)\varphi(j)} + \sum_{i\in N} a_{ii} b_{\varphi(i) \varphi(i)} \\ 
	        &=& \costf_{(A+C),B}(\varphi).
	    \end{eqnarray*}
	\end{proof}
	
	\begin{proposition}
	    Let $c_1$ and $c_2$ be arbitrary positive constants in $\mathbb{R}$. Then for all $\varphi \in S_n$
	    $$\costf_{c_1 A,c_2 B}(\varphi) = c_1 c_2 \costf_{A,B}(\varphi)$$
	    and hence $QAP(A,B)$ and $QAP(c_1 A, c_2 B)$ are equivalent QAP instances.
	\end{proposition}
        \begin{proof}
	    \begin{equation*}
	    \costf_{c_1 A,c_2 B}(\varphi)=\sum_{i\in N} \sum_{j \in N } c_1 c_2 a_{ij} b_{\varphi(i)\varphi(j)} = c_1 c_2 \sum_{i\in N} \sum_{j \in N}  a_{ij} b_{\varphi(i)\varphi(j)} = c_1 c_2 \costf_{A,B}(\varphi)
	    \end{equation*}
	\end{proof}
	
	\begin{proposition}
	    Let $c_1$ and $c_2$ be arbitrary constants in $\mathbb{R}$. Then 
	    for all $\varphi \in S_n$
	    $$\costf_{(A+c_1 I_n),(B+c_2 I_n)}(\varphi) = \costf_{A,B}(\varphi) + c_2 \trace{A} + c_1 \trace{B} + n c_1 c_2$$
	    and hence $QAP(A,B)$ and $QAP(A+c_1 I_n, B+c_2 I_n)$ are equivalent QAP instances.
	\end{proposition}
	\begin{proof}
	    \begin{eqnarray*}
	    \costf_{(A+c_1 I_n),(B+c_2 I_n)}(\varphi) &=& \sum_{(i,j)\in M^2} a_{ij} b_{\varphi(i)\varphi(j)} + \sum_{i \in N} \left[ (a_{ii} + c_1) (b_{\varphi(i)\varphi(i)} + c_2)  \right] \\
	    &=& \sum_{(i,j)\in M^2} a_{ij} b_{\varphi(i)\varphi(j)} + \sum_{i \in N} \left[ (a_{ii}b_{\varphi(i)\varphi(i)} + a_{ii} c_2 + b_{\varphi(i)\varphi(i)} c_1 + c_1 c_2  \right] \\
	    &=& \sum_{(i,j)\in M^2} a_{ij} b_{\varphi(i)\varphi(j)} + \sum_{i \in N} a_{ii}b_{\varphi(i)\varphi(i)} + \sum_{i \in N} a_{ii} c_2 +  \sum_{i \in N} b_{ii} c_1 + n c_1 c_2 \\
	    &=& \costf_{A,B}(\varphi) + c_2 \trace(A) + c_1 \trace(B) + n c_1 c_2
	    \end{eqnarray*}
	\end{proof}
	
	\begin{proposition}
	    Let $c_1$ and $c_2$ be arbitrary constants in $\mathbb{R}$. Then 
	    for all $\varphi \in S_n$
	    $$\costf_{A+c_1 (\ones_n - I_n),B+c_2 (\ones_n - I_n)}(\varphi) = \costf_{A,B}(\varphi) + \sum_{i\in N} \sum_{j \in N \setminus \left\{ i \right\} } a_{ij} c_2 + \sum_{i\in N} \sum_{j \in N \setminus \left\{ i \right\} } b_{ij} c_1 + n (n-1) c_1 c_2$$
	    and hence $QAP(A,B)$ and $QAP(A+c_1 (\ones_n - I_n), B+c_2 (\ones_n - I_n))$ are equivalent QAP instances.
	\end{proposition}
	\begin{proof}
	    \begin{eqnarray*}
	    \costf_{A+c_1 (\ones_n - I_n),B+c_2 (\ones_n - I_n)}(\varphi) &=& \sum_{(i,j)\in M^2} \left[ (a_{ij} + c_1) (b_{\varphi(i)\varphi(j)} + c_2) \right] + \sum_{i \in N} a_{ii} b_{\varphi(i)\varphi(i)}  \\
	    &=& \sum_{(i,j)\in M^2} \left[ a_{ij} b_{\varphi(i)\varphi(j)} + a_{ij} c_2 + b_{\varphi(i)\varphi(j)} c_1 + c_1 c_2 \right] + \sum_{i \in N} a_{ii} b_{\varphi(i)\varphi(i)} 
	    \\
	    &=& \costf_{A,B}(\varphi) + \sum_{(i,j)\in M^2} a_{ij} c_2 + \sum_{(i,j)\in M^2} b_{ij} c_1 + n (n-1) c_1 c_2
	    \end{eqnarray*} 
	\end{proof}

    \begin{proposition}
        For all $\varphi \in S_n$ $$\costf_{-A,-B}(\varphi) = \costf_{A,B}(\varphi)$$ and hence $QAP(A,B)$ and $QAP(-A,-B)$ are equivalent QAP instances.
    \end{proposition}
    \begin{proof}
        \begin{equation*}
            \costf_{-A,-B}(\varphi)=\sum_{i=1}^n\sum_{j=1}^n \left(-a_{ij}\right) \left(-b_{\varphi(i)\varphi(j)}\right)=\sum_{i=1}^n\sum_{j=1}^n a_{ij} b_{\varphi(i)\varphi(j)}=\costf_{A,B}(\varphi)
        \end{equation*}
    \end{proof}

\pagebreak
\section{Algorithm pseudocode}\label{app:algs}

\subsection{TRIPOD Score} \label{sec-app:tripod}

This algorithm measures the degree to which the input matrix (post normalisation) respects the triangle inequality, and hence resembles a distance matrix. If the triangle inequality is fully respected, the algorithm will output the maximum value of $1$.

Recall that $M^2$ is the set of all $(i,j)$ in $\left\{1, \dots, n\right\}^2$ satisfying the condition $i \neq j$. The strict inequality on line~\ref{algline:triangle-strict} is intended to avoid misidentification of flow matrices with many zero entries. 

\begin{algorithm}[ht]
\caption{TRiangle Inequality Property Of Distances (TRIPOD) score} \label{alg:triangle-inequality}
\begin{algorithmic}[1]
\State {\bf Assumption:} $A$ is a $n \times n$ matrix from a QAP instance in standard form 
\State $\mu \gets \textrm{mean of all elements of $A$ not on main diagonal}$
\State $\sigma \gets \textrm{standard deviation of all elements of $A$ not on main diagonal}$
\State $D \gets \textrm{empty $n \times n$ matrix}$
\ForAll{$(i,j) \in M^2$}

\State $d_{ij} \gets \max\left\{0, \frac{a_{ij} - \mu}{\sigma} + \alpha\right\}$

\EndFor

\State $t \gets 0$
\State $p \gets 0$

\ForAll{$(i,j) \in M^2$}
\ForAll{$k \in \left\{1, \dots, n\right\} \setminus  \left\{i,j\right\}$} 
\If{$d_{ij} > \max\left\{d_{ik}, d_{kj}\right\}$} \label{algline:triangle-strict} 

\State $t \gets t+1$
\State $p \gets p+\max\left\{0,\min(\beta, d_{ij} - d_{ik} - d_{kj})\right\}$

\EndIf
\EndFor
\EndFor

\State {\bf return} $\max\left\{0,\min\left\{1, \frac{t-p}{t}\right\}\right\}$

\end{algorithmic}
\end{algorithm}

\subsection{Diversity Measure} \label{sec-app:diversity}

This algorithm measures the degree of diversity in the input matrix. It will return $0$ if all values in the input matrix are identical, and return $1$ if all values in the input matrix are distinct.

\begin{algorithm}[ht]
    \caption{Diversity measure} \label{alg:diversity}
    \begin{algorithmic}[1]
    \State {\bf Assumption:} $A$ is a $n \times n$ matrix from a QAP instance in standard form 
    \State $D \gets \textrm{empty $n \times n$ matrix}$
    
    \ForAll{$(i,j) \in N^2$}
    
    \State $d_{ij} \gets \textrm{number of entries in A equal to $a_{ij}$}$
    
    \EndFor

    \State $\mu \gets \textrm{mean of all entries in $D$}$
    
    \State {\bf return} $\frac{1}{2} \left( 2 - \frac{\ln(\mu)}{\ln(n)} \right) $
    
    \end{algorithmic}
\end{algorithm}

\subsection{Near Similarity measure} \label{sec-app:nearsim}

This algorithm interprets the input matrix as a distance matrix, and measures the degree to which locations with a short distance between them do or do not have similar properties. The algorithm initially standardises the input matrix so that all values are between $0$ and $1$, and constructs the $D^{sym}$ matrix to determine which locations are close to each other in at least one direction. Next, for each location $i$ the algorithm constructs a vector $s_i$ representing which locations have similar properties to location $i$, and a vector $c_i$ representing which locations are closest to location $i$. The algorithm then returns the (appropriately normalised) dot products of $c_i$ and $s_i$. If this dot product is large, it signifies that location $i$ tends to have similar properties to other locations whose distance directly to or from $i$ is small.

    \begin{algorithm}[ht]
    \caption{Near Similarity measure} \label{alg:near-similarity}
    \begin{algorithmic}[1]
    \State {\bf Assumption:} $A$ is a $n \times n$ matrix from a QAP instance in standard form 
    \State $D \gets \reduced (A)$
    \State $m \gets \max \left\{D\right\}$
    \ForAll{$i \in \left\{1, \dots, n\right\}$} 
    \State $d_{ii} \gets 0$
    \ForAll{$j \in \left\{1, \dots, n\right\} \setminus \left\{ i \right\}$} 
    \State $d_{ij} \gets \frac{d_{ij}}{m}$
    \EndFor
    \EndFor
    \State $D^{sym} \gets \textrm{empty $n \times n$ matrix}$
    \ForAll{$i \in \left\{1, \dots, n\right\}$} 
    \ForAll{$j \in \left\{1, \dots, n\right\} \setminus \left\{ i \right\}$} 
    \State $d^{sym}_{ij} \gets \frac{1}{3} \left( d_{ij} + d_{ji} + \min\left\{d_{ij}, d_{ji}\right\} \right)$
    \EndFor
    \EndFor
    \ForAll{$i \in \left\{1, \dots, n\right\}$} 
        \State $c^i \gets \textrm{empty vector of length $n$}$
        \State $s^i \gets \textrm{empty vector of length $n$}$
        \ForAll{$j \in \left\{1, \dots, n\right\} \setminus \left\{ i \right\}$} 
        \State $v^i \gets \textrm{vector containing $d_{ik}$ and $d_{ki}$ for all $k \neq i, j$}$
        \State $v^j \gets \textrm{vector containing $d_{jk}$ and $d_{kj}$ for all $k \neq i, j$, in same order}$
        \State $w \gets \textrm{empty vector of length $2n-4$}$
        \ForAll {$k \in \left\{1, \dots, 2n-4\right\}$}
        \State $w_k \gets \max\left\{1 - v^i_k v^j_k,\frac{1}{n}\right\}$
        \EndFor
        \State $s^i_j \gets 1 - \frac{\sum_{k=1}^{2n-4}\left|v^i_k - v^j_k\right| w_k}{\sum_{k=1}^{2n-4} w_k}$
        \State $c^i_j \gets (d^{sym}_{ij} - 1)^2$
        \EndFor
    \EndFor
    \State \Return $\sum_{i=1}^n \frac{\sum_{j=1}^{n} s^i_j c^i_j }{\sum_{j=1}^{n} n c^i_j}$
    
    \end{algorithmic}
    \end{algorithm}

\subsection{Flow matrix generator with evolvable parameters} \label{sec-app:flow-evolve}

This algorithm defines a flow matrix generator with parameters defined so that any combination of parameter settings should result in a reasonable, non-degenerate problem. Lines \ref{algline:evoflow-clu1} to \ref{algline:evoflow-clu2} determine which facilities are assigned to which of the $C$ clusters. Lines \ref{algline:evoflow-in1} to \ref{algline:evoflow-in2} assign the flows within each cluster according to the parameters $C_{dens}$ and $C_{max}$; the loop between lines \ref{algline:evoflow-conn1} and \ref{algline:evoflow-conn2} guarantees that the flows within each cluster form a connected graph. Finally, lines \ref{algline:evoflow-out1} to \ref{algline:evoflow-out2} add additional one-way flows which may link the clusters together, according to the parameters $N_{freq}$ and $N_{max}$.

\begin{algorithm}[ht]
    \caption{Flow matrix generator with evolvable parameters} \label{alg:evoflow}
    \begin{algorithmic}[1]
    \State {\bf Input parameters:} $n \in \mathbb{N}$, $C \in \left\{1,\dots,\lceil \frac{n}{2} \rceil \right\}$, $C_{dens} \in [0,1]$, $C_{max} \in \left\{110, \dots, 1100\right\}$, $N_{freq} \in [0,\frac{1}{10}]$, $N_{max} \in [1,100]$
    \State $F \gets \textrm{empty $n \times n$ matrix}$
    \State $S_\sigma \gets \left\{2\sigma-1, 2\sigma\right\}$ for all $\sigma \in \left\{1, C-1\right\}$ \label{algline:evoflow-clu1}
    \If{$2C > n$}
         \State $S_{C} \gets \left\{n-1\right\}$
    \ElsIf{$2C = n$}
         \State $S_{C} \gets \left\{n-1, n\right\}$
    \Else
        \ForAll{$j \in \left\{2i+1, \dots, n \right\}$}
            \State $\sigma \gets \disunif{1}{C}$
            \State $S_\sigma \gets S_\sigma \cup \left\{j\right\}$
            \State $j \gets j + 1$
        \EndFor
    \EndIf \label{algline:evoflow-clu2}

    \ForAll{$\sigma \in \left\{1,\dots,C\right\}$} \label{algline:evoflow-in1}
        \State $L \gets \frac{C_{dens} |S_\sigma| (|S_\sigma| - 1)} {2}$
        \ForAll{$l \in \left\{1,\dots,L \right\}$}
            \State $(i,j) \gets \textrm{a random pair of facilities in } S_\sigma \textrm{ satisfying } f_{ij} = 0$
            \State $f_{ij} \gets \disunif{100}{C_{max}}$
            \State $f_{ji} \gets \disunif{100}{C_{max}}$
        \EndFor
        \While{the cluster of facilities represented by $S_\sigma$ has disconnected components} \label{algline:evoflow-conn1}
        \State $(i,j) \gets \textrm{a randomly selected pair of disconnected facilities in } S_\sigma$
            \State $f_{ij} \gets \disunif{100}{C_{max}}$
            \State $f_{ji} \gets \disunif{100}{C_{max}}$
        \EndWhile \label{algline:evoflow-in2}
    \EndFor \label{algline:evoflow-conn2}
    
    \State $L \gets N_{freq} \times n^2$ \label{algline:evoflow-out1}
    \ForAll{$l \in \left\{1,\dots,L \right\}$}
        \State $(i,j) \gets \textrm{a pair of facilities selected at random, without regard to cluster}$
        \State $f_{ij} \gets f_{ij} + \disunif{1}{N_{max}}$
    \EndFor \label{algline:evoflow-out2}
    
    \State {\bf return} $F$
    
    \end{algorithmic}
\end{algorithm}


\end{document}